\documentclass[12pt]{article}
\usepackage{amsfonts,amssymb,amscd}
\newtheorem{teo}{Theorem}[section]

\newtheorem{prop}[teo]{Proposition}

\newtheorem{lema}[teo]{Lemma}

\newtheorem{obs}[teo]{Remark}

\newtheorem{defnc}[teo]{Definition}

\newtheorem{coro}[teo]{Corollary}

\newcommand{\C}{{\mathbb C}}

\newcommand{\R}{{\mathbb R}}

\newcommand{\Q}{{\mathbb Q}}

\newcommand{\Z}{{\mathbb Z}}

\newcommand{\N}{{\mathbb N}}

\newcommand{\qed}{\hfill \fbox{}}

\newcommand{\fol}{{\mathcal F}}

\newcommand{\calk}{{\mathcal K}}

\newcommand{\wdc}{{\widetilde{\mathbb C}^2}}

\newcommand{\tilf}{{\widetilde{\mathcal F}}}

\newcommand{\wto}{{\widetilde{\omega}}}

\newcommand{\fracx}{{\partial /\partial x}}

\newcommand{\fracy}{{\partial /\partial y}}

\newcommand{\calh}{{\mathcal H}}

\newcommand{\Vv}{{\mathcal V}}

\newcommand{\Ww}{{\mathcal W}}

\newcommand{\supT}{{\rm supp}\, (T)}

\begin{document}

\title{On closed currents invariant by holomorphic foliations, {\sc I}}
\author{Julio C. Rebelo}
\date{}
\maketitle

\section{Introduction}

There are several motivations to pay attention to singular holomorphic foliations defined on complex surfaces
and admitting an invariant closed positive current. First the study of these foliations has a
natural ergodic theoretic interest. Additional motivation comes from observing that
this class of foliations captures and unifies special types of foliations such as those
having a compact leaf, foliations having a Zariski-dense leaf isomorphic to a quotient of $\C$ and Hilbert
modular foliations. It follows that the understanding of these foliations has consequences in a variety of
domains.

The purpose of this paper is to set up a dynamical method to study this type of foliations. The paper
is essentially elementary in that only the basic properties of this method are considered in detail.
Nonetheless, it will soon be apparent that the problem splits naturally into two cases corresponding
to the two sources of known examples of foliations admitting invariant closed currents, namely foliations having a compact leaf and
foliations carrying a transverse riemannian structure. Once this splitting will have been established in our context, we shall
pursue the case leading to the existence of a compact leaf. The other case, for which more subtle dynamical
arguments will be required, is going to be treated in the future.

Let us give a brief description of the point of view adopted here. The starting point is a construction appearing in
\cite{blm} which is revisited in Section~2.3. Essentially the
construction consists of observing that, inside the leaves of $\fol$,
there are real one-dimensional trajectories along to the which
the holonomy of $\fol$ ``tends to be contractive'' so long they
stay away from the singularities of $\fol$. These trajectories
can be viewed as defining a singular real one-dimensional foliation
denoted by $\calh$ (cf. Section~2.3 for details). Similar ideas already appeared in the context of foliated characteristic classes, see
\cite{ghys-bourbaki} and its references. The natural idea is then to ``follow'' these trajectories and to
understand how an invariant closed positive current can be compatible with the
resulting contraction of the holonomy of $\fol$. However, when this idea is further elaborated,
the above mentioned dichotomy manifests itself in the existence or absence of trajectories of
$\calh$ contained in the support of the current in question and possessing ``infinite length''
(cf. Section~5).
If trajectories of infinite length do exist, then they must yield
some ``definite amount of contraction'' in the holonomy pseudo-group of $\fol$. This clearly poses a
serious obstruction to the existence of the mentioned currents. The tension between
the existence of the invariant current and the existence of ``contractions'' in the holonomy pseudogroup is likely to
imply that the invariant current has a trivial nature: it is concentrated on a compact
leaf of $\fol$. On the other hand,
it may well happen that all the trajectories of $\calh$, or at least those
contained in the support of the current, are of ``finite length''. In vague terms,
this means that every trajectory has ``extremities'' or ``ends'' so that they cannot be followed ``for an arbitrarily large
period of time''. This prevents us from ensuring
the existence of ``contractions'' in the holonomy pseudo-group of $\fol$. Simple examples where this phenomenon is observed
arise in the context of transversely riemannian foliations whose discussion will be postponed to a subsequent work.

The contents of this article is then somehow bipartite. It begins with the general
definition of the foliation $\calh$ consisting of the mentioned real trajectories and it continues with the analysis of their
behavior near the singularities of $\fol$ that may exist in the support of a (invariant closed positive) current. Then
we turn to the global geometry of these trajectories. After some general reductions and definitions,
we shall finally be confronted with the basic dichotomy mentioned above: either
there are trajectories inside the support of the current having ``infinite length'' or
all these trajectories have finite length.

The remainder of the paper will then be devoted to the investigation of the first possibility, i.e. of the case in which
there are trajectories of infinite length contained in the support of our current. The corresponding results are ultimately summarized by Theorem~\ref{fim1}
(see also Theorem~A below).
The statement of this theorem is as follows: unless all the $\calh$-trajectories contained
in the support of the current have a uniformly bounded length, then this support contains
a compact leaf of $\fol$. Although the terminology is at this point imprecise, we can made
explicit the contents of this theorem as follows:
\begin{itemize}

\item at least as far as compact leaves are concerned, we only need to study those
foliations for which there is a compact invariant set where all the $\calh$-trajectories
have finite length.

\item if, by some reason, we can guarantee the existence of $\calh$-trajectories
of infinite length in the support of an invariant closed current $T$, then this support
contains a compact leaf.
\end{itemize}

As to the second item, some non-trivial examples
in which it is possible to prove {\it a priori}\, that no trajectory of $\calh$ on $M$ is of finite length will be supplied at the very
end of this article. These examples include some foliations on elliptic $K3$ surfaces having singularities that either are hyperbolic
or belong to the domain of Siegel.
It will be seen that for these examples, the cohomology class of an invariant closed current $T$ giving no mass to individual leaves must
have trivial self-intersection (and thus it
has the same cohomology of an elliptic fiber). Nonetheless it is not {\it a priori}\, clear that the foliation in question needs to have any compact leaf
at all. Another interesting class of examples includes foliations in the projective plane having singularities that may contribute non-trivially to the Lelong numbers
of $T$, so as to allow the foliated current $T$ to have strictly positive self-intersection. For example, consider a foliation $\fol$ on $\C P(2)$ having a radial
singularity $p$ and such that their remaining singularities belong to the Siegel domain or are hyperbolic. By a radial singularity, it is meant that $\fol$ is given in suitable coordinates about $p$ by the $1$-form $xdy -ydx$. If the degree of this foliation is at least~$3$, then we can arrange for the trajectories
of $\calh$ to have infinite length. Similarly if the degree of the foliation is at least~$4$ then the same result applies to foliations having up to~$2$
radial singularities (or more generally two singularities whose eigenvalues are $1, \lambda$ with $\lambda \in \R_+^{\ast}$).
Thus the results of this paper can be applied to these foliations to yield, in particular, the existence of algebraic
curves invariant for them. More details on the construction of examples can be found at the end of Section~7.

The study of the dynamics of these ``contractive trajectories'' is likely to have interest in different problems about holomorphic
foliations. For example, it may be useful to study the
dynamics of ``generic'' foliations such as in \cite{gabriel}. Similarly it should be mentioned that B. Deroin and V. Kleptsyn have
employed the foliated Brownian motion to study the transverse dynamics
of a conformal Riemann surface lamination, \cite{bertrand}. Roughly speaking
they show the evolution of ``most points'' under the Brownian
motion tends to give rise to a ``contractive holonomy''. In this sense
their work seems to be related to ours, i.e. the Brownian motion
evolution seems to be related to the foliation $\calh$. It would be
interesting to clarify possible relations between these approaches.
Along these lines, the first sections of this paper may also serve as an introduction to this circle ideas. In addition, we have included a preliminary section providing details
on some well-known facts that are usually not detailed in the literature. In particular the construction of \cite{blm} is reviewed.
We also explain in detail the role played by the condition of having an ambient surface that is algebraic as well as
the relation between Dirac masses for transversely invariant measures and compact leaves
for singular foliations. Hopefully this discussion will be useful for readers that are not experts in foliation theory.

Let us now state

\vspace{0.1cm}

\noindent {\bf Theorem A} (Main Theorem): {\sl Let $\fol$ be a singular holomorphic
foliation on an algebraic surface $M$.  Suppose that $\fol$ carries an invariant closed positive current $T$ and let
$\calk$ denote a (singular) minimal set for $\fol$ contained in the the support of $T$.
Suppose also that in $\calk$ there is one trajectory of $\calh$ having infinite length. Then
$\calk$ consists of an algebraic curve left invariant by $\fol$.}

\vspace{0.1cm}

Let us point out that the above theorem {\it does not}\, state that the current $T$ coincides
with the integration current over the mentioned algebraic curve (at least on a neighborhood of $\calk$).
A partial answer to this question is however supplied by
Theorem~B below. The reader can check Section~7 for the definition of Liouvillean integrability.

\vspace{0.1cm}

\noindent {\bf Theorem B} (Complement to Main Theorem): {\sl With the notations of Theorem~A, suppose that
$T$ is not locally given by integration over $\calk$. Then $\fol$ admits a Liouvillean
first integral on  a neighborhood of $\calk$.}

\vspace{0.1cm}

A by-product of the method developed here is the existence of a ``contractive'' element in its holonomy
pseudogroup provided that $\calh$ has a trajectory of infinite length on $M$. In principle, this contractive
element may be either a local hyperbolic diffeomorphism or a ``ramified'' super-attractive contraction
(cf. Section~$6$ for further details). Also a type of ``exponential super-contraction'' may appear in connection with saddle-nodes.
This is proven here under the additional assumption that the
singularities of $\fol$ yield no saddle-node singularities under the reduction procedure of Seidenberg (cf. Section~2, singularities verifying
this condition are sometimes called ``generalized curves''). The presence of saddle-node in the picture would actually not affect neither our
methods nor the validity of the conclusion. Nonetheless the elimination of this ``superfluous'' assumption would require us  to discuss the behavior of the
mentioned trajectories around a saddle-node and, in turn, this would lead us to a long detour in our way to the goals of the present paper. Yet, the study
of these trajectories near saddle-node singularities is an essential part of the analysis involved in the continuation of this work, so it seems
natural to let the general result about ``existence of contractions'' to be complemented there. We note however that, under the above assumption
concerning saddle-nodes singularities, the methods used here ensure the existence of contraction given either by a local hyperbolic diffeomorphism
or by a ``ramified'' super-attractive map, cf. Section~6.

To close this Introduction, let us briefly outline the contents
of this article. Section~2 contains background material on the subject.
In Section~2.1 we recall Seidenberg's reduction procedure for
the singularities of holomorphic foliations. Section~2.2 contains
precise definitions and a few basic facts regarding closed currents
invariant by a foliation including its relation with the notion of
transverse invariant measure. Finally, in Section~2.3,
the main ideas of \cite{blm} are presented. In particular
the foliation $\calh$ (associated to a given holomorphic foliation $\fol$)
is defined.

Section~3 is devoted to a detailed study of the behavior of $\calh$ on
a neighborhood of a singularity of $\fol$ belonging to the Siegel
domain. This includes a discussion of the ``Dulac transform'' defined
by means of $\calh$. Building on the material presented in Section~3, we develop in Section~4
a detailed analysis of the singularities of $\fol$ lying in the support
of $T$. The structure of these singularities will play an important role in the proof of Theorem~\ref{fim1} which is a slightly more
general version of Theorem~A.

Section~5 begins with a few global definitions, in particular the precise
definition of ``trajectory of $\calh$'' and of its corresponding
length. This section is technically very simple and its main result
is summarized by Propostion~\ref{4prop1}. Finally in Sections~6 and~7, we shall use
the preceding material to deduce the proofs of the above stated theorems.

\vspace{0.2cm}

\noindent {\bf Ackowledgements}: It is my pleasant duty to thank Emmanuel Paul for explaining me
his work \cite{paul} which allowed me to prove Theorem~B. I am also grateful to
B. Deroin and to A. Glutsyuk for the interest they have shown in this work.

\section{Preliminaries}

\subsection{Generalities about foliations on algebraic surfaces}

Let $M$ be a smooth compact complex surface. A {\it singular holomorphic foliation} $\fol$ on $M$ consists of the following data:

\begin{enumerate}
\item An open covering $\{ U_i \}$ of $M$.

\item A holomorphic vector field $Y_i$ with isolated singularities defined on
$U_i$ for each $i$.

\item For each pair $i, j$ such that $U_i \cap U_j \neq \emptyset$ a function
$g_{ij} \in {\mathcal O}^{\ast} (U_i \cap U_j)$ such that $Y_i = g_{ij} Y_j$.
\end{enumerate}
In particular it follows that the singularities of $\fol$ correspond to those of the
$Y_i$'s and are therefore isolated. The foliation $\fol$ can alternatively be defined
through holomorphic forms $\omega_i$ (rather than vector fields $Y_i$) subjected
to the relations $\omega_i = f_{ij}  \omega_j$ with $f_{ij} \in {\mathcal O}^{\ast}
(U_i \cap U_j)$. The transition functions $g_{ij}$ (resp. $f_{ij}$) satisfy the
natural cocycle relations and hence give rise to a line bundle $T^{\ast} \fol$
(resp. $N \fol$) on $M$ which is called the cotangent bundle of $\fol$
(resp. normal bundle of $\fol$). Furthermore, by virtue of the relations $Y_i = g_{ij} Y_j$,
the foliation $\fol$ can be interpreted as a global holomorphic section of $T^{\ast} \fol
\otimes TM$ with discrete zero set and modulo multiplication by non-vanishing
holomorphic functions. Thus we obtain.

\begin{lema}
\label{probablynotneeded}
A singular holomorphic foliation on an algebraic surface $M$ is always given by a globally
defined meromorphic form $\omega$ (resp. meromorphic vector field $Y$). Besides we can
suppose without loss of generality that the meromorphic form $\omega$ is not closed.
\end{lema}

\noindent {\it Proof}. The fact that $M$ is algebraic guarantees that every
line bundle over $M$ admits non-trivial meromorphic sections. In turn, this shows
that $\fol$ is generated by a global {\it meromorphic} vector field (or differential
$1$-form) in the obvious sense. Finally, if $\fol$ is given by a closed form $\omega$, then we just
need to replace $\omega$ by $f\omega$ for a generic meromorphic function $f$. In fact,
we have
$$
d(f\omega ) = df \wedge \omega + f d\omega = df \wedge \omega  \neq 0
$$
for a ``generic'' $f$. This proves the lemma.\qed

\begin{obs}
{\rm {\bf General Setting}: Throughout this work $M$ is supposed to be a complex surface equipped with a holomorphic foliation
$\fol$ which is given by a non-closed meromorphic form $\omega$. The preceding lemma shows that this
is always the case for $M$ projective algebraic. In the latter case, it would also be possible to choose a $1$-form $\omega$ satisfying
further ``generic conditions''. Suitable generic properties would simplify some parts of our discussion but we decided  not to use this.
The main reasons for our choice, besides having a slightly more general result, lies in the fact that some ``generic properties of $\omega$''
does not allow us, for example, to avoid a non-trivial intersection of $(\omega)_0$ and $(\omega)_{\infty}$ at a regular point of $\fol$. To
eliminate these intersection points a natural idea is to blow them up what, in turn, would bring us back to a situation where the corresponding transform of
$\omega$ is no longer ``generic''. Thus this transform would need to be replaced by a generic $1$-form on the blown-up surface and the final
construction would (even if successful) rely on constructions and arguments of algebraic geometry that appear to me as less elementary than the
approach chosen here. Indeed I also think that the treatment of all the ``degenerate situations'' that may arise for an arbitrary $\omega$ ends up making
the argument more ``concrete''.}
\end{obs}

The fact that $\fol$ is generated by global meromorphic differential
forms will be exploited in Paragraph~2.3. For the time being, we are going
to focus on local aspects such as the structure of the singularities of $\fol$. We can
then suppose that $\fol$ is given on a neighborhood of $(0,0) \in \C^2$ by a holomorphic
$1$-form $\eta = P dy + Q dx$ having an isolated singularity at the origin. Sometimes
it is also useful to think of $\fol$ as being given by the vector field $Y = P \fracx - Q
\fracy$. The
{\it order} of $\fol$ at $(0,0) \in \C^2$ is by definition the order of the first
non-zero jet of $\eta$ at $(0,0) \in \C^2$. This notion is well-defined since $\eta$ (or $Y$)
has isolated singularities. Similarly we define the eigenvalues of $\fol$ at
$(0,0)$ as the eigenvalues of the linear part of $Y$ at $(0,0)$. These eigenvalues
are therefore defined up to a multiplicative constant so that only their quotient
has an intrinsic meaning.

Let $\lambda_1 ,\lambda_2$ be the eigenvalues of $\fol$ at $(0,0)$. We say that
$(0,0)$ is a {\it hyperbolic} singularity if $\lambda_1 \lambda_2 \neq 0$ and
$\lambda_1 /\lambda_2 \in \C \setminus \R$.  If $\lambda_1 \lambda_2 \neq 0$
but $\lambda_1 /\lambda_2 \in \R_-$, then we say that $(0,0)$ is in the
{\it Siegel domain}. The singularity is said to be a {\it saddle-node} if $\lambda_1
\neq 0$ and $\lambda_2 =0$.

Singularities whose eigenvalues satisfy $\lambda_1 \lambda_2 \neq 0$ and
$\lambda_1 /\lambda_2 \in \R_+$ need some specific attention. Let us begin by saying that
a singularity (whose both eigenvalues are possibly zero) is {\it dicritical} if it
admits infinitely many {\it separatrices} i.e. anaytic curves passing through the
singularity and invariant under the foliation. Now consider a foliation whose eigenvalues
$\lambda_1 ,\lambda_2$ at the origin verify $\lambda_1 \lambda_2 \neq 0$
and $\lambda_1 =n\lambda_2$ for $n \in \N$ (or $\lambda_2 =n \lambda_1$).
This foliation is then conjugate to the foliation given either by
the $1$-form
$$
nx \, dy \; - \; y \, dx
$$
or by the $1$-form
$$
(nx + y^n) \, dy \; - \; y \, dx \, .
$$
In the first case $(0,0)$ is dicritical. In the second case it is
said to be a Poincar\'e-Dulac singularity (cf. \cite{arnold}). Next, if
$\lambda_1 \lambda_2 \neq 0$,  $\lambda_1 /\lambda_2 \in \R_+$ but
$\lambda_1 /\lambda_2$ is not an integer nor the inverse of an integer, then
Poincar\'e Theorem asserts that $\fol$ is linearizable (cf. \cite{arnold}). In other
words,
$\fol$ is conjugate to the foliation given by
$$
\lambda_1 x \, dy \; - \; \lambda_2 y \, dx \;  .
$$
The preceding implies that a singularity with eigenvalues $\lambda_1\lambda_2 \neq 0$ is dicritical if and only if
it is not a Poincar\'e-Dulac singularity and $\lambda_1 /\lambda_2 \in \Q_+$ and $(0,0)$. When
$\lambda_1 /\lambda_2 \in \R_+ \setminus \Q$ the resulting singularity is going to be called an
{\it irrational focus}.

Next we need to recall Seidenberg's reduction of singularities theorem \cite{seiden}.
Let $\pi : \wdc \rightarrow \C^2$ denote the blow-up of $\C^2$ at the origin. If $\fol$
is defined on a neighborhood $U$ of $(0,0) \in \C^2$, then $\pi^{\ast} \fol$ naturally
defines a holomorphic foliation on $\pi^{-1} (U)$. The foliation $\pi^{\ast} \fol
=\tilf_1$
is called the blow-up of $\fol$. Clearly this construction can be iterated: if $p$ is
a singularity of $\tilf_1$, then $\tilf_1$ can be blown up at $p$ to provde a new
foliation defined on an appropriate open surface. Seidenberg's theorem \cite{seiden}
then claims the existence of a finite sequence of blow-ups
$$
\fol=\fol_0 \stackrel{\pi_1}{\longleftarrow} \tilf_1 \stackrel{\pi_2}{\longleftarrow}
\cdots \stackrel{\pi_n}{\longleftarrow} \tilf_n
$$
such that the following holds:
\begin{itemize}
\item The irreducible components of the (total) exceptional divisor $(\pi_1 \circ \cdots
\circ \pi_n)^{-1} (0)$ are smooth rational curves $D_1 ,\ldots ,D_n$ of strictly
negative self-intersection.

\item The singularities of $\tilf_n$ are {\it reduced}\, i.e. they are of one of the
following
types: hyperbolic, in the Siegel domain, saddle-node or an irrational focus.

\end{itemize}

\noindent It should be noted that the exceptional divisor
$(\pi_1 \circ \cdots \circ \pi_n)^{-1} (0)$ need not be invariant by $\tilf_n$. In fact,
it
may contain irreducible components invariant under $\tilf_n$ along with irreducible
components that are not invariant by $\tilf_n$. If $D_i$ is an irreducible component
that is not invariant by $\tilf_n$, then the projection of the regular leaves of
$\tilf_n$
transverse to $D_i$ produces infinitely many separatrices for the initial foliation
$\fol =\fol_0$. In other words, $\fol$ is dicritical. Conversely if $\fol$ is dicritical
then, in the above situation, there must exist at least one irreducible component of
$(\pi_1 \circ \cdots \circ \pi_n)^{-1} (0)$ which is not invariant by $\tilf_n$. The
components
of $(\pi_1 \circ \cdots \circ \pi_n)^{-1} (0)$ that are not invariant by $\tilf_n$ are
also said to be {\it dicritical}.

\subsection{Closed currents and transverse invariant measures}

Let us now recall some standard definitions and results concerning closed invariant currents and
transverse invariant measures for singular foliations.

Given $M$ as before, we denote by $D^p (M)$ the Fr\'echet space of
$C^{\infty}$-differential
forms on $M$ of degree $p$. The space of currents of {\it dimension}\, $p$ is,
by definition, the topological dual $D'_p (M)$ of $D^p (M)$. The space of currents
possesses a natural differential ``$d$'' (as well as operators $\partial, \; \overline{\partial}$)
obtained by duality from the usual operators acting on differential forms.
In particular, we can talk about closed/exact currents.

Suppose now that $M$ is endowed with a (singular) foliation $\fol$. Consider
a current $T$ of dimension~$2$ and denote by $\supT
\subseteq M$ its support. The current $T$ is said to be {\it invariant}\, by $\fol$ if $T (\beta) =0$
for every $2$-form $\beta$ vanishing on $\fol$. An invariant current is sometimes also called a {\it foliated current}\,
or a {\it current directed by $\fol$}.
Because $M$ is a complex surface and $\fol$ is a holomorphic foliation, every current $T$ as above
is of type~$(1,1)$. In fact, on a neighborhood of a regular point, we can
choose coordinates $(x,y)$ in which $\fol$ is given by $dx =0$. Hence a $2$-form
$\beta$ vanishing on $\fol$ must be given by $\alpha_x \wedge dx +
\alpha_{\overline{x}} d\overline{x}$ where $\alpha_x, \, \alpha_{\overline{x}}$
are appropriate $1$-forms. Now the invariance of $T$ under $\fol$ becomes
encoded in the normal form
$$
T = T(x,y) dx \wedge d\overline{x}
$$
where $T(x,y)$ is naturally identified with a distribution. This shows that
$T$ is a $(1,1)$-current as claimed. For these currents, the notion of being {\it positive}\, becomes
especially transparent: $T$ is said to be positive if the local coefficient $T (x,y)$ is identified with a positive measure.
An equivalent condition consists of saying that for every smooth $(1,0)$-form $\alpha$ the wedge product
$T \wedge i\alpha \wedge \overline{\alpha}$ is a positive measure on $M$.

The most basic example of a foliation admitting a closed positive current
is provided by a compact leaf of a foliation $\fol$. More precisely,
let $C$ be a compact curve (smooth to simplify) which is invariant
by $\fol$. Consider then the current of integration over $C$, namely
the current $T$ given by
$$
T (\beta) = \int_C \beta \, .
$$
The fact that $C$ happens to be invariant by $\fol$ implies that $T$
is also invariant by $\fol$ in the sense mentioned above. Besides Stokes
formula shows that $T$ is, indeed, a closed current. Since $T$ is clearly
positive, we conclude that $T$ is a closed positive current invariant
by $\fol$.

Next we shall consider a more geometric point of view to study closed
currents invariant by a foliation $\fol$ as above. This point of view
relies on the notion of {\it transverse invariant measure}\, which is essentially due to
J. Plante. Since these transverse invariant measures are naturally defined
for {\it regular foliations}, we assume for the time being that ${\rm Sing} \, (\fol) =\emptyset$.

Since $M$ is compact and ${\rm Sing} \, (\fol) =\emptyset$, we can consider
a finite covering $\{ V_i \}$ of $M$ by foliated charts
$h_i : V_i \rightarrow \mathbb{D} \times \Sigma_i$ of $\fol$,
where $\mathbb{D}$ stands for the unit disc of
$\C$. The fact that the $h_i$'s define a foliated atlas implies
that the change of coordinates $h_j \circ h_i^{-1} (x,y)$ has the
special form $(f_{ij} (x,y), \gamma_{ij} (y))$.

\begin{defnc}
\label{plante?}
With the above notations, a transverse invariant measure for $\fol$ consists of
a collection $\mu_i$ of (positive) finite measures over the transverse sections $\Sigma_i$
which are invariant by change of coordinates. In other words,
for every pair $i,j$ and every Borel set $B \subset \Sigma_i$, one has $\mu_i (B)
= \mu_j (\gamma_{ij} (B))$.
\end{defnc}

Transverse invariant measures naturally provide closed (positive) currents
invariant by the foliation $\fol$ in question by means of the following
construction. Let $\{ \phi_i \}$ be a partition of the unity subordinated to
the finite covering $\{ V_i \}$. Given a $2$-form $\beta$, we define a current
$T$ by setting
\begin{equation}
T (\beta) = \sum_i \int_{\Sigma_i} \left ( \int_{\rm Plaque} (\phi_i \beta)
\right) d\mu_i \; , \label{sofoi}
\end{equation}
where the ``Plaque'' is naturally identified with the unit disc $\mathbb{D}$
through the coordinates $h_i$'s. It is easy to check that $T$ is a closed
positive current. In fact, it is a continuous linear operator on the Fr\'echet space of smooth $2$-forms. The conditions
of being closed, positive and invariant by $\fol$ can immediately be checked.
Conversely a closed positive current invariant by $\fol$ can be ``desintegrated''
to yield a transverse invariant measure so that the two objects turn out to be
equivalent as originally pointed out by D. Sullivan \cite{denis}.

Let us now go back to the case of a holomorphic foliation $\fol$ with
singularities. We then consider the open surface $M \setminus {\rm Sing}
\, (\fol)$ along with a covering $\{ V_i \}$, $i\in \N$, by foliated charts
$h_i : V_i \rightarrow \mathbb{D} \times \Sigma_i$ for the restriction
of $\fol$ to $M \setminus {\rm Sing} \, (\fol)$. Here
the covering $\{ V_i \}$ need not be finite. Again, if we are given a closed positive current $T$ invariant by $\fol$,
the procedure of ``desintegration'' mentioned above can still be carried out word-by-word to yield a transverse invariant
measure for the restriction of $\fol$ to $M \setminus {\rm Sing} \, (\fol)$ as in Definition~\ref{plante?}. Besides, from this transverse
invariant measure we can recover the current $T$ by means of Formula~(\ref{sofoi}). Whereas the ``summation over $i$'' (the indices
of foliated coordinates) may now be infinite, the series is naturally uniformly convergent so that it does define a current (i.e. a continuous
operator) that actually coincides with $T$.

Here it might be worth making a minor comment concerning the passage from an ``abstract'' transverse invariant measure for $\fol$ to a closed
positive current invariant by $\fol$. This remark however will not be used anywhere in this work since we always start with a current already defined
on $M$. Consider a transverse invariant measure for $\fol$ on the open set $M \setminus {\rm Sing} \, (\fol)$ as in Definition~\ref{plante?} and the
corresponding operator on smooth $2$-forms given by~(\ref{sofoi}). Since the summation over~$i$ is possibly infinite, it is necessary to make sure
that the operator in question is well-defined and continuous. This clearly amounts to bound the mentioned integral on a neighborhood of the singular
points of $\fol$. Whereas I believe that this bound always exist in the context of holomorphic foliations on complex surfaces, this problem cannot
{\it a priori}\, be reduced to an application of some Riemann extension or Hartogs theorem. The difficulty here being that we do not know {\it a priori}\,
whether or not the corresponding integration of $2$-forms is well-defined on a punctured neighborhood of the singularity in question.
We shall not elaborate on this discussion since, as mentioned, it is not necessary for our purposes.

Let us finish this paragraph with a well-known lemma that will often be used in the course of this work.

\begin{lema}
\label{atomicmass}
Let $T$ be a closed positive current invariant by $\fol$.
Assume that a point $p \in M \setminus {\rm Sing}\, (\fol)$ has positive
mass with respect to the transverse invariant measure for $\fol$ induced by $T$. Then the leaf $L_p$ of $\fol$ through
this point is contained in a compact curve.
\end{lema}

\noindent {\it Proof}. To prove the statement, let $\overline{L}_p$ denote the closure of
$L_p$. We just need to show that the set $\overline{L}_p \setminus L_p$ formed by the (proper) accumulation points of $L_p$
is contained in the singular set of $\fol$. Indeed, since
${\rm Sing}\, (\fol)$ has codimension~$2$, it follows from the classical theorem of Remmert-Stein
that $\overline{L_p}$ is an analytic set.

To check the claim, suppose for a contradiction that $q$ is a regular point of $\fol$ belonging to
$\overline{L}_p \setminus L_p$. Consider a trivializing coordinate around $q$. Since $q \in \overline{L}_p \setminus L_p$,
there exists a sequence of points $\{p_i \} \subset L_p$ such that $p_i \rightarrow q$. Besides for $i\neq j$, $p_i, p_j$ belong
to different plaques of the mentioned foliated chart. Denoting by $\Sigma$ the corresponding local transversal, the measure
on $\Sigma$ associated to each of these plaques is a positive constant. It then follows from Equation~\ref{sofoi}
that the corresponding current has ``infinite mass'', i.e. the integrals in~(\ref{sofoi}) diverge for a suitable choice of $\beta$.
The resulting contradiction establishes the lemma.\qed

\subsection{Brief review of Bonatti-Langevin-Moussu}

In this paragraph we shall expand on the method developed in \cite{blm}
to producing hyperbolic holonomy for certain holomorphic foliations (cf. also \cite{ghys-bourbaki} and references therein). The study of the oriented foliation
$\calh$ consisting of trajectories yielding ``contractive holonomy'' is the central object of this section.

Consider a surface $M$ endowed with a holomorphic foliation $\fol$ as before.
Let $\omega$ be a global non-closed meromorphic $1$-form defining $\fol$ on $M$. The existence of this form is guaranteed if $M$
is projective (cf. Section~2.1).
Also denote by $(\omega )_0$ (resp. $(\omega)_{\infty}$) the divisor of zeros (resp. poles) of $\omega$. Note that, in most
applications, the sets $(\omega )_0, (\omega )_{\infty}$ are viewed
as ordinary algebraic curves rather than as divisors (i.e. no multiplicity is associated to their components).

Next let $\omega_1$ be the $1$-form defined by
\begin{equation}
d\omega = \omega \wedge \omega_1 \, . \label{GV}
\end{equation}
To obtain a $1$-form $\omega_1$ satisfying the equation above it suffices to find a meromorphic vector field on $M$ such that $\omega (X) =1$.
In fact, for this vector field $X$
we have $d\omega = \omega \wedge \mathcal{L}_X (\omega)$, where $\mathcal{L}_X$ stands for the Lie derivative.
Note also that two $1$-forms satisfying the mentioned equation must differ by a multiple of $\omega$. In particular it follows that the values of $\omega_1$ on
vectors tangent to $\fol$ are well-defined even though $\omega_1$ is not so. This ambiguity however can be avoided if $\omega_1$ is regarded
as a {\it foliated $1$-form}\, (as opposed to an ``ordinary'' $1$-form). By a foliated $1$-form, it is meant a $1$-form that is defined only for vectors tangent to (regular)
leaves of $\fol$. In other words, a foliated $1$-form is not a usual $1$-form on $M$ since at a generic point of $p$ this form is not defined for vectors
in $T_pM$ that are transverse to the leaf of $\fol$ through $p$. Still another way of thinking of a foliated $1$-form consists of saying that it is a meromorphic
section of the cotangent bundle of $\fol$, cf. Section~2.1. The preceding discussion can then be summarized by stating that Equation~(\ref{GV}) unequivocally
defines a meromorphic foliated $1$-form on $M$. This foliated $1$-form will systematically be denoted by $\omega_1$.

The foliated $1$-form $\omega_1$ can explicitly be computed.
If $(x,y)$ are local coordinates about a regular point $p$ of $\fol$ in which $\omega = F(x,y) dy$ then $\omega_1$ is given by
$$
-\frac{\partial F/ \partial x}{F} dx \, .
$$
Clearly the above definition is compatible with foliated changes of coordinates so that it gives rise to a global (meromorphic) foliated
$1$-form $\omega_1$ on $M$ or, equivalently, to a global meromorphic section of the cotangent bundle of $\fol$.
This formula also shows that the form $\omega_1$ is holomorphic on a neighborhood of $p$ unless $p$ belongs to the
union of $(\omega )_0$ and $(\omega )_{\infty}$. A more accurate statement concerning the holomorphic nature of $\omega_1$ is given below.

\begin{lema}
\label{newversionSection2.11}
Let $p \in M$ be a regular point of $\fol$. Suppose that all the irreducible components of $(\omega )_0 \cup (\omega )_{\infty}$
passing through $p$ are invariant by $\fol$. Then $\omega_1$ is holomorphic at $p$.
\end{lema}

\noindent {\it Proof}. Consider foliated coordinates $(x,y)$ about $p$ so that $\omega$ becomes $F(x,y) dy$.
We can assume that $p$ belongs to $(\omega )_0 \cup (\omega )_{\infty}$, otherwise $\omega_1$ is holomorphic as already seen.
Nonetheless the assumption that all components of $(\omega )_0 \cup (\omega )_{\infty}$ passing through $p$ are invariant
by $\fol$ implies that there can be only one component which coincides in the coordinates $(x,y)$ with the axis $\{ y=0\}$. In other
words, we have $\omega = F(x,y) dy = y^k f(x,y) dy$ where $k \neq 0$ and for some holomorphic function $f$ satisfying
$f(0,0) \neq 0$. Now a direct computation yields
$$
\omega_1 = -\frac{\partial F /\partial x}{F} dx = -\frac{\partial f / \partial x}{f} dx \, .
$$
The statement follows since $f(0,0) \neq 0$.\qed

Conversely we have:

\begin{lema}
\label{newversionSection2.22}
Let $C \subset M$ be an irreducible component of $(\omega )_0 \cup (\omega )_{\infty}$ that is not invariant by $\fol$. Then $\omega_1$
has poles of order~$1$ over $C$.
\end{lema}

\noindent {\it Proof}. Let $p \in C$ be a regular point for $\fol$ which does not belong to any irreducible component of
$(\omega )_0 \cup (\omega )_{\infty}$ different from $C$ itself. It suffices to show that the divisor of poles of $\omega_1$ locally
coincides with $C$ with multiplicity equal to~$1$. As before we can choose foliated coordinates $(x,y)$ about
$p$ where $\omega = F(x,y) dy = x^k f(x,y) dy$, $k \neq 0$, for some holomorphic function $f$ satisfying
$f(0,0) \neq 0$. Now
$$
\omega_1 = -\frac{\partial F/ \partial x}{F} dx = - \frac{k}{x} -\frac{\partial f / \partial x}{f} dx \, .
$$
The statement follows.\qed

Consider now a regular leaf $L \subset M$ of $\fol$ where $\omega_1$ does not vanish identically. The restriction of $\omega_1$
to $L$ is a meromorphic $1$-form on the Riemann surface $L$ (we shall often say that it is an abelian form on $L$).
Therefore it induces a pair of (real one-dimensional) oriented singular
foliations on $L$, namely the foliations given by
$\{ {\rm Im}\, (\omega_1 )=0 \}$ and $\{ {\rm Re}\, (\omega_1 )=0 \}$. These foliations
will respectively be denoted by $\calh$ and $\calh^{\perp}$ and they are mutually orthogonal
for the underlying conformal structure of $L$. The orientation of $\calh$ (resp.
$\calh^{\perp}$)
is determined by the increasing direction of ${\rm Re}\, (\omega_1)$ (resp.
${\rm Im}\, (\omega_1 )$). More generally, the conformal structure of $L$ also allows us to
define the oriented foliation
$\calh^{\theta}$ whose trajectories form an angle $\theta$ with those of $\calh$ (where $\theta$ belongs to $(-\pi/2, \pi/2)$).
Finally by letting the leaf $L$ vary, the foliations
$\calh, \calh^{\perp}$, or more generally $\calh^{\theta}$, can also be viewed as singular foliations defined on $M$. We shall return to this point
when discussing the singularities of $\calh, \calh^{\perp}$.

Now the discussion in Lemma~\ref{newversionSection2.22} yields the following lemma borrowed from \cite{blm}.

\begin{lema}
\label{blm1}
Let $p \in M$ be a regular point of $\fol$ and denote by $L$ the leaf through $p$. Let $C \subset M$ be an irreducible component of $(\omega )_0 \cup (\omega )_{\infty}$
that is not invariant by $\fol$. Then we have:

\noindent 1.  $p \in C$ but $p$ does not belong to $(\omega)_0$.
Then $p$ is a source for $\calh$. Precisely there is a (complex one-dimensional) local
coordinate $\textsc{X}$ along $L$ where the restriction of $\omega_1$ to
$L$ becomes $\omega_1 = md\textsc{X}/\textsc{X}$, $m \in \N^{\ast}$. In particular
the leaves of $\calh$ are radial lines emanated from $p \in L$ (identified to $0 \in
\C$).

\noindent 2. $p \in C$ but $p$ does not belong to $(\omega)_{\infty}$.
Then $p$ is a sink for $\calh$. Precisely there is a (complex one-dimensional) local
coordinate $\textsc{X}$ along $L$ where the restriction of $\omega_1$ to
$L$ becomes $\omega_1 = - md\textsc{X}/\textsc{X}$, $m \in \N^{\ast}$. In particular
the leaves of $\calh$ are radial lines converging to $p \in L$ (identified to $0 \in
\C$).
\end{lema}

\noindent {\it Proof}. Consider the first case. Since $p$ is regular, we have
$\omega = d\textsc{Y}/f (\textsc{X}, \textsc{Y})$ where $f (0,0) =0$ for suitable
coordinates $\textsc{X}, \textsc{Y}$. The fact that $f(0,0) =0$ follows from the
assumption $p \in (\omega)_{\infty}$ and $p \not\in (\omega)_0$. Modulo
performing a further change of coordinates, we can assume without loss
of generality that $f (\textsc{X}, 0) = \textsc{X}^m$ for some $m \in \N^{\ast}$.
Now the equation $d\omega = \omega \wedge \omega_1$ yields the
desired form for $\omega_1$. The second case can analogously be treated.\qed

Naturally an analogous discussion applies to the foliations $\calh$ (with $\theta \in (-\pi/2, \pi/2)$).
Thus we already know that components of $(\omega )_0 \cup (\omega )_{\infty}$ that are not invariant by $\fol$ give rise to singularities of
$\calh, \,  \calh^{\perp}$ at regular points of $\fol$. Next the foliated form $\omega_1$ also have a divisor of zeros (resp. poles)
denoted by $(\omega_1)_0$ (resp. $(\omega_1)_{\infty}$). It follows from the combination of Lemma~\ref{newversionSection2.11}
and Lemma~\ref{newversionSection2.22} that $(\omega_1)_{\infty} \subset (\omega)_0 \cup (\omega)_{\infty}$. Besides no irreducible
component of $(\omega_1)_{\infty}$ can be invariant by $\fol$. Let us now consider a component $C$ of $(\omega_1)_0$ that is not invariant
by $\fol$. The following lemma is also borrowed from \cite{blm}.

\begin{lema}
\label{blm2}
Let $p$ be a regular point of $\fol$ which does not belong to $(\omega)_0 \cup
(\omega)_{\infty}$. Suppose that $p$ lies in a component $C$ of $(\omega_1)_0$ that is not invariant
by $\fol$ and denote by $L$
the leaf of $\fol$ containing $p$. Then the behavior of $\calh$ at $p$ is that of
a saddle with $2m$ separatrices. Precisely, in suitable coordinates $\textsc{X}$ along
$L$, the restriction of $\omega_1$ to $L$ becomes $\omega_1 = m \textsc{X}^{m-1}\,
d\textsc{X}$ for $m \geq 2$.
\end{lema}

\noindent {\it Proof}\,: Since $p$ is regular and $p \not\in (\omega)_0 \cup
(\omega)_{\infty}$,
there are local coordinates $\textsc{X}, \textsc{Y}$ around $p$ in which $\omega =
f (\textsc{X}, \textsc{Y}) d\textsc{Y}$ with $f$ holomorphic. Suppose first that
$f(\textsc{X}, 0)$ is not trivial. Then the restriction of $\omega_1$ to $\{ \textsc{Y}
=0\}$
is given by $-(\partial f /\partial \textsc{X}) d\textsc{X}/f$ where the functions are
evaluated at $(\textsc{X} ,0)$. The result then follows. On the other hand, if $\omega$
vanishes
identically on $\{ \textsc{Y} =0\}$ (or has poles over this leaf) then
$\omega =y^k f \textsc{Y}$ with $f$ as before. This still gives
$\omega_1 = -(\partial f /\partial \textsc{X}) d\textsc{X}/f$ so that the statement
follows.\qed

\begin{obs}
\label{2obs2}
{\rm Note that $m \textsc{X}^{m-1} d\textsc{X}$ is nothing but the lift of the regular
form $d \textsc{X}$ through the ramified covering $\textsc{X} \mapsto \textsc{X}^m$.
In particular  $m \textsc{X}^{m-1} d\textsc{X}$ has $2m$ separatrices (namely
the lifts of the separatrices $\R_+$ and $\R_-$ of $d\textsc{X}$) with alternate
orientation.}
\end{obs}

Let us split the divisor of zeros $(\omega)_0$ of $\omega$ into two divisors $(\omega)_0^{\fol}$ and
$(\omega)_0^{\perp \fol}$ as follows: an irreducible component $C$ of $(\omega)_0$ belongs to $(\omega)_0^{\fol}$ if and only if
it is invariant by $\fol$. Otherwise it belongs to $(\omega)_0^{\perp \fol}$ (the multiplicity of each component remaining unchanged).
Similarly we define the split of $(\omega)_{\infty}$ into $(\omega)_{\infty}^{\fol}$ and $(\omega)_{\infty}^{\perp \fol}$.

Let us now summarize the information so far obtained about possible singular points of $\calh, \, \calh^{\perp}$ (and of $\calh^{\theta}$).
Singular points for these foliation belong to the list below.
\begin{enumerate}
\item Singular points of $\fol$ (to be detailed later).

\item Irreducible components of $(\omega)_0^{\perp \fol}$. These points are sink singularities for $\calh$.

\item Irreducible components of $(\omega)_{\infty}^{\perp \fol}$. These points are source singularities for $\calh$.

\vspace{0.1cm}

Note that the foliated $1$-form $\omega_1$ is holomorphic away
from ${\rm Sing}\, (\fol) \cup (\omega)_0^{\perp \fol} \cup (\omega)_{\infty}^{\perp \fol}$.
In particular the support of the pole divisor of $\omega_1$ are the union of the components of $(\omega)_0^{\perp \fol}$ and $(\omega)_{\infty}^{\perp \fol}$.
Splitting also the zero divisor $(\omega_1)_0$ of $\omega_1$ into divisors $(\omega_1)_0^{\fol}$ and $(\omega_1)_0^{\perp \fol}$, consisting respectively
of components that are invariant by $\fol$ and components that are not invariant by $\fol$, we identify the locus of the possible additional singularities of
$\calh$ (resp. $\calh^{\perp}, \, \calh^{\theta}$), namely:

\vspace{0.1cm}

\item Irreducible components of $(\omega_1)_{0}^{\perp \fol}$. These components yield saddle singularities for $\calh$ (resp. $\calh^{\perp}, \,
\calh^{\theta}$).

\item Irreducible components of $(\omega_1)_{0}^{\fol}$. Over these compact leaves of $\fol$ the foliations $\calh, \, \calh^{\perp}$ and $\calh^{\theta}$
are not defined.

\end{enumerate}

We are now able to give the geometric meaning of the foliated form $\omega_1$. Consider a path $c$ contained in a leaf
$L$ of $\fol$ along with local transverse sections $\Sigma_0, \, \Sigma_1$ respectively through $c(0), \, c(1)$. The parallel transport
over the leaves of $\fol$ gives rise to a local diffeomorphism ${\rm Hol}\, (c)$ from $\Sigma_0$ to $\Sigma_1$ taking $c(0)$ to
$c(1)$ called the holonomy map of $\fol$ over $c$. It is well-known that the
derivative of ${\rm Hol}\, (c)$ is not intrinsically defined unless $c$ is a loop. These derivatives
however can be considered for a fixed choice of parametrizations for the transverse sections $\Sigma_0, \, \Sigma_1$
(and for fixed parametrizations they can be considered whether or not $c$ is a loop).
To be more precise,
suppose that $\Sigma_0, \, \Sigma_1$ are parameterized by $\omega$, i.e. consider local coordinates $\varphi_i : \Sigma_i \rightarrow \C$, $i=0,1$, defined
by
$$
\varphi_0 (p) = \int_{c(0)}^p \omega \; \; \, {\rm and} \; \; \, \varphi_1 (q) = \int_{c(1)}^q \omega
$$
where the integrals are well-defined modulo choosing $\Sigma_0, \, \Sigma_1$ simply connected. In these coordinates, the holonomy map
${\rm Hol}\, (c)$ can be identified to a local diffeomorphism of $(\C, 0)$. This local diffeomorphism satisfies
\begin{equation}
({\rm Hol}\, (c))' (c(0)) = \exp \left ( -\int_c \omega_1 \right) \, . \label{PLemma}
\end{equation}
This formula is sometimes referred to as Poincar\'e Lemma. Several comments are needed here to fully explain its meaning. First let us
fix a (finite) covering of a compact part $K$ of $M \setminus {\rm Sing}\, (\fol)$ by foliated coordinates $\varphi_i : U_i \rightarrow \C^2$ where
each $U_i$ is equipped with a local transverse section $\Sigma_i$ parametrized by $\omega$ as above. Setting $\varphi_i (U_i) = D \times {\bf T}_i$
where $D$ stands for the unit disc and where ${\bf T}_i$ is identified with $\Sigma_i$ parameterized as indicated, the Poincar\'e Lemma
becomes applicable to every path $c \subset K$ contained in a leaf of $\fol$ (modulo an obvious decomposition of $c$ into paths contained
in the open sets $U_i$). Some further remarks are needed:
\begin{itemize}
\item The neighborhoods $U_i$ are chosen so that $(\omega)_0^{\perp \fol} \cap \varphi_i^{-1} (\partial D \times {\bf T}_i) = \emptyset$.
Similarly $(\omega)_{\infty}^{\perp \fol} \cap \varphi_i^{-1} (\partial D \times {\bf T}_i) = \emptyset$ and
$(\omega_1)_0^{\perp \fol} \cap \varphi_i^{-1} (\partial D \times {\bf T}_i) = \emptyset$.

\item If $C$ is a component of $(\omega)_0^{\fol}$ (in particular $C$ is invariant by $\fol$), then the parametrization of $\Sigma_i$ is actually
ramified at the ``origin'' (it is a local ramified covering rather than a local diffeomorphism). An analogous conclusion (on a neighborhood of
infinity) applies to components of $(\omega)_{\infty}^{\fol}$.

\item If $K' \subset K$ is a compact part of $M \setminus ({\rm Sing}\, (\fol) \cup (\omega)_0^{\fol} \cup (\omega)_{\infty}^{\fol})$ then the parametrization
of $\Sigma_i$ restricted to $K'$ is ``equivalent'' to the parametrization induced by an auxiliary Hermitian metric on $M$ (here it is to be noted
that the first item above ensures that no $\Sigma_i$ intersects $(\omega)_0^{\perp \fol}$, $(\omega)_{\infty}^{\perp \fol}$).
In fact, every Hermitian metric on $M$
induces parametrization that are pairwise ``equivalent'' in the sense that ``corresponding lengths'' are mutually controlled, from below and by above,
by multiplicative constants.

\item If $c \subset K'$ is a path contained in a trajectory $l$ of $\calh$, $l \subset L$ where $L$ is a leaf of $\fol$, then the holonomy ${\rm Hol}\, (c)$ is such that
$({\rm Hol}\, (c))' (c(0))$ is strictly smaller than~$1$ (with respect to above fixed foliated coordinates). Indeed, by construction, the integral of $\omega_1$
over $c$ increases monotonically with the length of $c$, cf. Formula~(\ref{PLemma}).

\end{itemize}

Throughout the paper, we shall assume that $\fol$ is not a pencil, that is, not all the leaves of $\fol$ are compact. According to Jouanolou \cite{joa} this actually means
that $\fol$ leaves only finitely many algebraic curves invariant. In particular the support of $(\omega)_0^{\fol} \cup (\omega)_{\infty}^{\fol} \cup
(\omega_1)_0^{\fol}$ consists of finitely many algebraic curves (if not empty).

We are now ready to explain the fundamental observation of \cite{blm}. Let $K'$ be as above and consider a path $c \subset L \cap K'$ parametrizing
a trajectory of $\calh$ (i.e. $\omega_1 (c(t)). c'(t)$ is always a nonnegative real number). In particular the holonomy map ${\rm Hol}\, (c)$ (measured with
respect to the identifications fixed above) is such that $({\rm Hol}\, (c))' (c(0))$ decays exponentially with the {\it length of $c$}. The notion of
{\it length of $c$} can be identified with the length measured in $D$ for each coordinate $\varphi_i$. Alternatively this length can be measured with
respect to the fixed auxiliary Hermitian metric on $M$ (the two notions of lengths being equivalent up to multiplicative constants, i.e. the metrics induced
on $L$ are quasi-isometric). Also, because $c$ is contained in $K'$, the notions of distance in the transverse sections $\Sigma_i$ induced by the parametrization
through $\omega$ and by the mentioned Hermitian metric are mutually controlled by multiplicative constants.

Let then $c$ be defined on the interval $[0,t_0] \subset \R$. The corresponding holonomy map ${\rm Hol}\, (c)$
is then defined on a small disc $D_0 (r) \subset \Sigma_{i_0}$ (for some $i_0$). Naturally ${\rm Hol}\, (c)$
maps $D_0 (r)$ diffeomorphically onto a neighborhood of $c (t_0) \in \Sigma_{i_1}$ (for some $i_1$).
It is observed in \cite{blm} that the contractive character of the holonomy
along the oriented leaves of $\calh$ allows one to have a uniform bound on
the radius of $D_0 (r)$ regardless of the point $c (t_0)$ and of the length of
the $c$.  Denoting by $\calh_{\vert K'}$ the restriction of $\calh$ to $K'$ one has:

\begin{teo}
\label{blm}
{\rm ({\bf [B-L-M]})} There is a uniform $r > 0$ with the following properties:

\noindent 1. Let $l_p$ be an oriented trajectory of $\calh_{\vert K'}$ passing through
$p \in K'$. If $c: [0 , t_0] \rightarrow l_p \subset L_p \subset K$ is a parametrization
of (a segment of) $l_p$ ($p = c(0)$), then the corresponding holonomy map
${\rm Hol}\, (c)$ is defined on $D_p (r) \subseteq \Sigma_{i_0}$ (for some $i_0$). Besides ${\rm Hol}\, (c)$
maps $D_p (r)$ diffeomorphically onto its image in $\Sigma_{i_1}$ (for some $i_1$).

\noindent 2. Assume, in addition, that the distance of $l_p$ to the divisor
$(\omega_1)_0^{\perp \fol}$ is bounded from below by a positive constant $\delta$. Then
there are uniform constant $C >0$, $k >0$ ($k$ depending solely on $\delta$)
such that
\begin{equation}
\vert ({\rm Hol}\, (c)) (q) \vert \leq  C  \exp \, (-k \, {\rm length}\, (c) /2) \; , \label{Contraction1}
\end{equation}
for every $q \in D_p (r)$ and where ${\rm length}\, (c)$ stands for the length
of the path $c$.\qed
\end{teo}

In item~2 above, it is to be noted that the asymptotic exponential decay of the diameter of the set $({\rm Hol}\, (c)) (D_p(r))$ has an
intrinsic meaning since the length of $c$ (as well as the the notion of distance in the transverse sections $\Sigma_i$ restricted to $K'$)
vary in a way controlled by multiplicative constants as pointed out above. In particular if these metrics are changed, Formula~(\ref{Contraction1})
remains valid modulo changing the values of the constants $C, k$.

\begin{obs}
\label{2obs3}
{\rm The reader will check that the same statement remains true for the foliations
$\calh^{\theta}$ for a fixed $\theta$ in $(-\pi/2, \pi/2)$. All these statements will be revisited and sharpened later in this paper.}
\end{obs}

\section{The structure of $\calh$ around a singularity in the Siegel domain}

The local structure of $\calh$ around a regular point of $\fol$ was described in the preceding section. The next step is to discuss
the analogous problem on a neighborhood of a singularity $p$ of $\fol$ which belongs to the Siegel domain.
About this singularity there are coordinates $(u,v)$ ($p\simeq (0,0)$) in which $\omega$ becomes
\begin{equation}
\omega = h(u,v) [\lambda_1 u (1+r^1(u,v)) \, dv \; + \; \lambda_2 v(1+r^2(u,v)) \, du]
\label{siegel1}
\end{equation}
where $h$ is meromorphic and $r^1 ,r^2$ are holomorphic functions verifying
$r^1 (0,0) =r^2(0,0) =0$. Finally one also has $\lambda_1 \lambda_2 \neq 0$ and
$\lambda_1/\lambda_2 \in \R_+$ (as to the sign conventions, note that we are now using
differential forms, rather than vector fields, to represent a singularity in the Siegel
domain). In particular $\fol$ possesses exactly $2$ separatrices at $p \simeq (0,0)$
namely, those given by $\{ u=0\}$ and $\{ v=0\}$. Throughout this section we work under the following extra-assumption:

\noindent {\it Local invariance condition}: One has $h(u,v) =u^a v^b$ for some $a,b \in \Z$.

The contents of the local invariance condition is that, on a small neighborhood of $p$, the
curves $(\omega)_0$ and $(\omega)_{\infty}$ are invariant by $\fol$. As it will be
shown later, this assumption does not affect the generality of our
arguments since it can always be obtained by performing finitely many
blow-ups.

Because of the local invariance condition, the $1$-form $\omega$ can be written in the
coordinates $u,v$ as
\begin{equation}
\omega = u^a v^b [\lambda_1u (1+r^1(u,v)) dv + \lambda_2 v (1+r^2 (u,v)) du ] \, .
\label{siegel2}
\end{equation}
As seen in Section~2.3, the foliated form $\omega_1$ can be obtained by restriction to the leaves of $\fol$ of an actual
(locally defined) $1$-form $\Omega_1$ satisfying $d\omega = \omega \wedge \Omega_1$.
Setting $\Omega_1 = fdv + gdu$, it follows that
\begin{equation}
f\lambda_2 v (1+r^2) - g \lambda_1 u (1+r^1) = \lambda_1 (1+a)(1+r^1) + \lambda_1
ur^1_u - \lambda_2(1+b) (1+r^2) -\lambda_2 vr^2_v \label{siegel3}
\end{equation}
where $r^1_u$ (resp. $r^2_v$) stands for the partial derivative of $r^1$
(resp. $r^2$) with respect to $u$ (resp. $v$).

Let us first consider the behavior of $\calh$ on the separatrices $\{ u=0\}$ and $\{ v=0\}$.

\begin{lema}
\label{3lema1}
Suppose that $\lambda_1(1+a) - \lambda_2 (1+b) \neq 0$. Then the
behavior of $\calh$ over one of the separatrices is that of a sink (as in
Lemma~\ref{blm1}).
Besides, on the other separatrix, $\calh$ behaves like a source.
\end{lema}

\noindent {\it Proof}\,:  The restriction of $\Omega_1$ to $\{ v=0\}$ is the
Abelian form $g(u,0) du$. By letting $v=0$ in Equation~(\ref{siegel3}), we obtain
$$
-g(u,0)\lambda_1 u (1+r^1(u,0)) = \lambda_1 (1+a)(1+r^1(u,0)) + \lambda_1 ur^1_u (u,0)
-\lambda_2(1+b)(1+r^2(u,0)) \, .
$$
Hence
\begin{equation}
g(u,0) = - \frac{\lambda_1(1+a) -\lambda_2(1+b)}{\lambda_1 u} + \widetilde{s}_g (u)
\label{siegel4}
\end{equation}
where $\widetilde{s}_g (u)$ is holomorphic around $0 \in \C$. Similarly, on
$\{ u=0\}$, $\Omega_1$ becomes $f(0,v) dv$ and Equation~(\ref{siegel3}) yields
\begin{equation}
f(0,v) = \frac{\lambda_1(1+a) -\lambda_2(1+b)}{\lambda_2 v} + \widetilde{s}_f (v)
\label{siegel5}
\end{equation}
where $\widetilde{s}_f (v)$ is holomorphic around $0 \in \C$. Since $\lambda_1/\lambda_2
\in \R_+$, the statement follows from comparing Equations~(\ref{siegel4})
and~(\ref{siegel5}) and recalling that the restriction of $\Omega_1$ to the leaves of $\fol$ coincides with $\omega_1$.\qed

\begin{obs}
\label{restrictionholomorphic}
{\rm In the case where $\lambda_1(1+a) - \lambda_2 (1+b) = 0$ the calculation above shows that $\omega_1$ is holomorphic
over both separatrices of $\fol$ at $p$.}
\end{obs}

As a matter of fact we also need to control $\omega_1$ (or equivalently to understand the trajectories of $\calh$) on the leaves of $\fol$
distinct from the separatrices. To abridge notations, in the sequel $\fol$ is going to
be considered as a foliation defined on a neighborhood $U$ of $(0,0) \in \C^2$. The
corresponding arguments should also be understood modulo reducing this neighborhood.

Recall that $\fol$ is defined on $U$ by the $1$-form
$$
\eta = \lambda_1 u (1+r^1) dv + \lambda_2 v (1+r^2) du
$$
so that $\omega = u^a v^b \eta$. For the rest of this section we always suppose that
$\lambda_1(1+a) - \lambda_2 (1+b) \neq 0$.

\begin{lema}
\label{3lema2}
The form $\omega_1$ has no zeros on $U$. In fact, there is a positive constant
$C >0$ such that for every $p \in L \subset U$ and unit vector $\mathtt{v} \in T_pL$ one
has
$$
\Vert \omega_1 (p) . \mathtt{v} \Vert \geq C >0 \, .
$$
\end{lema}

\noindent {\it Proof}\,:  Given $\epsilon_1 ,\epsilon_2 >0$ sufficiently small, let us
denote by $\Sigma$ the local transverse section defined by
$$
\Sigma = \{ (u,v) \in \C^2 \;  ;  \; u =\epsilon_1 \; {\rm and} \; \vert v \vert <
\epsilon_2 \} \; .
$$
Denoting by $\Sigma_{\fol}$ the saturated of $\Sigma$ by $\fol$, it is proved
in \cite{mamo} (cf. also \cite{mat}, \cite{reis})
that $\Sigma_{\fol} \cup \{ u=0 \}$ contains a neighborhood
of $(0,0) \in \C^2$. Since we need a parametrization of the leaf $L$ in order to
estimate the restriction of $\omega_1$ to $L$, let us consider the set
$D_1^- = \{ (u,v) \in \C^2 \; ; \; \vert u \vert < \epsilon_1 \; {\rm and} \;
v \not\in \R_- \} $. We then define $W = \{ (u,v) \in \C^2 \; ; \;  u \in D_1^- \; {\rm
and} \;
\vert v \vert \leq \epsilon_2 \}$. Because $D_1^-$ is simply connected, the
restriction of $\fol$ to $W$ does not present the local holonomy associated to
the separatrix $\{ v=0 \}$. Thus fixed $y_0 \in \Sigma$, the leaf $L_0$ of $\fol$
restricted to $W$ through $y_0$ is the graph of a holomorphic
function $h$. Precisely the argument of \cite{mamo} shows the existence
of $h: D_{y_0} \subset D_1^- \rightarrow \C$ whose graph $\{ (u, h(u)) \}$, $u \in
D_{y_0}$, coincides with $L_0$.

Clearly to obtain estimates for the restriction of $\omega_1$ to the leaves of $\fol$,
it suffices to estimate $\omega_1$ over leaves $L_0$ as above since $\R_-$ can be
substituted by another semi-line in the definition of $D_1^-$. Now, when dealing with
$L_0$, we are allowed to use the parametrization $u \mapsto (u ,h(u))$. Fix a point
$q = (u_q ,h(u_q)) \in L_0$. The tangent space to $L_0$ at $q$ is spanned over $\C$
by the vector $(1, h'(u_q))$ whose norm is not uniformly bounded on $U$. In any
case $\omega_1 (q)$ evaluated over $(1, h'(u_q))$ coincides with the evaluation of $\Omega_1$ over the same vector. Thus we obtain
\begin{equation}
\omega_1 (q) . [1 , h' (u_q)] =  \Omega_1 (q) . [1, h'(u_q)] = f . h' (u_q) + g \; . \label{siegel6}
\end{equation}
On the other hand, $\omega (q) .  [1 , h' (u_q)] = 0$ so that  Formula~(\ref{siegel3})
provides
\begin{equation}
h' (u_q) [1+r^1 (u_q , h(u_q))] = \frac{\lambda_2 h(u_q)}{\lambda_1 u_q}
(1+ r^2 (u_q , h(u_q))) \; . \label{siegel7}
\end{equation}
Therefore
\begin{eqnarray}
f . h' (u_q) + g & \! = \! & - f\frac{\lambda_2 h(u_q)}{\lambda_1 u_q} \frac{1+r^2}{1+r^1} + g
\label{siegelpr1}\\
 & \! = \! & \frac{-1}{\lambda_1 u_q (1+r^1)} [\lambda_1 (1+a)(1+r^1) \!+\!
 \lambda_1 u_q r^1_u \! - \! \lambda_2(1+b) (1+r^2) \! - \! \lambda_2 h(u_q) r^2_v ]
\label{siegelpr2}
\end{eqnarray}
where the functions $r^1,r^2,r^1_u,r^2_v$ are evaluated at $(u_q , h(u_q))$, cf.
Formula~(\ref{siegel3}). Now recall that $\Vert u_q \Vert$ and
$\Vert h (u_q) \Vert$ are bounded by $\epsilon_1 , \epsilon_2$. It follows from
(\ref{siegel7}) that the norm of $(1 ,h' (u_q))$ is bounded by $\max \{ 1 , {\rm
const}/\Vert
u_q \Vert \}$ for a suitable constant const. The statement then results from
the condition $\lambda_1(1+a) - \lambda_2 (1+b) \neq 0$.\qed

We still need to describe the geometry of the leaves of
$\calh$ on $L_0$. According to Lemma~(\ref{3lema1}), we can suppose without
loss of generality that the oriented leaves of $\calh$ on $\{ v=0 \}$ converge to
$0 \in \{ v=0\} \subset \C^2$ (i.e. $0 \simeq (0,0)$ is a sink for $\calh$ over
$\{ v=0\}$). Similarly $0 \in \{ u=0\} \subset \C^2$ is a source for the leaves of
$\calh$ contained in $\{ u=0 \}$. Next we consider the (real $3$-dimensional)
set
$$
A = \{ (u,v) \in \C^2 \; ; \; \vert u \vert = \epsilon_1 \; {\rm and} \;
\vert v \vert < \epsilon_2 \} \, .
$$
If $\epsilon_1 ,\epsilon_2$ are appropriately chosen and sufficiently small, the the
oriented leaves of $\calh$ point inwards $A$, i.e. at a point $(u,v) \in A$ the leaf of
$\calh$ through this point is oriented in the decreasing direction of the absolute value
of $u$.

Now let $l$ be an oriented leaf of the restriction of $\calh$ to a small neighborhood
of $(0,0) \in \C^2$. As it will shortly be seen, $l$ is not closed. If $q_1 ,q_2 \in l$,
we say that $q_2 > q_1$ provided that one can move from $q_1$ to $q_2$ in the
sense of the orientation of $l$. We also denote by ${\rm dist}\, (q_1 ,q_2)$ the length
of the segment of $l$ whose extremities are $q_1 ,q_2$. Finally we are ready to state
the main result of this section.

\begin{prop}
\label{3prop1}
There is a neighborhood $V$ of $(0,0) \in \C^2$ with the following properties:

\noindent 1. Given $q_1 \in l \cap V$, there is $q_2 \in l \cap A$, with $q_1 > q_2$
and such that ${\rm dist}\, (q_1 ,q_2) < {\rm const}\, \epsilon_1$.

\noindent 2. Given $q_1 \in l \cap V$, there is $\overline{q} = (\overline{q}^1,
\overline{q}^2) \in l$, with $\overline{q} > q_1$ and ${\rm dist}\, (q_1 ,\overline{q})
< {\rm const}\, \epsilon_2$. Besides $\vert \overline{q}^2 \vert = \epsilon_2$ and
$\overline{q}^1 \in \pi_1 (V)$ where $\pi_1 (V)$ stands for the projection of
$V$ on the first coordinate.
\end{prop}

\noindent {\it Proof}\,: Let $B (\delta)$ be the bidisc $\{ (u,v) \in \C^2 \; ; \;
\vert u \vert < \delta \; {\rm and} \; \vert v \vert < \delta \}$. We are going to show
that
$B (\delta)$ satisfies the conditions in our statement provided that $\delta$
is sufficiently small. Consider $a \in B (\delta)$ and suppose without loss of
generality that the real part ${\rm Re}\, (q_1)$ of $q_1$ is positive. Let then $L$
(resp. $l$) be the leaf of the restriction of $\fol$ (resp. trajectory of the restriction of $\calh$) to $B (\delta )$
containing $q_1$. As already seen, $L$ is the graph of a holomorphic function
$h : D_q \subset D_1^- \rightarrow \C$. In the parametrization $u \mapsto (u , h(u))$,
the restriction of $\omega_1$ to $L$ becomes
\begin{equation}
fh' + g = \frac{\lambda_1 (1+a) - \lambda_2 (1+b) + \alpha}{\lambda_1 u}
+s (u) \label{siegel8}
\end{equation}
where $s$ is holomorphic and $\alpha$ can be made arbitrarily small by reducing
$\epsilon_1 ,\epsilon_2$. Indeed Formula~(\ref{siegel8}) is an immediate
reformulation of Formula~(\ref{siegelpr2}). In particular, one has $\lambda_1 (1+a)
-\lambda_2 (1+b) + \alpha \neq 0$. We now set $q_1 = (u_1 , h(u_1))$.
Recalling that $D_q \subset \C$, we denote by $R_q$ the radial line emanated from
$0 \in \C$ and passing through $u_1$. The intersection of $R_q$ with the circle
$\vert u \vert =\epsilon_1$ is $\epsilon_1 u_1 /\vert u_1 \vert$. Similarly let
$\pi_1 (l)$ be the oriented leaf of $\{ {\rm Im}\, (fh'+g) =0 \}$ containing
$u_1$ which is nothing but the projection of $l$ on the first coordinate.

\noindent {\it Claim}\,: There is a point $u_2 \in \pi_1 (l)$ such that
$\vert u_2 \vert =\epsilon_1$. Besides there is a uniform constant $C$ such that
$$
{\rm dist}\, \left( u_2 , \frac{\epsilon_1 u_1}{\vert u_1 \vert}  \right) < C
\epsilon_1^2 \, .
$$

\noindent {\it Proof of the Claim}\,: It is an elementary fact about
continuous/differentiable
dependence of the initial conditions for solutions of real ordinary differential
equations. The foliation associated to $\{ {\rm Im}\, [(\lambda_1 (1+a) - \lambda_2
(1+b))/
\lambda_1 u] =0 \}$ consists of radial lines through $0 \in \C$ so that the assertion
is trivial in this case. Nonetheless the foliation in which we are interested is given by
an Abelian form whose distance to $(\lambda_1 (1+a) - \lambda_2 (1+b))/
\lambda_1 u$ is less than $C \epsilon_1$ for an appropriate constant $C$. The
statement promptly follows.\qed

Combining the above claim with the fact that $\sigma_{\fol} \cup \{ u=0 \}$
contains a neighborhood of $(0,0) \in \C^2$, we conclude that $l$ intersects $A$
at a point $q_2$. Estimates in \cite{mamo} (see also
\cite{mat} and \cite{reis}) guarantee that $q_2$
satisfies the conditions in the statement. Analogously one proves that the continuation
of $l$ intersects the set $\vert v \vert =\epsilon_2$ at a point $\overline{q}$ with
the desired properties. For further details on these estimates we refer the reader
to the quoted papers.\qed

\begin{coro}
\label{hthetahperp}
Under the preceding conditions the trajectories of $\calh^{\perp}$ contained in the local separatrices of $\fol$ are
closed curves encircling the origin. For $\theta \in (-\pi/2 , \pi/2)$, the trajectories of $\calh^{\theta}$ contained in the
local separatrix of $\fol$ where $\calh$ has a sink singularity (resp. a source singularity) are
spiraling curves converging to the origin (resp. being emanated from the origin).

Furthermore, on a local leaf of $\fol$ different from its separatrices, the behavior of $\calh^{\perp}$ is essentially determined
by the local holonomy of the separatrices whereas the behavior of $\calh^{\theta}$, $\theta \in (-\pi/2 , \pi/2)$, is the combination
of the above described Dulac transform (cf. below) with a finite power of the mentioned local holonomy map.
\end{coro}

\noindent {\it Proof}\,: It follows immediately from the fact that the oriented trajectories of $\calh$ (resp. $\calh^{\theta}$)
form an angle of $\pi/2$ (resp. $\theta$) with the oriented trajectories of $\calh$.\qed

Let us close this section with a discussion of the so-called {\it Dulac
transform}\, associated to a singularity in the Siegel domain. Although this is
a local discussion formally independent of the structure of $\calh$, it naturally
involves definitions and results discussed above so that here seems to be a good place
to carry it out. The material below will also be used in Sections~4 and~6. Whereas
classical in nature, it is not easy to find a detailed exposition of this material
in the literature. First we resume some notations.

Recall that $\fol$ is defined on a neighborhood of $(0,0) \in \C^2$ by the
vector field
\begin{equation}
Y = \lambda_1 u (1 +r^1) \frac{\partial}{\partial u} - \lambda_2 v (1 +r^2)
 \frac{\partial}{\partial v} \; . \label{vectorfieldY}
\end{equation}
 Recall also that $A \subset \C^2$ was defined as $A = \{ (u,v) \in \C^2 \;
 ; \; \vert u \vert =\epsilon_1 \; \; {\rm and} \; \; \vert v \vert < \epsilon_2 \}$.
Similarly
 we set $B = \{ (u,v) \in \C^2 \;
 ; \; \vert v \vert =\epsilon_2' \; \; {\rm and} \; \; \vert u \vert < \epsilon_1' \}$
for certain
 $\epsilon_1' , \epsilon_2' >0$. Fixed $u_0$ with $\vert u_0 \vert =\epsilon_1$
 (resp. $v_1$ with $\vert v_1 \vert = \epsilon_2'$), we denote by $\Sigma_0^A$
 (resp. $\Sigma_1^B$) the set  $\{ (u,v) \in \C^2 \; ; \;  u =u_0 \; \; {\rm and} \; \;
\vert v \vert < \epsilon_2 \}$ (resp. $\{ (u,v) \in \C^2 \; ; \;   \vert u \vert <
\epsilon_1'
\; \; {\rm and} \; \; v =v_1 \}$). In the sequel $\epsilon_1 ,\epsilon_1', \epsilon_2'$
are
fixed and small whereas $\epsilon_2$ can be made smaller whenever necessary.

For $u_0 ,v_1$ as above, let us denote by $\fol_0^A, \fol_1^B$ the saturated
of $\Sigma_0^A , \Sigma_1^B$ by $\fol$. As already seen, both $\fol_0^A \cup \{ u=0\}
\cup \{ v=0\}$ and $\fol_1^B \cup \{ u=0\} \cup \{ v=0\}$ contain an open
neighborhood of $(0,0) \in \C^2$. Therefore, up to choosing $\epsilon_2$ very
small, for every $(u_0 ,v_0) \in \Sigma_0^A$, there exist paths $c: [0,1]
\rightarrow L_{(u_0 ,v_0)}$ such that $c(0) = (u_0 ,v_0)$ and
$c(1) \in \Sigma_1^B$ (where $L_{(u_0 ,v_0)}$ stands for the leaf of $\fol$
through $(u_0 ,v_0)$). If $c, c'$ are two paths as above and satisfying
$c(1) = (u ,v_1)$, $c'(1) = (u',v_1)$, then $u,u'$ belong to the same orbit of the local
holonomy of the axis $\{ u=0 \}$.

Now consider a simply connected domain $V_0 \subset \Sigma_0^A \setminus
\{ (u_0 ,0) \}$. Suppose we are given a point $(u_0 ,v_0) \in V_0$ and a path
$c_0 : [0,1] \rightarrow  L_{(u_0 ,v_0)}$ as before. For $(u,v)$ sufficiently close to
$(u_0 ,v_0)$, it is then possible to choose by continuity a path $c : [0,1]
\rightarrow L_{(u_0 ,v_0)}$ such that $c(0) = (u_0 ,v)$ and $c(1) \in \Sigma_1^B$. Since
$V_0$ is simply connected, we can extend this definition to the whole of
$V_0$. In this way, we obtain a holomorphic map ${\rm Dul}: V_0 \subset
\Sigma_0^A \setminus \{ (u_0 ,0) \}$ to $\Sigma_1^B$. This map is going to
be called the {\it Dulac transform}\, (which depends on the previously chosen
path $c_0$). Identifying $\Sigma_0^A$ with a neighborhood of $0 \in \C$, we
shall refer to a sector of angle $\theta$ and radius $r$ meaning the intersection
of the ball of radius~$r$ with a sector of angle $\theta$ (and vertex at $0 \in \C$).
In practice, $V_0$ will always be a sector of angle less than $2\pi$ and sufficiently
small radius. The choice of the initial path $c_0$ and of the semi-line in
question entirely determines the corresponding map ${\rm Dul}$.

The following lemma consists again of estimates that can be found
for example in \cite{mamo}, \cite{mat} or in \cite{reis}.

\begin{lema}
\label{3lema3}
Let $V_0 = \Sigma_0^A$ be a sector of angle less than $2\pi$ and sufficiently
small radius. Fix a path $c$ and consider the resulting Dulac transform
${\rm Dul} : V_0 \subset \Sigma_0^A \rightarrow \Sigma_1^B$. Then the following estimate
holds
$$
\Vert {\rm Dul}\, (v) \Vert \leq {\rm Const}\, \Vert v \Vert^{\lambda_1 /\lambda_2} (1 +
O
(\Vert v \Vert)) \; .
$$
\noindent \mbox{ }\qed
\end{lema}

In particular, if $\lambda_1 > \lambda_2$, the behavior of ${\rm Dul}$ is that
of a (strong) contraction provided that $\Vert v \Vert$ is small. When $\lambda_1
< \lambda_2$ then ${\rm Dul}$ behaves as an expansion for $\Vert v \Vert$ small.

Finally suppose that $\Sigma_0^A , \Sigma_1^B$ are endowed with measures
$\mu_0 , \mu_1$ which are part of a system (of transverse sections and measures)
defining a transverse invariant measure for a global realization of $\fol$ on some
complex surface (in the sense of Section~2.2). Note that, in general, ${\rm Dul}$
is not one-to-one on $V_0 \subset \Sigma_0^A$ (if the angle of $V_0$ is not
small) so that
$\mu_0 (V_0) \neq \mu_1 ({\rm Dul}\, (V_0))$. Nonetheless we have:

\begin{lema}
\label{3lema4}
With the preceding notations the following is verified.

\noindent 1. Suppose that $\lambda_1 > \lambda_2$ and let $V_0$
be a sector of angle slightly less than $2\pi \lambda_2 /\lambda_1$. Then, for
$\Vert v \Vert$ very small, ${\rm Dul}$ is one-to-one on $V_0$ and satisfies
$\mu_0 (V_0) = \mu_1 ({\rm Dul}\, (V_0))$.

\noindent 2. Suppose that $\lambda_1 < \lambda_2$ and let $V_0$
be a sector of angle slightly less than $2\pi$. Then, for
$\Vert v \Vert$ very small, ${\rm Dul}$ is one-to-one on $V_0$ and satisfies
$\mu_0 (V_0) = \mu_1 ({\rm Dul}\, (V_0))$.
\end{lema}

\noindent {\it Proof}\,: The proof consists of showing that for
$v \in W$, we can obtain
flow boxes containing  the corresponding paths $c: [0,1] \rightarrow
L_{(u_0 ,v)}$
so that the holonomy associated to these paths is well-defined and injective.
This is clear when $\fol$ is linearizable. In the general case it results again from
the asymptotic estimates already mentioned above.\qed

\begin{obs}
\label{holomorphiconaxes}
{\rm {\bf The case when the restriction of $\omega_1$ to the local separatrices is holomorphic}:
the reader has noted that the discussion of the behavior of $\calh$ (resp. $\calh^{\perp}$ and $\calh^{\theta}$) carried out in Proposition~\ref{3prop1}
was based on
the local invariance condition and on the assumption that $\lambda_1 (1+a) - \lambda_2 (1+b) \neq 0$. Now that we have already
introduced the notion of Dulac transform associated to a Siegel singularity, let us also consider the case where $\lambda_1 (1+a) - \lambda_2 (1+b) = 0$
(assuming that the local invariance condition is still satisfied). As mentioned this case is such that the
restriction of $\omega_1$ to the invariant axes $\{y=0\}$ and $\{x=0\}$ is holomorphic. Thus the restriction of $\omega_1$ to
$\{y=0\}$ (resp. $\{x=0\}$) either is regular or vanishes at the origin. For the time being we shall assume that this restriction is not identically zero,
though this is not strictly necessary for what follows (cf. Sections~4 and~5). Consider then the behavior of $\calh$ restricted to $\{y=0\}$ and suppose
there is a trajectory $l$ of $\calh$ that passes ``very close'' to the origin. The first remark to be made here is that $l$ can be ``deformed'' to avoid
a fixed neighborhood of the origin. These deformations are similar to deformations already performed when a singularity converges to a saddle-singularity
of $\calh$ occurring at a regular point of $\fol$, cf. Section~2 and/or \cite{blm}. In particular they can be done without destroying the ``contractive behavior''
of the holonomy of $\fol$ associated to the trajectories of $\calh$. Therefore, if needed, a Siegel singularity satisfying the condition
$\lambda_1 (1+a) - \lambda_2 (1+b) = 0$ can be avoided by the trajectories of $\calh$. In
other words, the singularity becomes ``invisible'' and thus it can be ignored.

However, even if these singularities can be avoided, we might want to take advantage of them by exploiting the (local) saddle-behavior of $\fol$.
In other words, it may be useful to let a $\calh$-trajectory to approximate the singularity so as to be continued ``through the other separatrix of $\fol$'', i.e.
the $\calh$ trajectory may go through the Dulac transform and then be continued in a different way. In this paper, if a trajectory of $\calh$ is about to entering some
(previously fixed) neighborhood of a Siegel singularity as above, we shall consider all possible continuations of it, namely those that actually ``avoid the singularity''
and those that passes through the Dulac transform associated to the singularity itself. We shall return to these cases
later in Sections~5 and~6.}
\end{obs}

\noindent {\bf An alternative point of view}: let us close this section by explaining an alternate way to see the above results on Dulac
transforms and their connections with the material of Section~2.3.

To begin with, let us make a simple remark concerning how the Dulac transform can be viewed in most of our applications.
With the preceding notations suppose that the orientation of the trajectories of $\calh$ is such that the origin is a sink for the restriction
of $\calh$ to $\{v=0\}$. Then $\epsilon_1, \epsilon_2$ can be chosen so that $\calh$ is transverse to $A\subset \C^2$. Besides every
$\calh$-trajectory intersecting $A$ points inward $A$ and, unless this intersection occurs at a point belonging to $\{v=0\}$,
it will eventually intersect $B$ with outward orientation. Thus we can define the
Dulac transform as being the map from $A$ to $B$ defined by the trajectories of $\calh$. Note that this map is locally holomorphic away from
$A \setminus \{v=0\}$. Besides, for $(u_0, v_0) \in A$, $v_0 \neq 0$, its image satisfy the estimates given in Lemma~\ref{3lema3}. Furthermore it is
not hard to adapt the contents of Lemma~\ref{3lema4} to this setting.

Naturally the preceding statements about the contractive or expansive character of the Dulac map can also be viewed in terms of Poincar\'e Lemma discussed in Section~2.3.
For this it is however necessary to work with (possibly) ramified
coordinates. Let us then consider a foliation $\fol$ defined on a neighborhood of $(0,0) \in \C^2$ by the vector field $Y$ in~(\ref{vectorfieldY}).
More precisely suppose that the $1$-form $\omega$ defining $\fol$ is actually
$\omega =  \lambda_1 u (1 +r^1) \, dv + \lambda_2 v (1 +r^2) \, du$. Consider also sections $\Sigma_0^A$ and $\Sigma_1^B$ as above and suppose
that the orientation of the trajectories of $\calh$ is such that they go from $\Sigma_0^A$ to $\Sigma_1^B$ (i.e. $\lambda_1 > \lambda_2$). To apply
Formula~(\ref{PLemma}) to this case, we need to consider the parametrizations of $\Sigma_0^A, \, \Sigma_1^B$ that are obtained through the integral
of $\omega$. It is then natural to set a coordinate $z_1$ on $\Sigma_0^A$ and a coordinate $z_2$ on $\Sigma_1^B$ such that
$$
z_1 = \lambda_2 v (1 + {\rm h.o.t.}) \; \; \, {\rm and} \; \; \, z_2 = \lambda_1 u (1 + {\rm h.o.t.}) \, .
$$
In these coordinates the derivative of the above introduced Dulac transform can be estimate by means of Formula~(\ref{PLemma}).
This amounts to estimating the integral of $\omega_1$ over a segment of trajectory of $\calh$ going from $\Sigma_0^A$ to $\Sigma_1^B$.
The latter estimate however is essentially equivalent to the calculations performed above.


\section{Singularities of $\fol$ and invariant measures}

Now we are going to begin the analysis of the global setting where $\fol$ is a singular holomorphic foliation defined on
a complex surface $M$. Throughout this section $\fol$ is
supposed to admit an invariant positive closed current $T$ whose associated
transverse measure does not give mass to points ($T$ is said to be diffuse).
Let $\supT \subseteq M$ be the {\it support}\, of $T$ which is obviously a
compact set invariant by $\fol$.

Modulo applying Seidenberg's theorem, we can suppose that all the
singularities of $\fol$ are reduced. Our first aim in this section is to
establish Proposition~(\ref{4.5prop1}) below.

\begin{prop}
\label{4.5prop1}
Let $p \in {\rm Sing}\, (\fol)$ be a singularity of $\fol$ lying in $\supT$.
Then $p$ is a singularity in the Siegel domain or it is an irrational focus.
Furthermore if $p$ belongs to the Siegel domain and has
eigenvalues with rational quotient, then $\fol$ is linearizable around $p$.
\end{prop}

\noindent Since $p \in {\rm Sing}\, (\fol) \cap \supT$ is reduced, the proof of
Proposition~(\ref{4.5prop1})
essentially consists of showing that $p$ is neither a hyperbolic singularity nor
a saddle-node. These are the contents of Lemmas~(\ref{4.5lema1})
and~(\ref{4.5lema2}) below.

\begin{lema}
\label{4.5lema1}
If $p \in {\rm Sing}\, (\fol) \cap \supT$, then $p$ is not hyperbolic.
\end{lema}

\noindent {\it Proof}\,: Suppose for a contradiction that $p$ is hyperbolic. Then
Poincar\'e Theorem ensures that $\fol$ is linearizable around $p$. In other
words, there are local coordinates $u,v$ in which $\fol$ is given by
$$
\eta = \lambda_1 u dv - \lambda_2 v du
$$
with $\lambda_1 /\lambda_2 \in \C \setminus \R$. Consider a local transverse
section $\Sigma$ passing through the point $(1,0)$. This section allows us
to identify the local holonomy of the separatrix $\{ v=0\}$ with a local diffeomorphism
$h$ fixing $0 \in \C$. The condition $\lambda_1 /\lambda_2 \in \C
\setminus \R$ implies that $h$ is hyperbolic, i.e. $\vert h'(0) \vert <1$. Now consider
a (local) leaf $L$ of $\fol$ contained in $\supT$ and intersecting $\Sigma$ at a
point $(1, z_0)$. Denote by $\mu_{\Sigma}$ a representative of $T$, viewed as
transverse invariant measure, over $\Sigma$ (cf. Section~2.2).  If $z_0 \neq 0$,
the orbit of $(1,z_0)$ under $h$ consists of infinitely many points converging
towards $(1,0) \in \Sigma$. Furthermore, if $V \subset \Sigma$ is a sufficiently
small neighborhood of $z_0 \simeq (1,z_0) \in \Sigma$, then the open sets
$V ,h(V), h^2 (V), \ldots$ are pairwise disjoint. Nonetheless they have all the
same $\mu_{\Sigma}$ measure for $h$ preserves $\mu_{\Sigma}$. In addition
$\mu_{\Sigma} (V) >0$ since $L$ is contained in $\supT$.
Together these facts imply that $\mu_{\Sigma} (\Sigma) =\infty$ what is impossible.
We then conclude that $\supT$ is locally contained in the separatrices of $\fol$
at $p$. Therefore $\mu_{\Sigma}$ has an atomic component which is necessarily
concentrated over an algebraic curve. Since this is impossible, the lemma
follows.\qed

Through a similar argument we are going to prove that $p \in {\rm Sing}\, (\fol)
\cap \supT$ cannot be a saddle-node either. A very complete reference for
saddle-node singularities is \cite{ramis}. The facts used below are however
well-known. If $\fol$ is a saddle-node singularity then it can be written
in Dulac Normal Form, i.e. in suitable local coordinates $u,v$, the foliation $\fol$
is given by the $1$-form $\eta$ satisfying
$$
\eta = [u(1+ \Lambda v^p) + R(u,v)] dv - v^{p+1} du \; \;  \; \; {\rm with} \; \; \; \;
\Lambda \in \C \; \; \; \; {\rm and} \; \; \; \; p \geq 1 \, .
$$
In particular $\{ v=0 \}$ is a separatrix of $\fol$ called the {\it strong invariant
manifold of $\fol$}. Considering a local transverse $\Sigma$ as in
Lemma~(\ref{4.5lema1}), we can identify the holonomy of the strong
invariant manifold to a (local) diffeomorphism $h$ fixing $0 \in \C$. However,
this time, $h$ has the form $h (z) = z + z^{p+1}  + {\rm h.o.t.}$, where as usual ${\rm h.o.t}$
stands for terms of higher order.

\begin{lema}
\label{4.5lema2}
If $p \in {\rm Sing}\, (\fol) \cap \supT$, then $p$ cannot be a saddle-node.
\end{lema}

\noindent {\it Proof}\,: Consider $\fol ,\Sigma$ and $\eta$ as above. Other than the
strong invariant manifold, a saddle-node may or may not possess another
separatrix (necessarily smooth and transverse to the former one) which is called
the weak invariant manifold. In particular a saddle-node possesses at least one
and at most two separatrices.

We now suppose that $\supT$ is not locally contained in the union of the
separatrices of $\fol$ since this would again lead us to a contradiction.
It follows from \cite{ramis} that the union of the saturated
$\fol_{\Sigma}$ of
$\Sigma$ by $\fol$ with the weak invariant manifold (if it exists) contains
a neighborhood of $p$. Thus there is a leaf $L \subset \supT$ of $\fol$ intersecting
$\Sigma$ at a point $(1,z_0)$ with $z_0 \neq 0$. The topological description
of the dynamics of $h (z) = z + z^{p+1}  + {\rm h.o.t.}$ is well-known (cf. for example
\cite{flower}) and it follows the existence of a small neighborhood $V \subset
\Sigma$ of $z_0 \simeq (1,z_0)$ such that $V, h (V)  ,h^2 (V) \ldots$ are pairwise
disjoint. By taking a representative $\mu_{\Sigma}$ of $T$ on $\Sigma$ as
in Lemma~(\ref{4.5lema1}) we conclude that $\mu_{\Sigma} (\Sigma) =\infty$.
This is however impossible and establishes the lemma.\qed

\vspace{0.1cm}

\noindent {\it Proof of Proposition~(\ref{4.5prop1})}\,: After Lemmas~(\ref{4.5lema1})
and~(\ref{4.5lema2}), we only need to prove that a Siegel singularity
with rational eigenvalues is linearizable. As already seen $\fol$ is locally given by
$$
\eta = \lambda_1 u (1 + {\rm h.o.t.}) \, dv  + \lambda_2 v (1 + {\rm h.o.t.}) \, du
$$
with $\lambda_1 /\lambda_2 \in \Q_+$. Denoting by $\Sigma$ a transverse
section passing through $(1,0)$, it was seen that the union of $\{ u=0\}$ with
$\fol_{\Sigma}$ (the saturated of $\Sigma$ by $\fol$) contains a neighborhood of $p$.
Thus, as before, there is a leaf $L \subset \calk$ intersecting $\Sigma$
at a point $(1 ,z_0)$. Without loss of generality we can suppose that $z_0 \neq 0$.

On the other hand, the linear part of the holonomy diffeomorphism $h$ associated
to $\{ v=0\}$ is precisely $e^{2\pi i\lambda_1 /\lambda_2}z$. Thus a power
of $h$ is tangent to the identity. According to a result of Mattei-Moussu \cite{mamo},
$\fol$ is locally linearizable if and only if the power of $h$ in question coincides
with the identity. Hence we suppose for a contradiction that this power
is tangent to the identity and different from the identity. In this case, however,
it has the form $z + cz^k + \cdots$ with $c\neq 0$. The final contradiction is then
obtained as at the end of Lemma~(\ref{4.5lema2}). The proposition is proved.\qed

Summarizing the preceding discussion we can suppose that the (reduced) singularities
of $\fol$ lying in $\supT$ are of one of the following types:

\noindent $\bullet$ a singularity in the Siegel domain.

\noindent $\bullet$ an irrational focus.

\noindent Note also that Poincar\'e Theorem still implies that an irrational
focus is automatically linearizable. Hence, in this case, $\fol$ is locally
given by the form
$$
\eta = \lambda_1 u dv - \lambda_2 v du
$$
with $\lambda_1/\lambda_2 \in \R_+\setminus \Q_+$. It is easy to work out the
structure of the foliation $\calh$ near to an irrational focus singularity. This is similar to
the discussion carried out in Section~3 whereas technically simpler since $\fol$ is always linearizable.
Again on a small neighborhood of $p$ the
curves $(\omega)_0$ and $(\omega)_{\infty}$ are supposed to be invariant by $\fol$ (local invariance condition).
This means that $\omega$ can be written in local coordinates $u,v$ as
\begin{equation}
\omega = h (u,v) u^a v^b [ \lambda_1 u dv - \lambda_2 v du]  \label{4.5eq1}
\end{equation}
where $h(0,0) \neq 0$. Setting $\Omega_1 = fdv + gdu$ the equation $d\omega =
\omega \wedge \Omega_1$ yields
\begin{equation}
h (u,v) ( \lambda_2 v f + \lambda_1 ug) = -h(u,v) (\lambda_1 (a+1) + \lambda_2 (b+1))
- u \frac{\partial h}{\partial u} - v\frac{\partial h}{\partial v} \, . \label{4.5eq2}
\end{equation}
In the sequel we suppose that $a,b$ are not simultaneously equal to~$-1$ so
that $\lambda_1 (a+1) + \lambda_2 (b+1) \neq 0$ (recall that $\lambda_1/\lambda_2 \in
\R_+ \setminus \Q_+$). By setting $u=0$ (resp. $v=0$) we conclude that the behavior
of $h$ over the separatrix $\{ u=0 \}$ (resp. $\{ v=0\}$) is either that of a sink or
that of a source according to whether $\lambda_1 (a+1) + \lambda_2 (b+1) >0$ or
$\lambda_1 (a+1) + \lambda_2 (b+1) <0$.

For the leaves of $\fol$ different from the separatrices, we can perform a
discussion similar to the one carried out in Section~3 by exploiting the presence of the ``multi-valued'' first integral
$u^{\lambda_2} v^{\lambda_1}$. The reader will easily check that the behavior
of $\calh$ over the separatrices is repeated over the general leaves. The result
is then summarized by

\begin{prop}
\label{4.5lema3}
Let $p \in {\rm Sing}\, (\fol ) \cap \supT$ be an irrational focus. Consider also local
coordinates $u,v$ defined on a bidisc of radius $\epsilon$
about $p$ and suppose that $\omega$ is given by~(\ref{4.5eq1}) where $a,b$ are not simultaneously equal to~$-1$. If $L$
is a leaf of $\fol$, then the restriction of $\calh$ to $L$ consists of lines of length
less than ${\rm Const}. \epsilon$ for an appropriate uniform constant ${\rm Const}$.
Furthermore these lines converge to $(0,0)$ if $\lambda_1 (a+1) + \lambda_2 (b+1) >0$
(i.e. the end of the leaf correponding to $(0,0)$ is a sink). Similarly
these lines are emanated from $(0,0)$ if $\lambda_1 (a+1) + \lambda_2 (b+1) <0$
(i.e. the end of the leaf correponding to $(0,0)$ is a source).\qed
\end{prop}

\begin{obs}
\label{4.5obs1}
{\rm An irrational focus $p \in {\rm Sing}\, (\fol ) \cap \supT$ is going to be
called a sink (resp. a source) if, with the notations of the lemma above, one
has $\lambda_1 (a+1) + \lambda_2 (b+1) >0$ (resp. $\lambda_1 (a+1) + \lambda_2 (b+1) <0$).
Sometimes we shall use the expressions sink-irrational focus or source-irrational
focus to emphasize that we are dealing with an irrational focus singularity.
This terminology also serves to distinguish between singularities of $\fol$ behaving as
sinks (or sources) for $\calh$ and sinks (or sources) of $\calh$ occurring at
regular points of $\fol$.}
\end{obs}

To close this section, we are going to introduce a sort of ``generalized Dulac
transform'' (or maybe ``compounded Dulac transform'')
for the singularities of the foliation $\fol$.
This material will be needed
in Section~6 since the singularities of the initial foliation $\fol$ (as in the
statement of Theorem~A in the Introduction) may be degenerate. Also it
should be pointed out that Proposition~(\ref{4.5prop1}) is not
used in the
following discussion although it will be necessary in Section~6.
In fact, the role played by
Proposition~(\ref{4.5prop1}) in Section~6 amounts to guaranteeing that the
situation considered in the discussion below always occurs. In particular, this will
enable us to consider the ``generalized Dulac transform'', cf. below.

To explain our concern with this ``generalized Dulac transform'', consider
the local situation given by a singularity of $\fol$ that belongs to the
Siegel domain. Let $\lambda_1, \lambda_2$ be the eigenvalues of $\fol$ at
$p$ and suppose that $\lambda_1 > \lambda_2$. Suppose in addition that $p$
lies away from the divisor $(\omega)_0 \cup (\omega)_{\infty}$ of zeros
and poles of $\omega$, where $\omega$ stands for a meromorphic $1$-form
defining $\fol$. Let $\textsc{S}_1, \textsc{S}_2$ denote the separatrices
of $\fol$ at $p$ that are respectively tangent to the eigendirections associated
to $\lambda_1, \lambda_2$. According to the discussion in Section~3, the
restriction of $\calh$ to $\textsc{S}_1$ consists of trajectories converging
to $p$. Similarly, the restriction of $\calh$ to $\textsc{S}_2$
consists of trajectories
emanated from $p$. Thus, if $l$ is a segment of $\calh$-trajectory passing near $p$,
the Dulac transform defined by means of $l$ behaves as a contraction
(cf. Section~3 and Lemma~\ref{3lema3}). The existence of this contraction is therefore
consistent with the principle of producing ``contractive holonomy'' by following the
trajectories of $\calh$. However, if $\textsc{S}_1, \textsc{S}_2$ are contained in the divisor
$(\omega)_0 \cup (\omega)_{\infty}$, then the orientation of $\calh$ around
$p$ may be ``unnatural'' in the sense that the Dulac transform induced by a
segment of $\calh$-trajectory as above actually behaves as an expansion (cf.
Lemma~\ref{3lema3}). The tension between contraction along the leaves of
$\calh$ and expansion for certain Dulac transforms would prevent us from guaranteeing the existence of a
contractive holonomy map in a suitable sense. It is to remedy this situation
that ``generalized Dulac transforms'' will be introduced. The aim of their study is show that
contraction eventually prevails.

Without loss of generality, we can assume that $\fol$ is a foliation
with reduced singularities defined on a certain compact surface. We also
fix a non-closed meromorphic $1$-form $\omega$ defining $\fol$ (which is supposed to
exist in our case). Let $(\omega)_0^{\perp\fol}$
(resp. $(\omega)_{\infty}^{\perp \fol}$) be the subdivisor of $(\omega)_0$
(resp. $(\omega)_{\infty}$) consisting of those irreducible components of $(\omega)_0$
(resp. $(\omega)_{\infty}$) that {\it are not}\, invariant by $\fol$. As before we set
$(\omega)_0^{\fol} = (\omega)_0 \setminus (\omega)_0^{\perp\fol}$ and
$(\omega)_{\infty}^{\fol} = (\omega)_{\infty}^{\fol} \setminus
(\omega)_{\infty}^{\perp \fol}$. Let $E$ be a connected component of
$(\omega)_0^{\fol} \cup (\omega)_{\infty}^{\fol}$. Modulo performing finitely
many blow-ups, we can assume without loss of generality that
that $(\omega)_0^{\perp \fol}$ (resp. $(\omega)_{\infty}^{\perp \fol}$) intersects $E$
only at regular points of $\fol$ (cf. Lemma~\ref{revision2} in Section~5 for a detailed explanation
of this procedure). The irreducible components of $E$ are going to be
denoted by $D_1, \ldots ,D_n$.

Let us now consider a leaf $L$ of $\fol$ that accumulates on
a singularity $\textsc{P}_0 \in
D_1 \subseteq E$. We suppose that $\textsc{P}_0$ belongs to the
Siegel domain and that $L$
is not locally contained in the separatrices of $\fol$ at $\textsc{P}_0$.
One of these separatrices, $\textsc{S}^{P_0}$, of $\fol$ at
$\textsc{P}_0$ is transverse to $E$ (and thus not contained in $E$).
The other
separatrix of $\fol$ at $\textsc{P}_0$ is obviously contained in $D_1
\subset E$. Next suppose
we are given a sequence of singularities of $\fol$ in $E$ verifying the two
conditions below:
\begin{enumerate}
\item Each singularity belongs to the Siegel domain.

\item Each singularity corresponds to the intersection of two irreducible
components of $E$ (recall that $E$ is already totally invariant by $\fol$).
\end{enumerate}

\noindent The above mentioned sequence of singularities will be denoted
by $\{ p_1, \ldots ,p_k \}$. We suppose that $p_k$ belongs to a component
$D_l$ of $E$ (note that $l$ may differ from $k$ since the Dynkin diagram
of $E$ is allowed to contain loops).
Finally one still has a singularity $\textsc{P}_1 \in
D_l$ belonging to the Siegel domain and having a
separatrix $\textsc{S}^{P_1}$ transverse
to $E$ (the other separatrix of $\fol$ at $\textsc{P}_1$ being contained in
$D_l \subset E$). Let $\Sigma_0 ,\Sigma_1$ be local transverse sections at points
$z_0 \in \textsc{S}^{P_0}$ and $z_1 \in \textsc{S}^{P_1}$, respectively. Denote
by $\mu_0 , \mu_1$ measures on $\Sigma_0 ,\Sigma_1$ representing $T$
over these transversals (as in Lemma~\ref{3lema4}).

We want to define the ``generalized Dulac transform'' ${\rm GDul}$ from a domain
$W \subset \Sigma_0$ to $\Sigma_1$. This can naturally be done by composing
the (ordinary) Dulac transforms associated to the singularities $\textsc{P}_0, p_1
, \ldots ,p_k, \textsc{P}_1$. Proposition~(\ref{4.5prop2}) below makes this
definition precise and collect the properties of ${\rm GDul}$ that are going to
be used in Section~6.

Keeping the preceding notations, we have two further assumptions.

\noindent 3. All the singularities $\textsc{P}_0, p_1, \ldots ,p_k,
\textsc{P}_1$ satisfy the condition $\lambda_1 (1+a) -\lambda_2 (1+b) \neq 0$ of Lemma~(\ref{3lema1})
and subsequent ones in Section~3.

\vspace{0.1cm}

\noindent 4. The trajectory $l_{z_0}$ of $\calh$ through $z_0 =\Sigma_0 \cap
\textsc{S}^{P_0}$ converges to $\textsc{P}_0$. It then continues to $p_1$
and from $p_1$ to $p_2$ and so on until it reaches $\textsc{P}_1$. From
$\textsc{P}_1$ this trajectory leaves $E$ (and thus a small tubular neighborhood of $E$)
by following the separatrix
$\textsc{S}^{P_1}$. This trajectory is also assumed
to pass through $z_1 = \Sigma_1 \cap \textsc{S}^{P_1}$.

\noindent For a detailed definition of the trajectories of $\calh$ ``passing
through singularities in the Siegel domain'', the reader is referred to the discussion carried out
in Section~5. The definition of ${\rm GDul}$ simply consists of the composition of
Dulac transforms associated to the singularities in question with ordinary holonomy
maps associated to the segments of the leaf of $\calh$ between two such singularities.
Finally we have:

\begin{prop}
\label{4.5prop2}
Under the preceding assumption, there is $1> \lambda >0$ with the following properties:
\begin{enumerate}
\item If $V_0 \subset \Sigma_0$ is a sector of angle less that $2\pi \lambda$ and
sufficiently
small radius, then ${\rm GDul}: V_0 \rightarrow \Sigma_1$ is well-defined and one-to-one.

\item For $v \in V_0$, one has $\Vert {\rm GDul}\, (v) \Vert \sim O (\Vert v
\Vert^{1/\lambda})$.
Therefore ${\rm GDul}$ is a contraction for $\Vert v \Vert$ small.

\item One has $\mu_0 (V_0) = \mu_1 ({\rm GDul}\, (V_0))$, provided that the
radius of $V_0$ is small enough.
\end{enumerate}
\end{prop}

\noindent {\it Proof}\,: The statement is clear if the divisor $E$ is empty
as an already mentioned consequence of the combination of
Lemmas~(\ref{3lema1}), ~(\ref{3lema3}) and~(\ref{3lema4}).

Consider now the case $k=0$, i.e. both $\textsc{P}_0$ and $\textsc{P}_1$
belong to $D_1$. Denote by $\lambda^0_1 ,\lambda_2^0$ (resp.
$\lambda^1_1 ,\lambda_2^1$) the eigenvalues of $\fol$ at $\textsc{P}_0$
(resp. $\textsc{P}_1$) where $\lambda^0_2$ (resp. $\lambda^1_2$) is the
eigenvalue associated to the eigendirection defined by $D_1$. Now let $b \in \Z$ be
the order if $D_1$ as a component of the divisor of zeros and poles
of $\omega$. Note that the separatrices $\textsc{S}^{P_0}, \textsc{S}^{P_1}$
are not locally contained in the support of this divisor since they are
transverse to $E$ (cf. the definition of $E$).
Hence there are local
coordinates $(u,v)$ (resp. $(w,v)$) around $\textsc{P}_0$ (resp. $\textsc{S}^{P_1}$)
in which $\omega$ can be written as
\begin{eqnarray}
\omega & =  & v^b [\lambda_2^0 u (1 + {\rm h.o.t.}) dv + \lambda_1^0 v (1 + {\rm h.o.t.})
du] \, ,
\label{aqui1} \\
\omega & =  & v^b [\lambda_2^1 w (1 + {\rm h.o.t.}) dv + \lambda_1^1 v (1 + {\rm h.o.t.})
dw]
\label{aqui2}
\end{eqnarray}
where $\{ v=0\} \subset D_1$. Since the trajectories of $\calh$ converge towards
$\textsc{P}_0$,
Lemma~(\ref{3lema1}) ensures that $\lambda^0_2 - \lambda^0_1 (b+1) <0$. Similarly,
because those trajectories also leave $\textsc{P}_1$ along $\textsc{S}^{P_1}$,
Lemma~(\ref{3lema1}) provides, in addition, that $\lambda_2^1 - \lambda_1^1 (b+1) >0$.
On the other hand,
recall that eigenvalues are defined only up to a multiplicative constant, so that
we can set $\lambda^0_2 = \lambda_2^1$. It then results that $\lambda_1^0 > \lambda_1^1$.
The corresponding generalized Dulac transform, however, clearly satisfies
$\Vert {\rm GDul}\, (1,u) \Vert = O \Vert u \Vert^{\lambda^0_1 /\lambda^1_1}$. Therefore
${\rm GDul}$ has the desired contracting behavior for $\Vert u \Vert$ small.
The second part of the statement can directly be checked. Indeed, there are only two
cases according to whether or not $\lambda_2^0 \geq \lambda_1^0$ (in any case
we have $\lambda_2^1 \geq \lambda_1^1$ provided that $b \geq 0$). This verification
is left to the reader.

Let us now consider the case $k=1$. The new element appearing in this situation
is the singularity $p_1 = D_1 \cap D_2$. Keeping similar notations, let $b_1$
(resp. $b_2$) denote the order of $D_1$ (resp. $D_2$) as a component of the divisor
of zeros and poles of $\omega$. The eigenvalues of $\fol$ at $\textsc{P}_1, \textsc{P}_2$
are still denoted as before. Finally let $\Lambda_1$ (resp. $\Lambda_2$) be the
eigenvalue of $\fol$ at $p_1$ associated to the eigendirection given by $D_1$
(resp. $D_2$). Around $\textsc{P}_0$, there are local coordinates $(u,v)$ where
$\omega$ is given as in~(\ref{aqui1}) (with $b=b_1$). Around $\textsc{P}_1$, we have
local coordinates $(w,t)$, $\{ t=0\} \subset D_2$, where $\wto$ becomes
\begin{equation}
\omega = t^{b_2}  [\lambda_2^1 w (1 + {\rm h.o.t.}) dt+ \lambda_1^1 v (1 + {\rm h.o.t.})
dw] \, .
\label{aqui3}
\end{equation}
Once again Lemma~(\ref{3lema1}) gives us that $\lambda_2^0 - \lambda_1^0 (b_1 +1)
<0$ and $\lambda_2^1 - \lambda_1^1 (b_2 +1) >0$. Finally, in the coordinates
$(v,t)$ around $p_1$, we obtain
$$
\omega = b^{b_1} t^{b_2} [\Lambda_1 t  (1 + {\rm h.o.t.}) dv+ \Lambda_2 v (1 + {\rm
h.o.t.}) dt] \, .
$$
Thanks to Lemma~(\ref{3lema1}), we know that $\Lambda_1 (b_2 +1) - \Lambda_2
(b_1 +1) >0$. To conclude, we first observe that we can set $\Lambda_1 = \lambda_2^0$
and $\Lambda_2 = \lambda^1_2$ since these eigenvalues are defined only
up to a multiplicative constant. Therefore one has
$$
\lambda_1^0 (b_1+1)(b_2+1) > \lambda^0_2 (b_2 +1) > \lambda_2^1 (b_1 +1)
> \lambda_1^1 (b_1+1)(b_2+1)
$$
so that $\lambda_1^0 > \lambda_1^1$. In other words, the generalized Dulac transform
has the contracting behavior indicated in the statement. Again the verification
of item~3 is left to the reader.

The general case of $k\in \N$ now follows easily by induction.\qed

Actually our proof yields a slightly more general result. To state it let us drop Condition~3 above, i.e. the singularities
$\textsc{P}_0, p_1, \ldots ,p_k, \textsc{P}_1$ need no longer to satisfy the condition $\lambda_1 (1+a) -\lambda_2 (1+b) \neq 0$. If
$p_i$ is a (Siegel) singularity at which we have $\lambda_1 (1+a) -\lambda_2 (1+b) = 0$, then the restrictions of $\omega_1$ to the local
separatrices of $\fol$ at $p_i$ are holomorphic on a neighborhood of $p_i$. This setting includes the case in which the restriction of $\omega_1$
to one (or to both) of these separatrices vanishes identically. Our purpose here is to allow the Dulac transform corresponding to $p_i =D_i
\cap D_{i+1}$ to be considered (with orientation going from $D_i$ to $D_{i+1}$) as a component in the constitution of the generalized Dulac
transform. The reader will note that the occasional use of the Dulac map in question is consistent with the contents of Remark~\ref{holomorphiconaxes}
and it will further be detailed in the next section. Naturally away from the singularities that fail to fulfill the condition $\lambda_1 (1+a) -\lambda_2 (1+b) \neq 0$
we shall always follow the trajectories of $\calh$. Then the proof of Proposition~\ref{4.5prop2} can be repeated word-by-word to provide:

\begin{coro}
\label{4.5prop2PRIME}
Under the preceding assumption the statement of Proposition~\ref{4.5prop2} still holds except that now
$1 \geq \lambda >0$. Besides if $\lambda=1$ then ${\rm GDul}$ is defined on every sector $V_0  \subset \Sigma_0$
of angle less than~$2\pi$ (and sufficiently small radius). In the latter case the generalized Dulac transform ${\rm GDul}$
is asymptotically flat at the ``origin of $\Sigma_0$''.
\end{coro}

\section{Topological dynamics of the trajectories of $\calh$}

In the preceding two sections, we have studied the local behavior of $\calh$
around singularities of $\fol$. It is now time to make global considerations on these trajectories.

In what follows we consider a holomorphic foliation $\fol$ given by a globally
defined meromorphic form $\omega$ on a compact
surface $M$. As always we suppose that $\omega$ is not closed and that
$\fol$ admits an invariant {\it diffuse}\, positive closed current $T$. Again $\supT$
will denote the support of $T$. Thanks to Seidenberg Theorem, we can assume without
loss of generality that the singularities of $\fol$ are all reduced. By virtue of
Proposition~(\ref{4.5prop1}) this, in fact, implies that the singularities of
$\fol$ in $\supT$ either belong to the Siegel domain or are irrational foci.
It is also known that a singularity of $\fol$ in $\supT$ belonging to the Siegel
domain is automatically linearizable provided that the quotient of its eigenvalues is rational.

As already explained, our strategy consists of following the trajectories of $\calh$
with the purpose of guaranteeing a ``contractive behavior for the corresponding holonomy maps''. If ``enough
contraction'' is obtained then we should be able to conclude that $T$ is
the current of integration over a compact leaf (cf. for example Lemma~\ref{atomicmass}).
It should be noted however that there are many paths, other than trajectories of
$\calh$, that tend to produce contraction for the corresponding holonomy maps of
$\fol$. These include, for example, the trajectories of $\calh^{\theta}$, $-\pi/2 < \theta < \pi/2$, or suitable combinations
of those. Therefore there is a large amount of flexibility to choose ``deformed trajectories'' when a trajectory
of $\calh$ approaches a singularity such as a saddle point.

Before giving precise definitions of what is meant by ``deformed trajectory'' or
by ``trajectory of finite length'', we shall perform a few reductions in our setting
so as to make the subsequent discussion more transparent. Let then $\omega, \, \fol$ be
as above. Denote by $(\omega)_0^{\fol}$ the sub-divisor consisting of those irreducible components of
$(\omega)_0$ that are invariant under $\fol$. Similarly set $(\omega)_0^{\perp \fol} = (\omega)_0
\setminus (\omega)_0^{\fol}$. Denoting by $(\omega)_{\infty}$ the
divisor of poles of $\omega$, the subdivisors $(\omega)_{\infty}^{\fol}$
and $(\omega)_{\infty}^{\perp \fol}$ are analogously defined and so are the divisors $(\omega_1)_0^{\fol}, \, (\omega_1)_0^{\perp \fol}$.
Let us remind the reader that $(\omega_1)_{\infty}$ is contained in $(\omega)_0^{\perp \fol} \cup (\omega)_{\infty}^{\perp \fol}$
so that it has no component invariant by $\fol$, cf. Lemma~\ref{newversionSection2.11}.

The next lemma allows us to assume some standard ``normalization'' conditions.

\begin{lema}
\label{revision2}
Modulo performing finitely many blow-ups, the conditions below are always
satisfied:
\begin{enumerate}

\item The singular set ${\rm Sing}\, (\fol)$ of $\fol$ is disjoint from
$(\omega)_0^{\perp \fol} \cup (\omega)_{\infty}^{\perp \fol}$ as well as
from $(\omega_1)_0^{\perp \fol}$.

\item Every irreducible component of $(\omega)_0, \, (\omega)_{\infty}$ and
of $(\omega_1)_0$ is smooth.

\item The divisor of zeros $(\omega)_0$ does not intersect the divisor of poles
$(\omega)_{\infty}$ at regular points of $\fol$.

\item Two distinct irreducible components of $(\omega)_0^{\perp \fol}$
(resp. $(\omega)_{\infty}^{\perp \fol}$, $(\omega_1)_0^{\perp \fol}$) are disjoint.

\item $\fol$ is transverse to every irreducible component of
$(\omega)_0^{\perp \fol} \cup (\omega)_{\infty}^{\perp \fol}$ or of
$(\omega_1)_0^{\perp \fol}$.
\end{enumerate}
\end{lema}

\noindent {\it Proof}. It is clear that the singularities of $\fol$ can be supposed to be
reduced (Seidenberg's theorem). Similarly the irreducible components of
$(\omega)_0, \, (\omega)_{\infty}$ and of $(\omega_1)_0$ can easily be made smooth.

To show that ${\rm Sing}\, (\fol)$ can be made disjoint from
$(\omega)_0^{\perp \fol} \cup (\omega)_{\infty}^{\perp \fol}$, let $\mathcal{C}$
be a local branch of
an irreducible component of $(\omega)_0^{\perp \fol} \cup (\omega)_{\infty}^{\perp \fol}$
passing through $p \in {\rm Sing}\, (\fol)$. By assumption $\mathcal{C}$ is not
invariant by $\fol$ so that it has a contact of finite order with the actual separatrices
of $\fol$ at $p$. By blowing-up $\fol$ at $p$, the new singularities appearing in
the exceptional divisor $\pi^{-1} (p)$ have their positions determined by the tangent
spaces at $p$ to the local separatrices of $\fol$. Therefore, after finitely many
repetitions of this procedure, the proper transform of $\mathcal{C}$ will no longer
pass through any of the resulting singularities of the blown-up foliation. A similar
argument applies to the divisor $(\omega_1)_0^{\perp \fol}$. Note also that, in the course of
performing the mentioned blow-ups, the ``new components'' of $(\omega)_0, \,
(\omega)_{\infty}$ and of $(\omega_1)_0$ that may have been introduced are all contained in the exceptional
divisor. Hence they are invariant by the corresponding foliation, i.e. they are not
contained in $(\omega)_0^{\perp \fol} \cup (\omega)_{\infty}^{\perp \fol} \cup (\omega_1)_0^{\perp \fol}$.

The remaining ``reductions'' are based on the following remark: if a regular
point of a foliation is blown-up, then the new foliation still leaves the
exceptional divisor invariant. Furthermore this exceptional divisor contains a
unique singularity of the blown-up foliation. This singularity is conjugate to
the linear singularity with eigenvalues $1,-1$.

Consider now a point $p \in M$ regular for $\fol$ where $(\omega)_{\infty}$
intersects $(\omega)_0$ and let $\fol$ be blown-up at $p$. As before, after finitely
many blow-ups, the proper transforms of $(\omega)_{\infty}$ and $(\omega)_0$ will
be separated. The components added by these blow-ups are all contained in the
exceptional divisor and thus are invariant by the corresponding foliation. In particular,
if we just wanted to ensure that $(\omega)_{\infty}^{\perp \fol}$ does not
intersect $(\omega)_0^{\perp \fol}$ at a regular point this would be enough. For
the general case, it suffices to note that the order of the exceptional
divisor resulting from a single blow-up is the difference between the orders
of the components of $(\omega)_{\infty}$ and of $(\omega)_0$ that pass through
the center of the blow-up. Thus after finitely many repetitions, there will appear
a exceptional divisor which is {\it regular for $\omega$}\, in the sense that it is
not contained in either $(\omega)_{\infty}$ or $(\omega)_0$. This leads to the verification
of item~3. The same reasoning allows us to obtain item~4 as well.

Finally, as to item~5, let $D$ be an irreducible component of
$(\omega)_0^{\perp \fol} \cup (\omega)_{\infty}^{\perp \fol}$ or of
$(\omega_1)_0^{\perp \fol}$. We need to check that $\fol$ can be made transverse to $D$.
Thanks to the preceding items, we can assume that $D \cap {\rm Sing}\, (\fol) =\emptyset$.
Next observe that the number of tangencies
between $\fol$ and $D$ is finite since$D$ is not invariant by $\fol$ (and of course every tangency has finite contact).
Thus once again we only need to blow-up tangency points sufficiently many times.
As always the exceptional divisors added in the procedure are all invariant by the
foliation and thus do not destroy the previous ``reductions''. This completes
the proof of the lemma.\qed

The behavior of $\calh$ near points in $(\omega_1)_{\infty} =
(\omega)_0^{\perp \fol} \cup (\omega)_{\infty}^{\perp \fol}$ is clear. However
the behavior of $\calh$ near points in $(\omega_1)_0$ needs further comments and,
in particular, leads to the notion of ``deformed trajectory''. It is convenient to identify three {\it critical regions}\,
where the trajectories of $\calh$ will be allowed to be deformed. These are as follows.

\noindent {\bf First critical region}: The divisor $(\omega_1)_0^{\perp \fol}$.

\noindent Let $C_j^0$, $j=1,\ldots ,l$, denote the irreducible components
of the zero divisor $(\omega_1)_0^{\perp \fol}$ of $\omega_1$. Recalling that every $C_j^0$
is smooth and transverse to $\fol$, we can find a small ``tubular neighborhood''
$\Vv_j$ of $C_j^0$ whose boundary $\partial \Vv_j$ is still transverse to $\fol$, for
every $j=1,\ldots ,l$. Besides, for $j$ fixed, we also assume that the intersection
of $\Vv_j$ with $(\omega_1)_0$ is reduced to $C_j^0$ (cf. item~4 of  Lemma~\ref{revision2}).
Finally set $\Vv = \bigcup_{j=1}^l \Vv_j$. If $p \in \partial \Vv$, we can suppose without
loss of generality that the leaf $L_p$ of $\fol$ through $p$ ``locally slices'' $\Vv$
into a connected disc. The collection of these ``discs'' form the fibers of a
differentiable submersion $\Vv \rightarrow (\omega_1)_0$.

Next let $p \in \partial \Vv$ and let $D_p \subset L_p$ be the above mentioned
disc, i.e. $D_p$ is the connected component containing $p$ of $L_p \cap \Vv$. The
structure of the trajectories of $\calh$ on $D_p$ is described by Lemma~\ref{blm2}.
Denoting by $\calh_{\vert D_p}$ the restriction of $\calh$ to $D_p$,
it follows the existence of $2m$ separatrices for $\calh_{\vert D_p}$
at $\textsc{Q} \simeq D_p \cap (\omega_1)_0$, $m \geq 2$. These separatrices are divided into
two groups. Namely there are $m$ separatrices over which one converges to $\textsc{Q}$
by moving in the sense of their orientation (these separatrices are said {\it to
approach $\textsc{Q}$}). The remaining $m$ separatrices are such that one converges to
$\textsc{Q}$ by moving in the sense opposite to their orientation (these separatrices are said {\it to
leave $\textsc{Q}$}). In addition, the total picture is symmetric by a rotation group of order $m$.

In this situation, if $l$ is for example a separatrix approaching $\textsc{Q}$, we
allow $l$ to be continued in the ``future'' by following one of
the separatrices of $\textsc{Q}$ that
leaves $\textsc{Q}$. The chosen separatrix can for example be one of the two separatrices
that are ``closest'' to $l$ but this is not necessary. More generally every trajectory of $\calh$ entering the neighborhood
$\Vv$ is allowed to be continued by ``following'' any of the separatrices of $\calh$ leaving $\textsc{Q}$.

\noindent {\bf Second critical region}: Siegel singularities of $\fol$ such that the restriction of $\omega_1$ to its local separatrices is holomorphic.

Fixed a Siegel singularity as above, let $V_0$ be a neighborhood of it similar to
the one defined before the statement of Lemma~\ref{3lema3}. In particular
we have sets $A = \{ (u,v) \in \C^2 \; ; \; \vert u \vert =\epsilon_1 \; \; {\rm and} \; \; \vert v \vert < \epsilon_2 \}$ and
$B = \{ (u,v) \in \C^2 \; ; \; \vert v \vert =\epsilon_2' \; \; {\rm and} \; \; \vert u \vert < \epsilon_1' \}$. To fix notations let $l$ be an oriented trajectory
of $\calh$ entering this neighborhood at an intersection point of $l$ and $A$ (i.e. $l$ is ``close'' to the axis $\{ v=0\}$). Denote by $L$ the leaf of
$\fol$ containing $l$. Since the restriction of $\omega_1$ to $\{v=0\}$ is holomorphic at the origin, the trajectory $l$ can be deformed {\it inside}\,
$L$ to avoid crossing the set $A$ (i.e. this trajectory can be deformed so as to stay away from the singularity itself). This deformation is similar
to the deformations performed in the case of saddle singularities of $\calh$ that appear in connection with the divisor $(\omega_1)_0^{\perp \fol}$.
In particular the continuation of $l$ will stay ``close to $\{v=0\}$'' during the procedure. In fact, this trajectory will leave the singularity by ``following''
one of the separatrices of $\calh_{\vert \{v=0\}}$ that leave the singularity in question (where $\calh_{\vert \{v=0\}}$ stands for the restriction of
$\calh$ to $\{v=0\}$). Another possibility to defined continuations for $l$ is to let $l$ enter the neighborhood of the mentioned singularity and then use
the corresponding Dulac transform to continue $l$ as a trajectory of $\calh$ that is now ``close to $\{u=0\}$''. In this case the desired continuation
of $l$ will be ``close'' to one of the separatrices of $\calh_{\vert \{u=0\}}$ oriented so as to leave the mentioned singularity
(where $\calh_{\vert \{u=0\}}$ stands for the restriction of $\calh$ to $\{u=0\}$). Summarizing it can be said that a trajectory $l$ of $\calh$ intersecting
the set $A$ admits all the above mentioned continuations.

\begin{obs}
\label{localregions}
{\rm In the preceding two types of critical regions the ``deformation'' of the trajectory $l$ consists of adding to it a ``small'' segment
of trajectory of $\calh^{\perp}$. By construction these pieces of $\calh^{\perp}$-trajectories have length bounded by a ``small constant''
and besides they are strictly comprised between two ``genuine'' segments of $\calh$-trajectories whose lengths are bounded from below
by positive constants depending solely on $\fol, \, M$. As a consequence these ``deformations'' do not disrupt the global contractive
nature of holonomy maps of $\fol$ defined by means of ``deformed trajectories of $\calh$''. We shall return to this point below.}
\end{obs}

\begin{obs}
\label{localregionsandirrationalfoci}
{\rm Besides singularities belonging to the Siegel domain also irrational focus singularities may be considered. Recall that an irrational focus singularity is
linearizable and hence it possesses exactly two separatrices. These separatrices are smooth and may be chosen as the coordinate axes
in the linearizing coordinates. The fact that the quotient between the eigenvalues
of these singularities cannot be a rational number implies that the only way in which the restriction of $\omega_1$ to these separatrices
may be holomorphic occurs when both separatrices are components with multiplicity~$1$ of $(\omega)_{\infty}$. This case will rarely occurs,
but if it does, the trajectories of $\calh$ will be deformed so as to avoid the singularity in the same way it may be done for analogous Siegel
singularities. Since irrational foci have no associated Dulac transforms only this type of continuation will be allowed in the present case.}
\end{obs}

\noindent {\bf Third critical region}: The divisor $(\omega_1)_0^{\fol}$.

By construction the support of the divisor $(\omega_1)_0^{\fol}$ consists of (irreducible) curves invariant by $\fol$. Let $C$ denote
one of these curves. Then the restriction of $\omega_1$ to $C$ vanishes identically so that it does not define any real foliation on
$C$. Nonetheless Poincar\'e Lemma can still be applied to this situation. In fact, let $c:[0,1] \rightarrow C$ be a path contained in $C$
and consider local transverse sections $\Sigma_{c(0)}, \, \Sigma_{c(1)}$ through $c(0), \, c(1)$ respectively. If the sections
$\Sigma_{c(0)}, \, \Sigma_{c(1)}$ are parameterized as indicated in Section~2.3, then the holonomy map ${\rm Hol}\, (c) :
\Sigma_{c(0)} \rightarrow \Sigma_{c(1)}$ obtained from $c$ and $\fol$ is such that $[ {\rm Hol}\, (c) ]'(0) =1$. In particular the usual
holonomy group associated to the ``leaf'' $C$ with respect to $\fol$ is entirely constituted by local diffeomorphisms tangent to the
identity.

Since the foliation $\calh$ is not defined on $C$, we shall allow every (``minimizing geodesic'') path joining two points of $C$ with length
less than the diameter of $C$ to be used to continue a given trajectory $l$ of $\calh$. Here both ``length'' of the path and ``diameter'' of $C$
arise from fixing once and for all some auxiliary Hermitian metric on $M$.

To better explain the above definition, consider two Siegel singularities $\textsc{P}, \textsc{Q}$ of $\fol$ lying in $C$. Denote by $S_{p}$
(resp. $S_{q}$) the local separatrix of $\fol$ transverse to $C$ at $\textsc{P}$ (resp. $\textsc{Q}$). Also fix neighborhoods $V_{p},
\, V_{q}$ of $\textsc{P}, \textsc{Q}$ as in the case of discussed in the second critical region. Since $\omega_1$ vanishes identically
on $C$, it follows that the restriction of $\omega_1$ to $S_{p}$ (resp. $S_{q}$) is holomorphic. If $l \subset S_{p}$
is a trajectory of $\calh$ that enters $V_{p}$, then $l$ can be continued as a trajectory $l'$ of $\calh$ {\it contained in
$S_{q}$ and oriented so as to leave the neighborhood $V_{q}$}. A similar convention applies to trajectories $l$ of $\calh$ that are
not contained in $S_{p}$ but that still enters the neighborhood $V_{p}$ of $\textsc{P}$. A continuation $l'$ for the trajectory $l$ will
be such that $l, l'$ are contained in the same global leaf $L$ of $\fol$ and $l'$ leaves the neighborhood $V_{q}$ of $\textsc{Q}$. In other words
the continuation of these trajectories can be pictured as if the curve $C$ were collapsed into a single point (heuristically imagined as a Siegel singularity
whose separatrices would be $S_{p}, \, S_{q}$). Then the mentioned continuation would be defined as in the case of the second
critical region discussed above.

\begin{obs}
\label{localregionsandirrationalfociPRIME}
{\rm In line with Remark~\ref{localregionsandirrationalfoci}, the use of Dulac transforms to follow a component $C$ of $(\omega_1)_0^{\fol}$
as above is only possible at a Siegel singularity i.e. no irrational focus singularity lying in $C$ will be associated with continuation of trajectories
by means of Dulac transforms.}
\end{obs}

We are now ready to define {\it global deformed trajectories of $\calh$}. Away from a fixed neighborhood of the three critical regions previously
discussed, a deformed trajectory must agree with an ordinary trajectory of $\calh$. However if a trajectory $l$ of $\calh$ enters a critical region,
then it possesses all the corresponding continuations mentioned above. As a consequence, every possible continuation
of $l$ will eventually leave the critical region
in question and become again an ordinary trajectory of $\calh$. In particular, given $p \in M$, the deformed trajectory of $\calh$ through $p$ is
in general not uniquely determined. It is convenient to think of it not as a ``single path'' but as a collection of paths that ramifies whenever one of its
branches enters a critical region. In other words, the $\calh$-trajectory through $p$ is in general not a single path but rather a collection of paths (or branches)
that are allowed to ramify at the critical regions.

Similar definitions apply if we decide to follow a (deformed) trajectory $l$ of $\calh$ in the direction opposite to its orientation (i.e. when the
``past'' of $l$ is considered). More generally for $\theta \in (-\pi/2 ,\pi/2)$ fixed, the deformed trajectories of $\calh^{\theta}$ are analogously defined
and the same remark concerning orientation can be done to define their ``continuations in the past''.

To fully define what will be understood by a deformed trajectory of $\calh$ (or $\calh^{\theta}$ we still need to clarify what to do when an ordinary
trajectory becomes close to the remaining ``singularities'' of $\calh$ (or of $\calh^{\theta}$). However before doing this, it is important to point out that
deformed trajectories as considered above are such that the corresponding holonomy maps of $\fol$ still keep the contractive behavior
characteristic of ordinary trajectories of $\calh$ (or $\calh^{\theta}$ for $\theta \in (-\pi/2 ,\pi/2)$). As in Section~2.3, recall that we have fixed an auxiliary
Hermitian metric on $M$ so that it is possible to consider the length of paths contained in $M$. For those paths whose images are contained in leaves
of $\fol$, their resulting lengths are also comparable with the sum of the lengths of their representatives in a fixed foliated atlas of $M$, cf. Section~2.3. Next
let $\theta \in (-\pi/2 ,\pi/2)$ be fixed and
consider a path $c:[0,1] \rightarrow L \subset M$ parameterizing a segment of deformed trajectory of $\calh$ (resp.
$\calh^{\theta}$ for fixed $\theta \in (-\pi/2 ,\pi/2)$) as above. Thus $c$ can be viewed as a concatenation
of paths $c^i$ that either are contained in a critical region or are segments of (ordinary) $\calh$-trajectories (resp. $\calh^{\theta}$-trajectories) away
from the critical regions and from the remaining singularities of $\calh, \, \calh^{\theta}$. In the latter case, the length of $c^i$ is bounded from below by a positive constant.
On the other hand a path $c^i$ whose image is contained in a critical region is such that $c^{i-1}$ and $c^{i+1}$ parameterize an ordinary segment of $\calh, \calh^{\theta}$
(lying in a compact part of the complement of the critical regions and of the remaining singularities of $\calh, \, \calh^{\theta}$). Furthermore these
paths $c^i$ are such that their length is uniformly bounded and, besides,
the holonomy map of $\fol$ obtained by means of $c^i$ (with respect to suitable transverse sections parameterized as indicated in Section~2.3)
is a holomorphic diffeomorphism whose linear part has modulo equal to~$1$. Then the preceding discussion can be summarized by Proposition~\ref{betterthanThmblm}
below, which ensures the exponential decay of the norm of the derivative at $c(0)$ with the length of $c$.

\begin{prop}
\label{betterthanThmblm}
Consider a path $c: [0,1] \rightarrow L$ that parametrizes a segment of deformed trajectory of $\calh$ or, more generally,
of $\calh^{\theta}$ ($-\pi/2 < \theta < \pi/2$). Then there are constants $C, k$ depending solely on $\theta$ (for $M, \fol, \omega$ and
the auxiliary Hermitian metric fixed) such that the estimate below holds
\begin{equation}
\vert ({\rm Hol}\, (c))' (0) \vert \leq  C  \exp \, (-k \, {\rm length}\,
(c) /2) \; , \label{4eq1}
\end{equation}
where ${\rm Hol}\, (c)$ stands for the holonomy map of $\fol$ induced by $c$.\qed
\end{prop}

Let us now complete the definition of deformed trajectories for $\calh$. First note that the ``remaining singularities'' of $\calh, \calh^{\theta}$
are provided by either singular points of $\fol$ (different from Siegel singularities since these were already taken into account) and
by the divisors $(\omega)_0^{\perp \fol}$ and
$(\omega)_{\infty}^{\perp \fol}$. For the purposes of this paper however, Proposition~\ref{4.5prop1} allows us to rule out hyperbolic
singularities as well as saddle-node singularities from the discussion below.

Unless otherwise stated, in what follows we shall simply say {\it trajectory of
$\calh$ (resp. $\calh^{\theta}$)}\, instead of ``deformed trajectory of $\calh$ (resp. $\calh^{\theta}$)''. Hence for $p \in \calk$,
let $l_p$ denote the trajectory of $\calh$ through $p$ (in precise words, this
means a deformed trajectory of $\calh$ through $p$) and consider the leaf $L$ of $\fol$ containing $l_p$.
Let us first introduce the notion of {\it endpoint}\, for $l_p$.
The trajectory $l_p$ is said to have an {\it endpoint}\,
at a point $q \in M$ if one of the following possibilities hold:
\begin{itemize}
\item $q \in (\omega)_0^{\perp \fol}$ is a sink of $\calh_{\vert L}$ and $\overline{l}_p^+ = q$.

\item $q \in (\omega)_{\infty}^{\perp \fol}$ is a source of $\calh_{\vert L}$ and $\overline{l}_p^- = q$.

\item $q$ is a sink-irrational focus (resp. source-irrational focus) singularity of $\fol$
to which $l_p$ converges (resp. from which $l_p$ is emanated, cf. Lemma~\ref{4.5lema3} and Remark~\ref{4.5obs1}).

\end{itemize}
Similar definitions apply to the case of $\calh^{\theta}$-trajectories, $\theta \in (-\pi/2 ,\pi/2)$.
Sometimes we shall use the expressions {\it future end}\, (resp. {\it past end}) to refer to the cases above
which are concerned with a sink-like (resp. source-like) endpoint of $l_p$. To define the trajectory of $l_p$ of
$\calh$ through $p$ we start with the ordinary trajectory of $\calh$ through $p$. Whenever this trajectory enters one of
the above described critical regions, all its resulting ramifications are considered together as its continuations. Thus
it is possibly more convenient to
speak about {\it branches} of $l_p$. In this case a branch of $l_p$ has a future endpoint if its converge to a point $q$ of $M$ that behaves
locally as a sink for $\calh$. In view of the preceding $q$ either belongs to $(\omega)_0^{\perp \fol}$ or it is a sink-irrational focus singularity
of $\fol$. Naturally we can also follow the trajectory $l_p$ in the sense opposite to its standard orientation. In this case, we shall denote the
resulting semi-trajectory by $l_p^-$ (in certain cases where the context might be unclear, the semi-trajectories through $p$ with the usual orientation
will also be denoted by $l_p^+$). Past endpoints for a branch of $l_p^-$ is then analogously defined and so are future and past endpoints for
branches of the trajectories $l_p^{\theta} = l_p^{\theta, +}$ and $l_p^{\theta, -}$ of $\calh^{\theta}$ through $p$.

The length of a branch of the semi-trajectory $l_p^+$ is defined in natural differential geometric
terms for the auxiliary Hermitian metric fixed from the beginning provided that the branch is finite. Otherwise the branch is said
to be of infinite length.
Now the semi-trajectory $l_p^+$ of $\calh$ (resp. $\calh^{\theta}$) through $p$ is said to be {\it finite}\, if and only the supremum of the
length of all its branches is finite. In this case the number of branches of $l_p^+$ is itself finite so that the supremum is also attained.
Once again the definition of length for the semi-trajectory $l_p^-$ can analogously be given. Finally the deformed trajectory $l_p$ of
$\calh$ (resp. $\calh^{\theta}$) through $p$ will be called finite if both semi-trajectories $l_p^+$ and $l_p^-$ are so. The length of $l_p$ will
then be the maximum of the lengths of all branches contained in $l_p$.

\begin{obs}
{\rm If a branch of a $\calh$-trajectory $l_p$ (resp. $\calh^{\theta}$-trajectory $l_p^{\theta}$) consists of a loop, possibly passing through
critical regions, then this branch contains neither future nor past endpoints. It then follows that its length is infinite.
This means that one is allowed to go around the loop
infinitely many times what explains why the length of the loop must be considered as infinite. This is very
natural as definition since the holonomy of $\fol$ associated to one of these trajectories is clearly hyperbolic.

In fact, with our
terminology, the length of $l_p$ is finite if and only if all branches of $l_p$ possess both future and
past ends and, in addition, the supremum of the lengths of these branches is finite. It is also
clear that the above definition is invariant by blow-ups/blow-downs. Therefore
the length of the trajectories of $\calh$ (resp. $\calh^{\theta}$) can be considered whether or not the foliation
$\fol$ has reduced singularities. More generally, this definition makes sense whether
or not the normalizing conditions of Lemma~\ref{revision2} are satisfied.}
\end{obs}

With the above terminology, the contents of Proposition~\ref{betterthanThmblm} can be complemented
by the following
simple generalization of Theorem~\ref{blm} that is better adapted to our needs.
Let $K$ be a compact part of the complement of the singular set of $\fol$ and consider a path $c: [0,1]\rightarrow K$
parameterizing a segment of deformed $\calh$-trajectory  (resp. $\calh^{\theta}$-trajectory, $\theta \in (-\pi/2 ,\pi/2)$).
Finally let ${\rm Hol}\, (c)$ denote the holonomy map of $\fol$ induced by $c$ and recall that these maps are identified
with local diffeomorphisms of $\C$ by means of transverse sections $\Sigma_{c(0)}, \, \Sigma_{c(1)}$ parameterized by
$\omega$ (cf. Section~2.3). Then we have:

\begin{teo}
\label{Totalversionofblm}
With the preceding notations there is $\delta >0$ (depending only on $K$) and constants $C, k > 0$ such that
the following holds:
\begin{enumerate}
\item ${\rm Hol}\, (c)$ is defined on the transverse disc $B_{c(0)} (\delta )$ of radius $\delta > 0$ about $c(0)$.

\item The image $( {\rm Hol}\, (c)) (B_{c(0)} (\delta ))$ of $B_{c(0)} (\delta )$ by ${\rm Hol}\, (c)$ is contained in a
transverse disc $B_{c(1)} (r)$ of radius $r$ about $c(1)$ where
$$
r \leq C  \exp \, (-k \, {\rm length}\, (c) /2) \, .
$$
\end{enumerate}
More generally if $c$ parameterizes a segment of deformed trajectory of $\calh^{\theta}$, $\theta \in (-\pi/2 ,\pi/2)$,
then the statement still holds, only the values of the constants $\delta, C, k$ will depend further on $\theta$.
\end{teo}

\noindent {\it Proof}\,: Fix a finite covering of $K$ by foliated coordinates of $\fol$ along with transverse sections parameterized
by $\omega$ as indicated in Section~2.3. According to Proposition~\ref{betterthanThmblm}, there are constants $C_1, k_1$ such that
the absolute value of the derivative of ${\rm Hol}\, (c)$ at $c(0)$ satisfies the estimate
$$
\vert ({\rm Hol}\, (c))' (0) \vert \leq  C_1  \exp \, (-k_1 \, {\rm length}\, (c) /2)  \, .
$$
As in \cite{blm}, this estimate allows us to show that ${\rm Hol}\, (c)$ is defined on a uniform domain $B_{c(0)} (\delta )$. To check
the rest of the statement, note that ${\rm Hol}\, (c)$ is univalent on its domain of definition. Thus modulo reducing this domain,
K\"oebe's theorem (cf. \cite{unival}) can be applied to ensure that ${\rm Hol}\, (c)$ has ``bounded distortion'' so that the diameter
of its image can be estimate from the value of its derivative at $c(0)$. This completes the proof of the theorem for $\theta=0$. The general
case is however totally analogous.\qed

From now to the rest of the paper we fix a closed set $\calk \subseteq \supT$ that is minimal
for $\fol$. Denote by $\calh_{\calk}$ (resp. $\calh^{\theta}_{\calk}$)
the restriction of $\calh$ (resp. $\calh^{\theta}$) to $\calk$. As it is usually the case, in the
sequel the word ``trajectory'' actually means ``deformed trajectory''.
Let us close this section with the following proposition:

\begin{prop}
\label{4prop1}
The following alternative holds:
\begin{itemize}
\item There is a uniform constant $C$ (resp. $C^{\theta}$) such that the length of every trajectory of
$\calh_{\calk}$ (resp. $\calh^{\theta}_{\calk}$) is less than $C$ (resp. $C^{\theta}$).

\item There is a non-empty compact set $\calk^0 \subseteq \calk$, invariant by $\calh$ (resp. $\calh^{\theta}$),
where all the corresponding trajectories of (the restriction of) $\calh$ (resp. $\calh^{\theta}$) have infinite
length.
\end{itemize}
\end{prop}

\begin{obs}
\label{calhinvariant}
{\rm The invariance of $\calk^0$ by $\calh$ (resp. $\calh^{\theta}$) means that through each point of $\calk^0$ there
passes {\it a branch of trajectory of $\calh$ (resp. $\calh^{\theta}$)}\, which is entirely contained in $\calk^0$. Since, in general,
the trajectory of $\calh$ (resp. $\calh^{\theta}$) through $p$ is constituted by several branches, it may well happen that some of them
are not fully contained in $\calk^0$.}
\end{obs}

\vspace{0.2cm}

\noindent {\it Proof of Proposition~(\ref{4prop1})}\,: In the sequel we suppose that the
conditions of Lemma~\ref{revision2} are satisfied. It suffices to check the statement for the foliation
$\calh=\calh^0$.

Denote by ${\rm Irr}_+ (\fol)$ (resp. ${\rm Irr}_- (\fol)$) the irrational focus
singularities
of $\fol$ in $\calk$ that behave as a sink (resp. source) for $\calh$ in the sense
of Lemma~(\ref{4.5lema3}). Given one such singularity $p$, let $B_p (\epsilon)$ be
the real $3$-dimensional ball of radius $\epsilon>0$ about $p$. It is easy
to see that $\fol$ is transverse to this ball for $\epsilon$ sufficiently small. In fact,
also the foliation $\calh$ is transverse to this ball as it easily follows from
Proposition~(\ref{4.5lema3}). Let then ${\rm Irr}_+^{\epsilon} (\fol)$ be the union
of these balls about the points in ${\rm Irr}_+ (\fol)$. The set
${\rm Irr}_-^{\epsilon} (\fol)$ is analogously defined.

Suppose that all the (deformed) $\calh$-trajectories
contained in $\calk$ are of finite length. We are going to show the existence
of a uniform bound for all the corresponding lengths. If this bound did not exist, then there
would be a sequence $\{ l_i\}_{i \in \N}$ of branches of $\calh$-trajectories contained in $\calk$
such that the sequence formed by their corresponding lengths goes off to infinity. For each
$i$, let $c_i (a_i , b_i) \subset \R \rightarrow M$ be a parametrization of $l_i$.
Naturally
$c_i (a_i)$ belongs to $(\omega)_{\infty} \cup {\rm Irr}_- (\fol)$
whereas $c_i (b_i)$ belongs to
$(\omega)_0 \cup {\rm Irr}_+ (\fol)$. Thus we conclude that,
in fact, $c_i (a_i) \in (\omega)_{\infty}^{\perp \fol} \cup {\rm Irr}_- (\fol)$
and $c_i (b_i) \in
(\omega)_0^{\perp \fol} \cup {\rm Irr}_+ (\fol)$. Modulo passing to a subsequence, we can
suppose that the $c_i (a_i)$ (resp. $c_i (b_i)$) converge to a point $a \in \calk$
(resp. $b \in \calk$). One of the following two possibilities
must occur:

\noindent 1. $a$ (resp. $b$) belongs to $D_{\infty} \subseteq
(\omega)_{\infty}$ (resp. $b \in D_{0} \subseteq (\omega)_0$) where $D_{\infty}$
(resp. $D_0$) stands for an irreducible component of $(\omega)_{\infty}$
(resp. $(\omega)_0$).

\noindent 2. $a$  (resp. $b$) belongs to ${\rm Irr}_- (\fol)$ (resp.
${\rm Irr}_+ (\fol)$). In this case, modulo shortening the length of $c_i$ by
a uniform small constant (cf. Proposition~\ref{4.5lema3}) we can replace $a_i$ (resp. $b_i$)
by $a_i'$ (resp. $b_i'$) such that $c_i (a_i') \in {\rm Irr}_-^{\epsilon} (\fol)$
(resp. $c_i (b_i) \in {\rm Irr}_+^{\epsilon} (\fol)$). Therefore we can consider
without loss of generality that $a \in {\rm Irr}_-^{\epsilon} (\fol) \cup
{\rm Irr}_+^{\epsilon} (\fol)$.

Consider the leaf $L_a$ of $\fol$ through $a$ and note that the restriction of $\calh$
to $L_a$ is well-defined (it is not fully constituted be a critical region of third type). In particular, there is a trajectory
$l_a^+$ of $\calh$ being emanated from $a$. Although this trajectory possibly consists of several branches, due to ramification
at critical regions, it contains one special branch defined as follows: whenever the (branch of the) trajectory in question
enters a critical region, its continuation is dictated by the continuations of
the $l_i$'s (for $i$ large enough). The resulting branch $l_a^+$ clearly has infinite length.
Otherwise leaves emanated from points sufficiently close to $a$, and choosing appropriate
ramifications at critical regions, would have bounded
length and this would contradict our assumption.

The preceding discussion also shows the existence of semi-trajectories of infinite lengths provided that the first case
in the statement of the proposition does not occur. Let then $l^+$ denote a branch of infinite length contained in some deformed
trajectory in $\calk$. Since $\calk$ is compact, the closure $\overline{l}^+ \subset \calk$ of $l^+$ is not empty. Besides every
semi-trajectory through a point of $\overline{l}^+$ is infinite, or in more accurate terms, it contains a branch of infinite length. In fact,
if all branches of a deformed trajectory through a point $p \in \overline{l}^+ \subset \calk$ were of finite length, then the above
argument would imply that $\overline{l}^+$ intersects $(\omega)_0^{\perp \fol}
\cup {\rm Irr}_+^{\epsilon} (\fol)$. This is however impossible since it contradicts
the infinite length of $l^+$. In other words, $\calk^0 = \overline{l}^+$ satisfies the condition in
the second alternative of our statement. The proposition is proved.\qed

Since $\calk^0=\overline{l}^+$ and $(\omega)_0^{\perp \fol}
\cup {\rm Irr}_+^{\epsilon} (\fol)$ are compact disjoint,
there is a positive distance between them. Also, by construction, $l^+$ cannot
accumulate (in the future) on $(\omega)_{\infty}^{\perp \fol}
\cup {\rm Irr}_-^{\epsilon} (\fol)$ thanks to Lemma~(\ref{blm1}) and
Lemma~(\ref{4.5lema3}).
Thus we obtain:

\begin{coro}
\label{4coro1}
Suppose that the first alternative in Proposition~(\ref{4prop1}) is not
verified. Then there is a small open neighborhood $V$ of $(\omega)_0^{\perp \fol} \cup
(\omega)_{\infty}^{\perp \fol} \cup {\rm Irr}_-^{\epsilon} (\fol) \cup
{\rm Irr}_+^{\epsilon} (\fol)$
such that $\calk^0 \cap V =\emptyset$. In particular all the singularities of
$\fol$ lying in $\calk^0$ are in the Siegel domain unless are irrational foci as in Remark~\ref{localregionsandirrationalfoci}.
\end{coro}

\section{Invariant currents vs. infinite trajectories of $\calh$}

The remaining two sections are devoted to proving the theorems stated in the Introduction.
We keep the context and the notations of Section~5. Recalling that
$\calk$ stands for a minimal set of $\fol$ contained in the support of $T$, the restriction of
$\calh$ to $\calk$ is going to be denoted by $\calh_{\calk}$. Let us begin by rephrasing Theorem~A:

\begin{teo}
\label{fim1}
Let $\fol$ and $T$ be as above. If $\calk$ does not contain a compact leaf of $\fol$,
then all deformed $\calh$-trajectories (resp. $\calh^{\theta}$-trajectories with fixed $\theta \in (-\pi/2, \pi/2)$) in $\calk$ have length smaller than some
positive constant ${\rm Const}$.
\end{teo}

It suffices to prove the statement for $\calh =\calh^0$ since the generalization to $\calh^{\theta}$
is very straightforward. Thus let us suppose that
the lengths of the $\calh$-trajectories in $\calk$ are not uniformly bounded. Our aim will then be to ensure the existence of an
algebraic curve contained in $\calk$.
Since the lengths of the $\calh$-trajectories in $\calk$ are not uniformly bounded,
we can consider a compact set $\calk^0 \subset \calk$ satisfying the
conclusions of Proposition~\ref{4prop1} and Corollary~\ref{4coro1}. Moreover, by applying Zorn Lemma,
we can assume without loss of generality that $\calk^0$ is {\it minimal}\, for $\calh$
i.e.  through every point in $\calk^0$ there
passes a branch of (deformed) trajectory of $\calh$ which is dense in $\calk^0$ (the branch being obviously contained in $\calk^0$).
Here it is worth pointing out that the assumption that $\calk^0$ is minimal is not indispensable for our discussion (and so the use of Zorn Lemma
can also be avoided). In fact, it would be enough
to consider an accumulation point of a $\calh$-trajectory $l^+ \subset \calk$ of infinite length that happens to be regular for $\calh$. Because of
Theorem~\ref{Totalversionofblm} the holonomy maps of $\fol$ induced by (segments of) $l^+$ is defined on a uniform domain. Thus if $l^+$ accumulates
on a regular point of $\calh$ (and of $\fol$) this trajectory will be captured by the holonomy maps to which it gives rise (modulo a slight local deformation of
$l^+$). The latter statement would be sufficient for our purposes.
Yet it is simpler to assume that $\calk^0$ is minimal so
that a self-accumulating $\calh$-trajectory $l^+$ can be selected.

Next let us perform on $\fol$ the normalizations described in Lemma~\ref{revision2}. Since these transformations include the blowing-up of points,
they give rise to {\it compact leaves}\, contained in the closure of the transform of $\calk$. The curves obtained in this way however
{\it do not form loops}\, and this will be exploited in the sequel, cf. below. More generally let $A_{\calk}$ denote the union of all algebraic curves contained in
transform of $\calk$. Actually, by an abuse of notation, the closure of the proper transform of $\calk$ will still be denoted by $\calk$.
In this sense, $\calk$ is no longer minimal for $\fol$ but it satisfies the following condition: every leaf of $\fol$ in $\calk$ that is not dense in $\calk
\setminus A_{\calk}$ is necessarily contained in $A_{\calk}$ itself. This condition is going to be used in the sequel. To give more accurate statements,
consider all irreducible compact leaves of $\fol$ contained in $\calk$.
Obviously we can assume there are finitely many $D_1, \ldots ,D_r$ of those. Besides $A_{\calk} = D_1 \cup \cdots \cup D_r$.
We shall say that these curves {\it contain a loop}\, if there are pairwise distinct points $p_{i_1}, \ldots ,
p_{i_s}$ with $p_{i_j} \in D_{i_j} \cap D_{i_{j+1}}$ for $1\leq j <s$ and $p_{i_s} \in D_{i_s}
\cap D_{i_1}$. We can now state a sharper form of Theorem~(\ref{fim1}).

\begin{teo}
\label{fim2}
Let $\fol$ and $T$ be as in Theorem~\ref{fim1}. Suppose that $\calk$ contains only
finite many irreducible compact curves $D_1, \ldots ,D_r$ invariant by $\fol$ and that these curves
do not contain loops. Suppose also that the remaining leaves of $\fol$ are dense in $\calk \setminus A_{\calk}$.
Then all (deformed) $\calh$-trajectories in $\calk$ have length smaller than some
positive constant ${\rm Const}$. An analogous statement is valid for the trajectories of
$\calh^{\theta}$-trajectories, $\theta \in (-\pi/2, \pi/2)$.
\end{teo}

Theorem~\ref{fim1} is an immediate consequence of Theorem~\ref{fim2}. In fact, in the context of Theorem~\ref{fim1} (before performing the normalizations
associated with Lemma~\ref{revision2}), we assume that the lengths of all deformed $\calh$-trajectories contained in $\calk$ are not uniformly bounded.
In particular there is a compact set $\calk^0 \subset \calk$ and a self-accumulating infinite brach $l^+$ of a deformed $\calh$-trajectory that is contained in $\calk^0$.
To pass from this situation to the context of Theorem~\ref{fim2}, let us now perform the normilzations of
Lemma~\ref{revision2}. in view of Theorem~\ref{fim2} the irreducible components of $A_{\calk}$ must form a loop.
Since the components of $A_{\calk}$ introduced in the course of the normalization procedure in question are rational
curves contained in pairwise disjoint tree-like arrangements, the only way for the components of $A_{\calk}$ to form a loop arises from the existence
of algebraic curves in the initial minimal set $\calk$. Thus Theorem~\ref{fim1} follows.

The rest of this paper is ultimately devoted to the proof of Theorem~\ref{fim2} for the case $\calh=\calh^0$. The extension to
$\calh^{\theta}$ for $\theta \in (-\pi/2, \pi/2)$ will be left to the reader. For this we
assume that $\fol$, $\omega$, $\omega_1$ and so on, satisfy all the conditions in
the statement of Lemma~\ref{revision2}. In particular the so-called ``local invariance condition''
of Section~3 concerning singularities of $\fol$ that belong to the Siegel domain is
verified. Also, modulo fixing a neighborhood $\Ww$ of the singular set of $\fol$, holonomy maps of $\fol$ obtained by
means of (segments of) deformed $\calh$-trajectories contained in the complement of $\Ww$ must satisfy the conclusions of
Theorem~\ref{Totalversionofblm}. Besides $l^+$ and $\calk^0$ will always be as indicated above. Now we have:

\begin{lema}
\label{6lema2}
We have $\calk^0 \cap {\rm Sing}\, (\fol)  \neq \emptyset$.
\end{lema}

\noindent {\it Proof}\,: It is a simple application of Theorem~\ref{Totalversionofblm}.
Suppose that the statement is false. Thus modulo reducing $\Ww$ we can assume that $\calk$ lies entirely in the complement
of $\Ww$. Now consider
a parametrization $c$ for a trajectory of $\calh$ such that $c(0) =p \in \calk^0$ (where $p$ does not belong to any critical region).
Fix a local transverse section $\Sigma_p$ and a disc $B_p (r) \subset \Sigma_p$
as in item~1 of Theorem~\ref{Totalversionofblm}.
Finally, for a fixed $t_0 \in \R_+$, let ${\rm Hol}\, (c_{t_0})$ be the holonomy map
associated to the restriction of $c$ to $[0,t_0]$.

By construction, ${\rm Hol} \, (c_{t_0})$ is defined on $B_p (r)$ for every $t_0$.
On the other hand, since the leaves of $\calk$ are dense in $\calk^0$, there is a
sequence of times $t_0^1 , t_0^2 ,\ldots$ going to infinity and such that
$\{ c (t_0^i) \}$ converges to $p$ when $i \rightarrow \infty$. Because $\Sigma_p$
is not a transverse section for $\calh$ at $p$, we cannot ensure {\it a priori}\, that
$c(t_0^i)$ can be chosen in $\Sigma_p$ for every $i \in \N$. Nonetheless, modulo
performing a slight modification of the trajectories of $\calh$ on a neighborhood of $p$
(similar in spirit to the ``deformed'' trajectories arising from the first critical regions)
this assumption can be made without loss of generality.

Now, thanks to the second part of the statement of Theorem~\ref{Totalversionofblm}, it follows that
the image of $B_p (r)$ under ${\rm Hol} (c_{t_0})$ is contained in a disc of radius
$r/10$ about $c(t_0^i) \in \Sigma_p$ provided that $i$ is large enough. In other
words, for $i$ very large ${\rm Hol} \, (c_{t_0})$ takes the disc $B_p (r)$ inside
itself. This actually implies the existence of a loop with hyperbolic holonomy for
$\fol$. As already seen, this gives a contradiction in the present case since $T$ is a diffuse,
cf. Lemma~\ref{atomicmass}.\qed

Recall that, in principle, singularities of $\fol$ lying in $\calk^0$ either are Siegel singularities (possibly associated to critical regions
of second type) or are irrational foci as in Remark~\ref{localregionsandirrationalfoci}. The latter singularities are however necessarily avoided by
the trajectories of $\calh$ as it was discussed in Section~5. Therefore irrational focus singularities can be ignored in our context and we can
suppose that all singularities of $\fol$ lying in $\calk^0$ are, in fact, Siegel singularities.

Let us now state a proposition that plays a key role in the proof of Theorem~\ref{fim2}.

\begin{prop}
\label{l+closed}
The above mentioned trajectory l$^+$ and set $\calk^0$ can be chosen so that $l^+ \subset \calk^0$ is closed.
\end{prop}

In view of the assumption about minimality of $\calk^0$ with respect to $\calh$, the preceding proposition actually says that $\calk^0$
is reduced, in a suitable sense, to a closed trajectory $l^+$. On the other hand, recall that our definition of ``closed trajectory'' allows $l^+$ to
pass through singularities of $\fol$ lying in $\calk^0$ (which are necessarily Siegel singularities as already pointed out). Indeed a closed trajectory
must necessarily go through singularities of $\fol$ since otherwise it gives rise to a holonomy map of $\fol$ possessing a hyperbolic fixed point.
As already seen this forces the current $T$ to be concentrated over an algebraic curve contradicting its diffuse nature. On the other hand, the fact that a closed trajectory goes through
a singularity of $\fol$ implies the existence of at least one saddle connexion for $\fol$.

In the rest of this section Proposition~\ref {l+closed} is going to be proved. In the next section we shall use it to derive the proof of
Theorem~\ref{fim2} and of Theorem~B in the Introduction.

As in Section~4, we consider the connected components $E=E_1, E_2, \ldots$
of the compact curves invariant by $\fol$ and contained in $\calk$. Hence each
$E_i$ consists of a number of irreducible curves $D_{i_k} \subset \{ D_1, \ldots ,D_r\}$,
Note that our
terminology allows $E_i$ to be empty i.e. reduced to a Siegel singularity that does
not belong to any compact curve invariant by $\fol$. To make the subsequent discussion
more transparent we shall first consider the following special situation:

\noindent {\bf First Case}: there is a unique connected component $E$.

The general case can easily be deduced from our discussion as it will be shown at the end of the section. Let us then fix a
singularity $p \in E \cap \calk^0$ (lying away from the critical regions). The next lemma allows us to suppose in addition that
a segment of $l^+$ is contained in a local separatrix of $\fol$ at $p$.

\begin{lema}
\label{6lema3}
There is a (deformed) semi-trajectoryof $\calh_{\calk^0}$ contained
in a separatrix of $\fol$ and having $p$ as an accumulation point (where $\calh_{\calk^0}$
stands for the restriction of $\calh$ to $\calk^0$). This trajectory will still be denoted by $l^+$.
\end{lema}

\noindent {\it Proof}\,: Fix local coordinates $(u,v)$ around $p \simeq (0,0)$ in which
the $1$-form $\omega$ defining $\fol$ satisfies Equation~(\ref{siegel2}). We choose $u,v$
so that the trajectories of $\calh$ in $\{ v=0 \}$ converge to $p \simeq (0,0)$. Next let
$\Sigma_{\theta}$ be a local transverse section passing through the point $(
e^{2\pi i \theta} ,0)$, $\theta \in [0, 2\pi)$. Since $p \in \calk^0$,
Proposition~(\ref{3prop1}) implies the existence of a sequence of points
$(\theta_i , v_i)$ such that $(e^{2\pi i \theta_i} , v_i) \in \calk^0$ and
$\vert v_i \vert \rightarrow 0$. Since $\calk^0$ is closed, it follows the existence of
$\theta_{\infty}$ such that $(e^{2\pi i \theta_{\infty}} ,0) \in \calk^0$. The trajectory
of $\calh$ through this point then satisfies the required conditions.\qed

\begin{obs}
\label{vamos}
{\rm Without loss of generality we can suppose that $\theta_{\infty} =0$ so that
$(1,0) \in \calk^0$. We also set $\Sigma =\Sigma_0$ and denote by $l$ the trajectory
of $\calh$ through $(1,0)$ which is obviously contained in $\calk^0$.}
\end{obs}

\begin{lema}
\label{6lema3.5}
$l^+$ is not entirely contained in $E$.
\end{lema}

\noindent {\it Proof}\,: Suppose for a contradiction that $l^+$ is entirely contained
in $E$. Suppose that $D_1$ is the irreducible component of $E$ containing $p$.
We can assume without loss of generality that $D_1$ contains the whole
of $l^+$. Indeed, suppose that by following $l^+$ one passes from $D_1$ to
another irreducible component $D_2$. This passage is then made through a
singularity $q_{1,2} =D_1 \cap D_2$ which belongs to the Siegel domain. By assumption
the orientation of the trajectories of $\calh$ around $q_{1,2}$ (always given by
Lemma~\ref{3lema1}) is such that they go from the separatrix contained in $D_1$
to the separatrix contained in $D_2$. Hence the trajectory $l^+$ cannot return to $D_1$
through $q_{1,2}$. Because $D_1 \cap D_2 = q_{1,2}$, this trajectory cannot return
to $D_1$ through any point in $D_2$. Since the ``graph of irreducible components''
associated to $E$ contains no loop, we
conclude that $l^+$ will never return to $D_1$. In other words, if $l^+ \subset E$,
then $l^+$ will eventually be ``captured'' by an irreducible component of $E$
that can be supposed to be $D_1$.

On the other hand, the trajectory $l^+$ cannot approach any singularity $q \in D_1$ where the
leaves of $\calh$ restricted to $D_1$ approach $q$.  Otherwise $l^+$ would
leave $D_1$ by means of the separatrix of $\fol$ at $q$ which is transverse
to $D_1$. This, in fact, implies that the complement of a compact part of
$l^+$ does not accumulate on any singularity. Since $D_1$ is compact and $l^+$ is of infinite length,
Theorem~\ref{Totalversionofblm} can be employed to ensure that the holonomy group of $D_1 \setminus {\rm Sing}\, (\fol)$,
w.r.t. the foliation $\fol$, contains a hyperbolic element. This is however impossible
since it would force $T$ to be concentrated over $D_1$.\qed

Recalling Remark~(\ref{vamos}), we can suppose that $l^+$ arrive to $E$
through $p$. In other words, around $p$ there is a separatrix of $\fol$
transverse to $E$ (given by $\{ v=0\}$ and denoted by $S_p$)
in the local coordinates $(u,v)$
used in the proof of Lemma~(\ref{6lema3}), with $\{ u=0\} \subset E$.
Since $E$ is a connected component of the set of all compact curves invariant by $\fol$ that are contained in
$\calk$, it follows that this separatrix
is contained neither in the divisor of zeros and poles of $\omega$ nor in $(\omega_1)_0$.
Similarly, in coordinates $(w,t)$ around $q$, $\{ w=0\} \subset E$, there is a separatrix
of $\fol$ at $q$ which is transverse to $E$ (given by $\{ t=0\}$ and
denoted by $S_q$. Again this separatrix is not contained in $(\omega)_0 \cup (\omega)_{\infty} \cup (\omega_1)_0$.

To prove Proposition~\ref {l+closed} let us suppose for a contradiction the existence of a $\calh$-trajectory $l^+$ in $\calk^0$
satisfying the above conditions but which is not a closed trajectory passing through singular points of $\fol$. Recall that
$l^+$ accumulates on itself, i.e. it has non-trivial recurrence. By using
the local coordinates $(u,v)$, $(w,t)$ introduced above, the non-trivial recurrence
of $l^+$ implies the existence of points $(u_n ,v_n) = (e^{2\pi i \theta_n} ,v_n)$
satisfying the following:
\begin{enumerate}
\item $(u_n ,v_n) = (e^{2\pi i \theta_n} ,v_n)$ belongs to $l^+$ for every $n \in \N$.

\item Both sequences $\{ \theta_n \}, \, \{ v_n \}$ converge to zero.
\end{enumerate}
Actually we can be more precise. Let $U_E$ be a small ``tubular'' neighborhood
of $E$. Let us consider the ``full'' sequence of ``first returns'' of $l^+$ to $U_p$ which will
still be denoted by $(e^{2\pi i \theta_n} ,v_n)$. If $U_E$ is appropriately chosen, then
the local connected component of $l^+$ through $(e^{2\pi i \theta_n} ,v_n)$ satisfies the
conclusions of Proposition~(\ref{3prop1}). We have:

\noindent {\bf Claim}: We can assume that $\theta_n =0$ for every $n$.

The above assumption is not really needed from a strict point of view. It simply allows
us to shorten our discussion which applies equally well to the general case. It can
also be formalized by again locally deforming the leaves of $\calh$ on a neighborhood
of the circle $(e^{2\pi i \theta} ,0)$, $\theta \in [0,2\pi)$. This deformation is
essentially given by the local holonomy of $\{ v=0\}$ and does not affect neither the global
dynamics of $l^+$ nor the estimates involved in the holonomy of (compact pieces of)
$l^+$.

Summarizing, in what follows we assume that $l^+$ enters the
``tubular neighborhood''
$U_E$ by means of a sequence of points having the form $P_n= (1, v_n)$.
Besides the sequence $\{ \vert v_n \vert \}$ converges to $0 \in \C$. In particular,
these points belong to $\Sigma_{\rm in}$, a local transverse section through the point
$(1,0)$ in $(u,v)$-coordinates.
In the sequel we sometimes identify the point $(1,v) \in \Sigma_{\rm in}$
with the point $v \in \C$, thus identifying $\Sigma_{\rm in}$ itself with a neighborhood of $0 \in \C$.

Let us briefly review what is the nature of the holonomy map of $\fol$ associated to a trajectory of $\calh$ as
above. We begin with a definition
on $U_E$ involving the generalized Dulac transform introduced in Section~4.
It is clear that the segment of $l^+$ delimited by the points $p,q$ above verifies the
condition discussed in Section~4 in connection with the generalized Dulac transform.
Let then $\Sigma_{\rm out}$ be a transverse section through the point $(1,0)$ in
$(w,t)$-coordinates (on a neighborhood of $q$) so that the corresponding generalized
Dulac transform is well-defined. The last
statement can be made precise as follows: Let $V_0 \subset \Sigma_p$ be a simply connected domain containing
a point $(1,z) \in \calk^0 \cap \Sigma_{\rm in}$. According to Sections~3,~4, the oriented trajectory
$l_{(1,z)}$ of $\calh$ through $(1,z)$ intersects $\Sigma_{\rm out}$ at a point $(z' ,1) \in
\calk^0 \cap \Sigma_q$.
Parameterizing by $c: [0,1] \rightarrow l_{(1,z)}$ the segment of $l_{(1,z)}$ delimited by
$(1,z)$ and $( z' ,1)$, we ask the generalized Dulac transform ${\rm GDul} : V_0
\rightarrow
\Sigma_{\rm out}$ to be well-defined w.r.t. the path $c$ (in the sense of Sections~3, 4).
Obviously
${\rm Dul} (1,z) = ( z', 1)$. Modulo reducing $U_E$, we can suppose that
$U_E$ is the saturation of $\Sigma_{\rm out}$ by $\fol$. Finally let $\lambda$ denote the exponent associated with this
generalized Dulac transform as in the context of Proposition~\ref{4.5prop2}.

Now consider the compact set $\calk^0 \setminus U_E$ where the foliation $\fol$ is regular.
Here the holonomy associated to $\fol$ and to a segment of the trajectory $l^+$ contained in $\calk^0 \setminus U_E$ has a clear meaning.
Besides if $c: [0,1]   \rightarrow \calk^0 \setminus U_E$ is a path parameterizing a segment of $l^+$, then
${\rm Hol}\, (c)$ satisfies the conclusions of Theorem~\ref{Totalversionofblm}, in particular
${\rm Hol}\, (c)$ is defined on a transverse disc of uniform radius $\delta >0$ (regardless of the length of $c$).

Resuming the notations of Proposition~(\ref{4.5prop2}), there are two
cases to be considered according to whether or not the value of $\lambda$ is rational.

\vspace{0.1cm}

\noindent $\bullet$ Let us first suppose that $\lambda \in \R \setminus \Q$.

Recall that $\lambda$ is less than or equal to~$1$.
Let $l^+_{p,q}$ denote the segment of $l^+$ delimited by $p,q$. By choosing a point
in this segment, we can consider the holonomy of $\fol$ generated by the
local holonomy maps associated to the separatrices of the singularities of $\fol$ lying in
$l^+_{p,q}$ (in particular $p,q$). This group is Abelian. Otherwise, a non-trivial
commutator in this group would be ``parabolic'' in the sense that it is tangent to
the identity. As in Lemma~(\ref{4.5lema2}), the existence of this local diffeomorphism
would yield a contradiction since $l^+ \subset \calk$ and $T$ is diffuse having $\calk$
as its support.

Clearly the local holonomy maps associated to the above mentioned singularities can be
identified with elements of ${\rm Diff}\, (\C ,0)$ by appropriately choosing
local transverse section and (ordinary) Dulac transforms. The fact that $\lambda
\in \R \setminus \Q$ then implies the existence of a holonomy map $h$ in this group that, under the
above identification, is an element of ${\rm Diff}\, (\C ,0)$ whose linear part is an irrational rotation.
Then the following elementary statement holds:

\noindent \textsc{Fact 1}: Given an arbitrary $l \in \N^{\ast}$, there is a constant $C = C(l ,2\pi \lambda)$
such that, for every $r \in \R_+$ sufficiently small and $z \in \C$ with
$\vert z \vert = r$, the sets $B_z (C .r), h(B_z (C.r)), \ldots , h^{l-1} (B_z (C.r))$
are pairwise disjoint (where $B_z (C.r)$ stands for the ball about $z$ of radius $C.r$).

Now consider a strictly monotone sequence $\{ r_j \} \subset \R_+$ converging
to {\it zero}. For each $j$, denote by $B (r_j) \subset \Sigma_{\rm in}$ the ball of radius
$r_j$ about $0 \simeq (1,0) \in \Sigma_{\rm in}$. Next recall
that $(1, v_n) \in \Sigma$ is the sequence of the ``first returns'' of $l^+$ to $U_E$.

\vspace{0.2cm}

\noindent {\it Proof of Proposition~\ref{l+closed} when $\lambda \in \R \setminus \Q$}\,: We begin by fixing $l \in \N$
larger than $2\pi/ \lambda$. Next we consider a constant $C =
 C(l ,2\pi \lambda)$ as in Fact~1. Finally let us denote by
$\mu ,\nu$ suitable measures
representing the current $T$ on the transverse sections $\Sigma_{\rm in}, \Sigma_{\rm out}$ respectively.

Consider a ball $B(r) \subset \Sigma_{\rm in}$ with $r$ very small. Let $V_{\alpha} (r) \subset
B_{r}
\subset \Sigma_p$ be a sector of angle $\alpha < 2\pi \lambda$ and radius
$r_j$. The invariance of $\mu$ under the local diffeomorphism $h$ obtained as
the holonomy map of $\{ v=0\}$ implies that, for every $\epsilon >0$ fixed, one has
\begin{equation}
\frac{\mu (V_{\alpha} (r))}{\mu (B(r))} \geq (1-\epsilon) \lambda
\label{lafoi1}
\end{equation}
provided that $r$ is very small and that $\alpha$ is sufficiently close to
$2\pi \lambda$.

Now let us choose $j_0 \in \N$ very large and denote by $v_{n_0}$ the smallest
positive integer $n$ such that $v_n \in B(r_{j_0})$. Set $r_0 = \vert v_{n_0} \vert$
and denote by $B(r_0)$ (resp. $B(2r_0)$) the ball of radius $r_0$ (resp. $2r_0$)
about $0 \simeq (1,0) \in \Sigma_{\rm in}$. For $\alpha$ very close to
$2\pi \lambda$, let us denote by $V_{\alpha} (2r_0)$ the sector
of angle $\alpha$ and radius $2r_0$ which is divided into two equal parts by the
semi-line joining $0$ to $v_n$. Modulo taking $r_0$ sufficiently small (i.e. $j_0$ large
enough), we can consider a generalized Dulac transform ${\rm GDul}$ which is well-defined
on $V_{\alpha} (2r_0)$ (w.r.t. some path $c$ fixed once and for all).

The image $W(2r_0)$ of  $V_{\alpha} (2r_0)$ under ${\rm GDul}$ is contained
in a disk $B' (r_0') \subset \Sigma_{\rm out}$ of radius $r_0' \simeq r_0^{1/\lambda}$
(cf. Proposition~\ref{4.5prop2}). Again by choosing $r_0$ small enough, this estimate
yields $r_0' < 2C r_0$ where $C$ is the constant fixed above. Also we know that
$\nu(W(2r_0)) > (1- \epsilon) \mu (V_{\alpha} (2r_0))$.

Finally let $c: [0,1] \rightarrow \calk^0$
be a parametrization of the segment of trajectory $l^+$ going from
the first point in which $l^+$ intersects $\Sigma_q$ to the point $(1,v_{n_0}) \in
\Sigma$.
Note that $c$ remains away from the neighborhood $U_E$.  In fact, its distance to $U_E$
is bounded below by a uniform constant times $r_0$. In particular, for $r_0$ small enough,
this implies that the holonomy map ${\rm Hol}\, (c(t))$ is well-defined on $W(2r_0)$
for every $t \in [0,1]$. Besides Theorem~\ref{Totalversionofblm} also ensures that ${\rm Hol}\, (c) (W(2r_0))$ is contained in a disc of radius less
than $2r_0'$ about $v_{n_0} \simeq (1,v_{n_0}) \in \Sigma$. Since $\vert v_{n_0} \vert
=r_0$, this ensures that ${\rm Hol}\, (c) (W(2r_0)) \subset B(2r_0)$. Furthermore we have $\mu [{\rm Hol}\, (c) (W(2r_0))] = \nu (W(2r_0))$
since $\mu, \, \nu$ are local transverse measures representing the invariant current $T$. Thanks to
Estimate~(\ref{lafoi1}), it follows that
\begin{equation}
\mu [{\rm Hol}\, (c) (W(2r_0))] \geq (1-\epsilon) \lambda
\mu (B(2r_0)) \, . \label{lafoi2}
\end{equation}
Finally, since $r_0' < 2C r_0$, the set ${\rm Hol}\, (c) (W(2r_0))$ has $l$ images
pairwise disjoint under the holonomy map $h$. Since all these images have the same
$\mu$ measure, the total measure of their union is $l .\mu [{\rm Hol}\, (c) (W(2r_0))]
> \mu (B (2r_0))$ in view of the choice of $l$ and for $\epsilon$ very small.
This yields the desired contradiction since the union of these images is
contained in $B (2r_0)$. Therefore our statement is proved
in the present case.\qed

\vspace{0.2cm}

\noindent $\bullet$ Let us now suppose that $\lambda \in \Q$.

This case is somehow similar to the preceding one. The main difference is that
Fact~1 no longer holds. On the other hand, we know
that $\fol$ is linearizable around every singularity in $l^+_{p,q}$
(Proposition~\ref{4.5prop1}). Yet the
local holonomy maps associated to these singularities generate an Abelian group which is therefore finite and hence conjugate
to a (finite) group of rotations. This will allow us to make precise asymptotic
calculations so as to dispense with the ``$\epsilon$ margin'' involved in preceding discussion.
In fact, standard arguments involving this group and the nature of the singularities of $\fol$ contained
in $E$ shows that the restriction of $\fol$ to $U_E$ admits a non-constant holomorphic first integral (see \cite{mamo},
\cite{paul}). The existence of this integral however will not be necessary in what follows.

Let us resume the notations of the case where $\lambda$ was not rational. Since the above mentioned holonomy group associated
to the curve $L$ is conjugate to a finite group of rotations, we can find a local diffeomorphism $h$ in this group which is itself a
rotation of angle $2\pi \lambda$. Let $V_{2\pi \lambda} (r) \subset \Sigma_{\rm in}$ be a sector of angle
$2\pi \lambda$ and radius $r$. If $B (r) \subset \Sigma_p$ is the corresponding
ball of radius $r$, we obviously have $\mu [V_{2\pi \lambda_2/\lambda_1}
(r) ] = \lambda \mu (B(r))$.

\vspace{0.2cm}

\noindent {\it Proof of Proposition~\ref{l+closed} when $\lambda \in \Q$}\,:  Again
choose $j_0 \in \N$ very large and denote by $v_{n_0}$ the smallest
positive integer $n$ such that $v_n \in B(r_{j_0})$. Set $r_0 = \vert v_{n_0} \vert$.
Next let $V_{2\pi \lambda} (2r_0)$ be the sector
of angle $2\pi \lambda$ and radius $2r_0$ which is divided
into two equal parts by the semi-line joining $0$ to $v_n$.
Modulo taking $r_0$ sufficiently small (i.e. $j_0$ large
enough), we can consider a generalized
Dulac transform ${\rm GDul}$ which is well-defined
on $V_{2\pi \lambda} (2r_0)$ (w.r.t. some path $c$ fixed once and for all).
The image $W(2r_0)$ of  $V_{2\pi \lambda} (2r_0)$ under ${\rm GDul}$ is contained
in a disc $B' (r_0') \subset \Sigma_{\rm out}$ of radius $r_0' \simeq r_0^{1/\lambda}$
(cf. Proposition~\ref{4.5prop2}). We note that the possibility of having $\lambda =1$ is not {\it a priori} excluded
so that we first suppose $\lambda < 1$. As in the proof of the case ``$\lambda$ irrational'', the above mentioned disc is taken by
the holonomy of $l^+$ to a set ${\rm Hol}\, (c) (W(2r_0)) \subset B(2r_0) \subset
\Sigma_{\rm in}$.
Furthermore we have
$$
\lambda \mu (B(2r_0))  = \mu [V_{2\pi \lambda} (2r_0)]
= \nu [W(2r_0)] = \mu [{\rm Hol}\, (c) (W(2r_0))] \, .
$$
However the set ${\rm Hol}\, (c) (W(2r_0))$ has the ``denominator of $\lambda \in \Q$''
images pairwise disjoint
under the iterations of the rational rotation $h$ of angle $2\pi \lambda$.
Since they are all contained in $B(2r_0)$, it follows that $\lambda$
is the inverse of an integer.
In this case $\mu [\bigcup_{i=1}^{1/\lambda} h^i ({\rm Hol}\, (c) (W(2r_0)))] = \mu (B (2r_0))$.
Finally because $0 \simeq (1,0) \in \Sigma_{\rm in}$ is not contained in the closure of
the set $\bigcup_{i=1}^{\lambda_1} h^i ({\rm Hol}\, (c) (W(2r_0)))$, it follows that
a sufficiently small neighborhood of $0 \in \Sigma_{\rm in}$ has $\mu$-measure {\it zero}.
Since $\calk^0$ is in the support of $[T]$, it must not intersect the neighborhood in
question. This yields the desired contradiction since $\calk^0$ accumulates on $p$.

Finally if we have $\lambda=1$, then ${\rm GDul}$ yields an identification of neighborhoods of the origin
in $\Sigma_{\rm in}, \,\Sigma_{\rm out}$ which takes $\mu$ to $\nu$. In other words, the segments of $l^+$ passing through
$U_E$ behaves as they remained ``away from the singular set of $\fol$''. In this case the conclusion follows simply from the contractive
behavior of the holonomy of $\fol$ defined with the help of the segments of $l^+$ contained in the complement of $U_E$. The proof of
Proposition~\ref{l+closed} in the case corresponding to the connectedness of $E$ is now over.\qed

\vspace{0.1cm}

\noindent {\it Proof of Proposition~\ref{l+closed} in the general case}\,: To finish this section, let us now show how the previous arguments
can naturally be adapted to yield the proof of Proposition~\ref{l+closed} when $E$ contains more than one connected component.
Suppose first that, instead of a single connected component $E$, there are two
connected components $E_1, E_2$. We consider a trajectory $l^+$ of $\calh$ arriving
to $E_1$ through a singularity $p_1$ and leaving $E_1$ through another singularity $q_1$
(as in Lemmas~\ref{6lema3} and~\ref{6lema3.5}). The first possibility that
may occur is a ``saddle-connection'' between $q_1$ and a singularity $p_2 \in E_2$.
More precisely, it may happen that the separatrix of $\fol$ at $q_1$ which is
transverse to $E_1$ coincides with a separatrix of another singularity $p_2 \in
E_2$ of $\fol$. Then $l^+$ will arrive to $E_2$ through $p_2$.

The solution of this first difficulty is provided by the proof of Proposition~\ref{4.5prop2}.
In fact, by using the leaf of $\fol$ joining $q_1$ to $p_2$, we can define a
new ``generalized Dulac transform'' encompassing both $E_1, E_2$ as they were
a single connected component. This definition is straightforward and the
proof of Proposition~\ref{4.5prop2} shows that the resulting ``Dulac transform''
still satisfies the conditions given in the statement in question.

Suppose now that there is no ``saddle-connection'' in the sense described above
between $E_1 ,E_2$. The difficulty here arises from the fact that $l^+$ may accumulate
on $E_2$ before returning to $E_1$. Again we suppose that $l^+$ arrives to $E_1$
(resp. $E_2$) through a singularity $p_1$ (resp. $p_2$) and leaves it through
a singularity $q_1$ (resp. $q_2$). Let
$\textsc{S}_{p_2}$ denote the separatrix of $\fol$ at $p_2$ which is transverse
to $E_2$.

Again we keep the notations used in the course of this section. Let us then
consider the image
$W(2r_0)$ of  $V_{\alpha} (2r_0)$ under the generalized Dulac transform associated
to $E_1$, ${\rm GDul}_1$.
As $l^+$ continues from $q_1$ to $p_2$, we let ${\rm Hol} \, (W(2r_0))$
denote the image of $W(2r_0)$ by the corresponding holonomy map. The special
difficulty here is that, when approaching $p_2$, $l^+$ may become ``very close'' to
$\textsc{S}_{p_2}$. In particular, with the obvious identifications, it may happen
that ${\rm Hol} \, (W(2r_0))$ contains $\textsc{S}_{p_2}$. This situation prevents us
from considering the ``Dulac transform'' associated to $E_2$, ${\rm GDul}_2$,
as being defined on all
of ${\rm Hol} \, (W(2r_0))$. To deal with this case, we proceed as follows. First we
note that, in the hard case, this phenomenon should occur for
``every sequence of returns'' of
$l^+$ to the transverse section in which ${\rm GDul}_1$
is defined. Under this assumption, we substitute $l^+$ by a trajectory
$l^+_2$ of $\calh$ which has properties analogues to those of $l^+$ and, furthermore,
is contained in $\textsc{S}_{p_2}$. The existence of this trajectory is clear since $l^+$
accumulates on $\textsc{S}_{p_2}$. We then start our argument with $l^+_2$ so
that the Dulac transform associated to $E_2$ is automatically well-defined. We claim
that, for $l^+_2$,
the Dulac transform associated to $E_1$ will also be well-defined on the appropriate
domain. It is clear that the desired statement results easily from this claim.

To check the claim,
we note that $l^+$ is close to $l^+_2$ at the same order of the diameter of
$W(2r_0)$ near $p_2$. In turn, this diameter is small when compared to $r_0$
(cf. Propostion~\ref{4.5prop2}). If $l^+, l^+_2$ remain close to each other for
all time, then $l^+_2$ reaches the domain
of definition of ${\rm GDul}_1$ at a point ``very close''
to a point of return of $l^+$ to this domain (i.e. the distance between these two points
is small to order superior to the distance of the return point to $E_1$). Since the
diameter of the corresponding $W(2r_0), \, {\rm Hol} \, (W(2r_0))$ is also small, the claim
follows at once. This therefore completes the argument modulo the assumption that
$l^+, l^+_2$ remain close to each other for all time. This assumption however can
always be made. In fact, we first observe that the leaf $L_2$ of $\fol$ containing
$l_2^+$ approaches the trajectory $l^+$ due to the contracting behavior of the
holonomy along $l^+$. By a simple argument of continuous dependence for solutions
of differential equations, $l^+, \, l_2^+$ remain close for an {\it a priori}\, fixed
period of time. But at the end of this period of time, we can modify $l^+_2$ inside
$L_2$ by adding a short line (transverse to $\calh$) so as to bring the modified
trajectory $l_2^+$ close again to $l^+$. This new trajectory $l_2^+$ satisfies all
the previous requirements and establishes the claim.

Now it is clear that the existence of several $E_1, E_2, \ldots$ connected components
does not pose any new intrinsic difficulty. The proof of Proposition~\ref{l+closed} is finally completed.\qed

\section{Proofs for the main results}

To be able to prove the theorems stated in the Introduction, we shall need to consider the closed
trajectory $l^+$ whose existence is ensured by Proposition~\ref{l+closed}. This trajectory will be referred to as a {\it singular}\, closed trajectory since it passes
through the singularities of $\fol$. In what follows, we shall keep the notations and the terminology of the
preceding section.

Denote by $L_0, L_1, \ldots, L_m$ the leaves of $\fol$ that contain a non-trivial segment of $l^+$. For
$i=0,\ldots, m-1$, $L_i$ intersects $L_{i+1}$ at a singularity of $\fol$ belonging to the Siegel domain. Also
$L_k$ intersects $L_0$ at a Siegel singularity so that the leaves $L_0, L_1, \ldots, L_m$ form a loop
by means of their saddle-connexions.

Fix a base point $p \in l^+ \cap L_0$ and consider a local transverse section $\Sigma$ at $p$. To prove our
main results we are going to consider the pseudogroup $\Gamma$ of transformations of $\Sigma$ obtained by
the collection of first return maps over paths contained in the leaves $L_0, L_1, \ldots, L_m$. More generally denote
by $\mathcal{L}$ the union of the (finitely many) leaves de $\fol$ defined by the following rules:
\begin{itemize}

\item $L_0, L_1, \ldots, L_m$ belongs to $\mathcal{L}$.

\item If $L$ belongs to $\mathcal{L}$ and $L$ (locally) defines a separatrix for a Siegel singularity $p$ of
$\fol$, then the global leaf obtained from the other separatrix of $\fol$ at $p$ must also belong to $\mathcal{L}$.
\end{itemize}

\noindent The pseudogroup $\Gamma$ is then obtained by means of first return maps defined over all paths contained
in $\mathcal{L}$.

One first element $f$ of $\Gamma$ corresponds of course to the singular loop $l^+$. As already seen, we can suppose
that $f$ is a (ramified) map of the form $f(z) = z^{\lambda} (1 + u(z))$ where $\lambda > 1$ and $u(0)=0$ ($u$ being defined
on a neighborhood of $0 \in \C$). Note that $\lambda$ need
not be an integer so that $f$ should be thought of as a ``ramified'' map. Yet, in sectors of angle slightly smaller than
$2\pi /\lambda$, $f$ is well-defined and one-to-one onto its image.

Another element of $\Gamma$, denoted by $g$, corresponds to the local holonomy map arising from the singularities
of $L_0, \ldots, L_m$. Observe that at least one of the Siegel singularities of $\fol$ lying in $l^+$ must have eigenvalues
different from $1,-1$. In fact, otherwise $f$ would have no ramification and thus it would consist of a hyperbolic contraction
defined on a neighborhood of $p \in \Sigma$. Therefore every invariant measure on $\Sigma$ would automatically be concentrated
at $p$ what it is impossible. Summarizing, we conclude that $\Gamma$ contains
an element $g$ defined about $p \simeq 0 \in \C$ and having the form
$$
g(z) = e^{2\pi i \alpha} z + {\rm h.o.t.}
$$
with $\alpha \in (0,1)$. This leads us to study the dynamics of a pseudogroup $\Gamma$ containing elements $f,g$
as above.

Recall that $\lambda >1$. If $\lambda$ were an integer, then the classical theorem of B\"ottcher would provide coordinates
where $f(z) =z^{\lambda}$. Hovewer, in general, the map $f$ is ramified so that it is well-defined only on suitable sectors. The
``position'' of the sectors where we want to define a particular determination of $f$ are naturally permuted by means of the
local holonomy maps associated to the Siegel singularities of $\fol$ lying in $l^+$. Indeed, the ambiguity in the definition of $f$
as compositions of ordinary holonomy maps and suitable Dulac transforms is precisely codified by the local holonomy
maps arising from the mentioned singularities. In particular two different determinations of $f$ commute in the obvious sense with the corresponding local holonomy
maps. These elementary facts will freely be used in what follows.

Now, even though $\lambda$ is not an integer, the method of B\"ottcher still provides a conjugacy between $f$ and
$z \mapsto z^{\lambda}$ over appropriate sectors. Whereas the conjugacy map is clearly not defined about $0 \in \C$, it
has all natural asymptotic properties at $0 \in \C$. By using one of these coordinates, we can suppose that
$\Gamma$ contains the map $f(z) =z^{\lambda}$ along with at least one map $g$ of the form $g(z) = e^{2\pi i \alpha} z
+r(z)$ where $\alpha \in (0,1)$ and $\Vert r(z) \Vert \leq C \Vert z \Vert^2$ for some constant $C$. Furthermore different
determinations of $f$ are naturally permuted by $g$.

Before continuing let us make some elementary remarks about the function ``$k^{\rm th}$-root''. More precisely
let $k \in \R$, $k > 1$, be fixed. Consider the map $z \longmapsto (1 + z)^k$
which is well-defined for $\Vert z \Vert < 1/2$. The corresponding derivative is simply $(1+z)^{(1-k)/k} /k$.
In particular, for $\Vert z \Vert < 1/2$, the norm of its derivative is uniformly bounded by
\begin{equation}
\frac{1}{k} 2^{(k -1)/k} \leq \frac{2}{k} \, . \label{elementarybound}
\end{equation}
Next consider the element $h_1$ of $\Gamma$ defined by $h_1 (z) = f^{-1} \circ g \circ f(z)$
and note that $h_1$ is well-defined on a uniform sector (slightly smaller than the sector in which $f$ was defined). More generally
different determinations of $h_1$ are naturally permuted by $g$ since so are the determinations of $f$. Similarly we define
$$
h_n (z) =  f^{-n} \circ g \circ f^n(z) \; .
$$
Our first task is to show that the elements $h_n$ are defined on a uniform domain and that they converge to the identity on this domain.
For this let us set $g(z) = e^{2\pi i \alpha} z + c_2 z^2 + c_3 z^3 + \cdots$. Now note that $h_1$ admits the form
$$
h_1 (z) = e^{2\pi i \alpha} z (1 + c_2z^{\lambda} + c_3 z^{2\lambda} + \cdots)^{1/\lambda} \; .
$$
The expression $c_2z^{\lambda} + c_3 z^{2\lambda} + \cdots = r(z^{\lambda})/z^{\lambda}$ is clearly less than $1/2$
for $\Vert z \Vert$ sufficiently small. In particular $h_1 (z)$ is actually holomorphic on a neighborhood of $0 \in \C$. In addition,
Estimate~\ref{elementarybound} yields
\begin{equation}
\Vert h_1 (z) - e^{2\pi i \alpha/\lambda} z \Vert \leq \frac{2}{\lambda} C \Vert z \Vert^{\lambda} \label{elementarybound2}
\end{equation}
on the same domain. A direct inspection shows that $h_n (z)$ is holomorphic on the same neighborhood of $0 \in \C$
and that it satisfies the following estimate
\begin{equation}
\Vert h_n (z) - e^{2\pi i \alpha/\lambda^n} z \Vert \leq \frac{2}{\lambda^n} C \Vert z \Vert^{\lambda^{2n-1}} \; .\label{elementarybound3}
\end{equation}
Because $\lambda >1$, we obtain

\begin{lema}
\label{lemma7.A}
For $\tau >0$ sufficiently small, all the $h_n$ are holomorphic and well-defined on the disc $B_{\tau} (0)$ of radius $\tau$
about $0 \in \C$. Furthermore these diffeomorphisms converge uniformly to the identity on $B_{\tau} (0)$.\qed
\end{lema}

Recall that a vector field $Y$ defined on a neighborhood $U$ of $0 \in \C$ is said to {\it belong to the closure of $\Gamma$
(relative to $U$)}\, if for every $V \subset U \subset \C$
and $t_0 \in \R_+$ so that $\phi_Y^{t} (V)$ is well-defined for all $0 \leq t \leq t_0$, the map $\phi_Y^{t_0} : U \rightarrow \phi_Y^{t_0}
(U) \subset U$ is a uniform limit of elements of $\Gamma$ defined on $V$. Here $\phi_Y^{t}$ stands for the local flow generated by $Y$. From the
definition it follows that $\phi_Y^{t_0}$ is holomorphic as a uniform limit of holomorphic maps (contained in $\Gamma$).
Next we have:

\begin{prop}
\label{prop7.A}
The vector field whose local flow consists of rotations about $0 \in \C$ belongs to the closure of $\Gamma$
(relative to the disc $B_{\tau/2} (0)$). In other words, every rotation $R_{\beta} : z \mapsto e^{2\pi i \beta} z$
is a uniform limit on $B_{\tau/2} (0)$ of actual elements of $\Gamma$.
\end{prop}

\noindent {\it Proof}. Fix a rotation $R_{\beta}  (z) = e^{2\pi i \beta} z$. We need to find a sequence of elements in $\Gamma$
that approximate $R_{\beta}$ on $B_{\tau/2} (0)$. This sequence can explicitly be obtained as follows. For $n$ large
enough let $k_n$ be the integral part of $\beta \lambda^n/\alpha$. Clearly the linear part of $h_n$ at $0 \in \C$ is a rotation
of angle $[\beta \lambda^n/\alpha]\alpha/\lambda^n = k_n \alpha/\lambda^n$. In particular the difference $\vert
\beta - k_n \alpha/\lambda^n \vert$ is bounded by  $\alpha/\lambda^n$ which, in turn, tends to zero when $n \rightarrow
\infty$ (since $\lambda >1$). Therefore to establish the proposition it suffices to check that the sequence $\{ h_n^{k_n}\}_{n \in \N}
\subset \Gamma$ satisfies the two conditions below.
\begin{enumerate}

\item For $n$ very large, $h_n^{k_n}$ is well-defined on $B_{\tau/2} (0)$.

\item On $B_{\tau/2} (0)$, $h_n^{k_n}$ converges uniformly towards its own linear part at $0 \in \C$.

\end{enumerate}
These conditions will simultaneously be verified as consequences of Estimate~(\ref{elementarybound3}). To abridge notations, denote
by $R_n$ the rotation of angle $\alpha/\lambda^n$ about $0 \in \C$. The linear character of $R_n$ gives $D_zR_n =R_n$ for every
$z \in \C$. In particular the norm $\Vert D_zR_n \Vert$ is constant equal to~$1$. Next observe that, for $\Vert z \Vert$ sufficiently small,
Estimate~(\ref{elementarybound3}) yields
\begin{eqnarray*}
\Vert h_n^2 (z) - R_n^2 (z) \Vert & = & \Vert h_n^2 (z) - R_n \circ h_n (z) + R_n \circ  h_n (z) - R_n^2 (z) \Vert \\
& \leq & \Vert (h_n - R_n) \circ h_n (z) \Vert + \Vert R_n (h_n (z)) -R_n (R_n (z)) \Vert \\
& \leq & \frac{2C}{\lambda^n} \Vert h_n (z)\Vert^{\lambda^{2n-1}} + \sup_{B_{\tau/2} (0)} \Vert DR_n\Vert . \Vert
h_n (z) - R_n (z) \Vert \\
& \leq & \frac{2C}{\lambda^n} \Vert h_n (z)\Vert^{\lambda^{2n-1}} + \frac{2C}{\lambda^n} \Vert z \Vert^{\lambda^{2n-1}} \\
& = &  \frac{2C}{\lambda^n} (\Vert h_n (z)\Vert^{\lambda^{2n-1}} + \Vert z \Vert^{\lambda^{2n-1}}) \, .
\end{eqnarray*}
If we set $\Vert h_n^3 (z) -R_n^3 (z) \Vert \leq \Vert h_n (h_n^2(z)) - R_n (h_n^2 (z)) \Vert + \Vert   R_n (h_n^2 (z))
-R_n (R_n^2(z)) \Vert$ and repeat the above procedure, it follows that
$$
\Vert h_n^3 (z) -R_n^3 (z) \Vert \leq \frac{2C}{\lambda^n} \Vert h_n^2 \Vert^{\lambda^{2n-1}} +
\frac{2C}{\lambda^n} (\Vert h_n (z)\Vert^{\lambda^{2n-1}} + \Vert z \Vert^{\lambda^{2n-1}}) \; .
$$
By induction, if $l$ is such that all the iterates $h_n (z), h_n^2 (z), \ldots , h_n^{l-1} (z)$ remain in the disc of radius
$\tau$ for every $z$ with $\Vert z \Vert < \tau/2$, we derive the following estimate
\begin{eqnarray}
\Vert h_n^l (z) -R_n^l (z) \Vert & \leq & \frac{2C}{\lambda^n} \left( \Vert z \Vert^{\lambda^{2n-1}} + \Vert h_n (z)\Vert^{\lambda^{2n-1}}
+ \cdots + \Vert h_n^{l-1} (z)\Vert^{\lambda^{2n-1}} \right) \label{elementarybound4} \\
& \leq & \frac{2lC}{\lambda^n}  \Vert \tau \Vert^{\lambda^{2n-1}} \; . \label{elementarybound5}
\end{eqnarray}

For $\tau> 0$ small and fixed, we see that $l > k_n$. In fact, since $k_n < \beta \lambda^n/\alpha$, we obtain for $n$ sufficiently
large
$$
\Vert h_n^{k_n} (z) -R_n^{k_n} (z) \Vert \leq \frac{2\beta C}{\alpha} \Vert \tau \Vert^{\lambda^{2n-1}} \leq \frac{\tau}{2} \; .
$$
Furthermore $\Vert R_n^{k_n} (z) \Vert < \tau/2$ provided that $\Vert z \Vert < \tau/2$. This remark combines with the preceding
estimate to guarantee that $h_n^{k_n}$ is well-defined on $B_{\tau/2} (0)$ for $n$ large as above. Since the right
hand side in~(\ref{elementarybound5}) tends to zero as $n \rightarrow \infty$ ($\lambda > 1$), we can also conclude that
$\{ h_n^{k_n} \}$ converges uniformly towards $R_n^{k_n}$ on $B_{\tau/2} (0)$. This finishes the proof of the proposition.\qed

Let us denote by $\overline{\Gamma}$ the closure of $\Gamma$ (relative to $B_{\tau/2} (0)$). Naturally the contents of
Proposition~(\ref{prop7.A}) is that all the rotations $R_{\beta} = e^{2\pi i \beta} z$ belong to $\overline{\Gamma}$. With this
information in hand, let us go back to our original setting where $\Gamma$ is supposed to preserve a measure $\mu$
on $B_{\tau/2} (0)$ which, in addition, has no Dirac components. It is immediate to check that $\mu$ must also be preserved
by all elements lying in $\overline{\Gamma}$. In particular $\mu$ is preserved by the group of rotations $z \mapsto
e^{2\pi i \beta} z$.

Consider polar coordinates $r,\theta$ for $B_{\tau/2} (0)$. Since the only measures on the circle that are preserved by the
group of rotations are the constant multiples of the Haar measure, Fubini's theorem provides:

\begin{lema}
\label{lemma7.B}
The measure $\mu$ is given in polar coordinates by $T(r) drd\theta$ where $T$ is naturally identified with a
$1$-dimensional distribution.\qed
\end{lema}

Clearly all measures $\mu$ having the form indicated in the preceding lemma are automatically invariant by the group
of rotations. The fact that $\mu$ is also preserved by $f(z) =z^{\lambda}$ can then be translated into the functional equation
\begin{equation}
T (r) = \lambda^2 r^{\lambda -1} T(r^{\lambda}) \; . \label{elementarybound6}
\end{equation}

We are now able to prove Theorem~A.

\vspace{0.2cm}

\noindent {\it Proof of Theorem~A}. Recall that $l^+$ is contained in a closed set $\calk$ that is minimal
for $\fol$. Let us point out that the condition of having $\calk$ minimal has not been used so far. This condition however
is going to play a role in the sequel. To prove the theorem we are going to show that $\calk$ is itself an algebraic curve. To do this consider the
collection $\mathcal{L}$ of leaves of $\fol$ as defined in the beginning of this section. Clearly we have $l^+
\subset \mathcal{L}$. Next denote by $\overline{\mathcal{L}}$ the closure of $\mathcal{L}$ and consider the dimension
of the set $\overline{\mathcal{L}} \setminus \mathcal{L}$ of the proper accumulation points of $\mathcal{L}$. According to
the classical Remmert-Stein theorem, if the codimension of $\overline{\mathcal{L}} \setminus \mathcal{L}$ is at least two,
then $\overline{\mathcal{L}}$ is itself an analytic set so that the statement follows at once.

Thus we can suppose that the codimension of $\overline{\mathcal{L}} \setminus \mathcal{L}$ is strictly less than
two. In particular $\overline{\mathcal{L}} \setminus \mathcal{L}$ cannot be contained in the singular set of $\fol$.
Thus we can consider a point $p \in \overline{\mathcal{L}} \setminus \mathcal{L}$ that is regular for $\fol$. By considering
a plaque of $\fol$ containing $p$, we see that $\mathcal{L}$ must non-trivially accumulate on this plaque. Since
$\calk$ is minimal, it then follows that $\mathcal{L}$ has non-trivial recurrence. In other words, on $\Sigma$
(identified with the disc $B_{\tau} (0) \subset \C$), we can consider a point $q \in B_{\tau /2} (0)$, $q \neq 0$,
belonging to $\mathcal{L}$.

To establish the statement we shall derive a contradiction from the preceding with the fact that $q$
belongs to the support of the (transverse) invariant measure $\mu$. To do this, consider the pseudogroup
$\Gamma'$ of first return maps defined over paths in $\mathcal{L}$ but based at $q$. The group
$\Gamma'$ is conjugate to the group $\Gamma$. The desired contradiction arises as follows. Recall that
the structure of $\mu$ on $B_{\tau /2} (0)$ was already clarified by Lemma~\ref{lemma7.B} and Equation~(\ref{elementarybound6}).
Because $\Gamma'$ is conjugate to $\Gamma$, the analogous conclusions have to apply to a neighborhood of $q$
as well. In particular $\mu$ is ``constant'' over suitable closed loops about $q$. Since $\mu$ is also ``constant'' over the initial
circles about $0 \in B_{\tau/2} (0) \subset \C$, it follows that $\mu$ should be ``constant'' on a neighborhood of $q$, i.e.
on a neighborhood of $q$ the measure $\mu$ must be a constant multiple of the Lebesgue measure. This however contradicts
the analogous of Equation~(\ref{elementarybound6}) corresponding to the point $q$. The theorem is proved.\qed

Let us close this paper with the proof of Theorem~B. The method employed here is to large extent borrowed from
\cite{paul} to which we refer for further details. The prototype of an equation admitting a Liouvillean first integral
(integrable in the sense of Liouville) is the $1$-dimensional equation $y' = a(x) y + b(x)$ for which an explicit solution
involving two integrals can be obtained. In the complex domain, these integrals are in general multivalued so that, loosely
speaking, we can say that the equation admits a first integral that is ``twice multivalued''. Let us make this notion
precise.

Keeping the preceding notations, let us denote by $\calk$ the algebraic curve obtained from Theorem~A. In particular
$\calk$ coincides with $\overline{\mathcal{L}}$ where $\mathcal{L}$ was defined in the beginning of the section.
In the sequel consider meromorphic $1$-forms $\eta$ inducing $\fol$ but being defined only on  a neighborhood of $\calk$ in $M$.
So, unlike the previously used form $\omega$, $\eta$ need not be globally defined. Consider also a collection of local
representatives $\{ (U_a, \eta_a) \}$ for $\eta$ on a neighborhood of $\calk$. The compatibility condition among the local
representatives $\eta_a$ being given by the condition $\eta_a = u_{ab} \eta_b$ where $u_{ab} \in \mathcal{O}^{\ast} (U_a \cap
U_b)$. A {\it holomorphic integrating factor}\, for $\fol$ on the mentioned neighborhood consists of a collection
$\{ g_a\}$ of holomorphic functions, $g_a =u_{ab} g_b$, vanishing on $\calk$ and verifying
$$
d\left( \frac{\eta_a}{g_a} \right) =0 \; .
$$
The conditions above ensure that the local forms $\eta_a/g_a$ can be glued together to yield a closed meromorphic form
defining $\fol$ on a neighborhood of $\calk$. Therefore every primitive $H$ of the latter global form produces a multivalued
first integral for $\fol$, in fact, one has $\eta_a \wedge dH =0$ for all $a$.

A Liouvillean first integral for $\fol$ on a neighborhood of $\calk$ as above consists of going one step further into the
preceding discussion. A natural definition taken from \cite{paul} is as follows. Consider the universal covering
$\Pi: \mathcal{U} \rightarrow M \setminus \calk$ of $M \setminus \calk$. The sheaf $\mathcal{O}_{\mathcal{U}}$ induces a sheaf
over $M$ corresponding to its direct image by $\Pi$ and by the natural inclusion $M \setminus \calk \hookrightarrow M$. The restriction of the
latter sheaf to $\calk$ is going to be denoted by $\widetilde{\mathcal{O}}$. By construction an element belonging to the fiber of
$\widetilde{\mathcal{O}}$ over a point $q \in \calk$ is represented by a holomorphic function on $\Pi^{-1} (V \setminus V \cap \calk)$
where $V$ stands for a neighborhood of $q \in M$. The property of unique lift of functions through $\Pi$ allows us to identify
$\mathcal{O}_{\calk}$ to a subset of $\widetilde{\mathcal{O}}$ whose elements are, in addition, invariant by the local covering
automorphisms. With analogous constructions, we also define over $\calk$ the sheaves of (germs of) multivalued vector fields/holomorphic
forms.

Clearly the exterior differential $d$ can naturally be lifted to all above mentioned sheaves. An element $H$ of
$\widetilde{\mathcal{O}} (V)$ is said to be a {\it primitive}\, if $df$ is a $1$-form invariant by the local covering automorphisms and
admitting a meromorphic extension to $\calk$. Similarly $H \in \widetilde{\mathcal{O}}$ is said to be an {\it exponential of
primitive}\, if $df/f$ is a $1$-form invariant by the local covering automorphisms and
admitting a meromorphic extension to $\calk$. Let $\mathcal{S}^+ (V)$ (resp. $\mathcal{S}^{\times} (V)$) be the additive
(resp. multiplicative) subgroup of primitives (resp. exponential of primitive) of $\widetilde{\mathcal{O}} (V)$. The {\it first Liouvillean
extension}\, of $\mathcal{O} (V)$, denoted by $\mathcal{S}^1 (V)$ is the subring of $\widetilde{\mathcal{O}} (V)$ generated by
$\mathcal{S}^+ (V), \, \mathcal{S}^{\times} (V)$. The resulting presheaf turns out to be a sheaf over $\calk$. This construction can be
continued by induction to yield higher order Liouvillean extensions of $\mathcal{O} (V)$ but we shall not need those here (see \cite{paul}).

A Liouvillean (or $1$-Liouvillean) integrating factor for $\fol$ on  a neighborhood of $\calk$ consists of a collection $\{ g_a\}$
of elements in $\mathcal{S}^1 (U_a)$ such that
$$
g_a = u_{ab} g_b \; \; \; \, {\rm and} \; \; \; \, d \left( \frac{\eta_a}{g_a} \right) = 0 \; .
$$
A Liouvillean integrating factor is called {\it distinguished}\, if $g_a$ belongs to $\mathcal{S}^{\times} (V)$, ie. if $dg_a /g_a$
is a meromorphic closed form. By using this terminology we can state a slightly more accurate version of Theorem~B.

\begin{teo}
\label{strengthenedTheoremB}
Under the assumptions of Theorem~B the foliation$\fol$ admits a distinguished integrating factor. More precisely $\fol$
is given by a (Liouvillean) meromorphic closed form of  type $dg_a/g_a$ where $g_a \in  \mathcal{S}^{\times} (V)$.
\end{teo}

To prove Theorem~\ref{strengthenedTheoremB} consider the local transverse $\Sigma$ through $p \in l^+$
along with the pseudogroup $\Gamma$ of first return maps over paths contained in
$\mathcal{L}$ (recalling that $\overline{\mathcal{L}} = \calk$). Recall that $\Sigma$ is endowed with a coordinate
$z$ in which the first return $f$ over $l^+$ becomes $f(z) =z^{\lambda}$ on suitable sectors. The new ingredient
leading to the proof of the mentioned theorem is the following proposition.

\begin{prop}
\label{projectiveinvariance}
The vector field $\mathcal{X} = z \partial/\partial z$ is projectively invariant by $\Gamma$. In other words, if $h
\in \Gamma$ then
$$
h^{\ast} \mathcal{X} = c_h \mathcal{X}
$$
for a constant $c_h$ and whenever both sides are defined.
\end{prop}

To not interrupt the discussion we shall prove this proposition later. In the sequel we shall derive Theorem~\ref{strengthenedTheoremB}.
To make the argument more transparent, suppose first that $\mathcal{Y}$ were a vector field on $\Sigma$ fully invariant by $\Gamma$,
ie. satisfying $h^{\ast} \mathcal{Y} =  \mathcal{Y}$ for every $h \in \Gamma$. If this vector field exists, then it induces a vector field
(or rather a $1$-parameter subgroup of automorphisms) on the leaf space of $\fol$ (restricted to a neighborhood of $\calk$ as it will always be
the case in what follows). More precisely, on a neighborhood of $\calk$ consider the sheaf $\Theta_{M \fol}$ consisting of germs
of holomorphic vector fields tangent to $\calk$ and preserving the foliation $\fol$. If $\mathcal{Z}(V)$ is an element of $\Theta_{M \fol}$
then it verifies $L_{\mathcal{Z}(V)} \eta_a \wedge \eta_a=0$, where $L_{\mathcal{Z}(V)}$ denotes the Lie derivative, as a consequence
of the fact that $\mathcal{Z}(V)$ preserves $\fol$. Similarly we denote by ${\rm Tang}_{M \fol}$ the subsheaf of $\Theta_{M \fol}$
constituted by those germs of vector fields that are tangent to $\fol$. The sheaf ${\rm Symm}_{M \fol}$ of ``symmetries'' of $\fol$ is then defined by means
of the following exact sequence
$$
0 \longrightarrow {\rm Tang}_{M \fol} \longrightarrow \Theta_{M \fol} \longrightarrow {\rm Symm}_{M \fol} \longrightarrow 0 \, .
$$
With this terminology, it is clear that the vector field $\mathcal{Y}$ defined on $\Sigma$ and invariant by the holonomy of $\fol$
naturally induces a section of ${\rm Symm}_{M \fol}$ which is still denoted by $\mathcal{Y}$.

Next let a function $g_a$ be defined on each open set $U_a$ by the equation $g_a = \eta_a (\mathcal{Y})$. Because $\mathcal{Y}$
is identified with a (global) section of ${\rm Symm}_{M \fol}$, it follows that $g_a =u_{ab} g_b$ so that the collection $\{ g_a\}$
forms a holomorphic integrating factor for $\fol$. In fact, the condition $L_{\mathcal{Z}(V)} \eta_a \wedge \eta_a=0$ combines with the
Cartan formula $L_{\mathcal{Y}} =di_{\mathcal{Y}} + i_{\mathcal{Y}} d$ to yield
$$
d \left( \frac{\eta_a}{g_a} \right) = 0 \, .
$$
In other words, the collection of $1$-forms $\{ \eta_a/g_a\}$ defines a closed meromorphic form $\eta$ defined on a neighborhood of the {\it regular
part} of $\calk$ in $M$. However, the extension of this form to the singular points of $\fol$ lying in $\calk$ poses no difficulties
since all these singularities must have two non-vanishing real eigenvalues as a consequence of the discussion in Section~4 (cf. also
\cite{paul}). Finally a multivalued meromorphic first integral for $\fol$ can be obtained by means of the (multivalued) integral
$$
\int  \eta \; .
$$
In particular it follows that the ambiguity in the definition of the mentioned first integral is precisely determined by the
group of periods of $\eta$.

\vspace{0.2cm}

\noindent {\it Proof of Theorem~\ref{strengthenedTheoremB}}. Since a vector field $\mathcal{Y}$ on $\Sigma$
that is invariant by the pseudogroup $\Gamma$ need not exist, we shall generalize the preceding discussion to the
projectively invariant vector field $\mathcal{X}$ whose existence in ensured by Proposition~\ref{projectiveinvariance}.

Unlike the previous discussion, we are now going to exploit the fact that $\fol$ is given by a globally defined meromorphic
form $\omega$. Thus, if we consider the collection of local forms $\{ \eta_a\}$ obtained from the restrictions of $\omega$,
we conclude that the transition functions $u_{ab}$ are all constant equal to~$1$. Now on each open set $U_a$ we consider
a projectively invariant vector field $\mathcal{X}_a$. Thanks to Proposition~\ref{projectiveinvariance} this collection of vector fields
can be chosen so that $\mathcal{X}_a =c_{ab} \mathcal{X}_b$ where all the (transition) functions $c_{ab}$ are constant.
Again on the collection of open sets $\{ U_a\}$, we define the local functions $g_a = \eta (\mathcal{X}_a )$. Thus $g_a =c_{ab} g_b$.

Now on each $U_a$ consider the (local) meromorphic $1$-form $\Omega_{1a} = dg_a/g_a$ (with simple poles over $\calk$).
Since $g_a =c_{ab} g_b$ we conclude that $\Omega_{1a} = \Omega_{1b}$ so that these local forms glue together into a
meromorphic form $\Omega_1$ defined on a neighborhood of $\calk$. It is straightforward to check that the following
relations hold:
$$
\eta_a \wedge \Omega_{1a} = d\eta_a \; \; \; \, {\rm and} \; \; \; \, d\Omega_{1a} = 0 \, .
$$
In turn these relations are equivalent to the fact that the local functions $\{ g_a\}$ provide a distinguished Liouvillean integrating
factor for $\fol$ on  a neighborhood of $\calk$, since $u_{ab} =1$ and $\Omega_{1a} = \Omega_{1b}$ cf. \cite{paul}.
For the same reason mentioned above the extension of these
factors to the singularities of $\fol$ lying in $\calk$ poses no additional difficulty. The theorem is proved.\qed

The reader will note that the condition $d\omega = \omega \wedge \Omega_1$ means that the restriction of $\Omega_1$ to the
leaves of $\fol$ actually coincides with the foliated form $\omega_1$.

As before we can use the collection $\{ g_a\}$ to obtain a multivalued meromorphic first integral for $\fol$. For this we consider
the (multivalued) integrating factor $h_a = \int (\eta_a /g_a)$. Since the monodromy of the collection $\{ (U_a, g_a\}$ amounts
to multiplication by a constant ($g_a =c_{ab} g_b$), the same type of monodromy acts on the collection $\{ h_a\}$ by affine transformations. Using these
multivalued functions $\{ h_a \}$ we can define a (further multivalued) first integral analogous to the one found in the first case
discussed above. In slightly vague terms, the resulting first integral is given by
$$
\int \frac{\omega}{\int \Omega_1} \, .
$$
In particular the ambiguity in the definition of the first integral above lies in the period groups of $\omega$ and $\Omega_1$.

Let us now prove Proposition~\ref{projectiveinvariance}.

\vspace{0.2cm}

\noindent {\it Proof of Proposition~\ref{projectiveinvariance}}. Denote by $\phi_{\mathcal{X}}$ (resp. $\phi_{i\mathcal{X}}$)
the real (resp. purely imaginary) flow generated by $\mathcal{X}$. Note that $\phi_{i\mathcal{X}}$ is the real $1$-parameter
group consisting of rotations about $0 \in \C$. Proposition~\ref{prop7.A} then says that $\phi_{i\mathcal{X}}$ is contained in
$\overline{\Gamma}$. Next suppose for a contradiction that $h \in \Gamma$ does not preserve $\mathcal{X}$ up
to a multiplicative constant.

By construction $\Gamma$ admits a finite generating set whose elements either are defined on a neighborhood
of $0 \in \C$ or are defined on an appropriate sector $W$ (with vertex at $0 \in \C$). In the latter case, the element
in question is holomorphic on $W$ and of the form $z \mapsto z^a (1 + u(z))$ where $a \in \C$ and $u$ is defined on
a neighborhood of $0 \in \C$ with $u(0) =0$. In fact, the generators of $\Gamma$ that are not defined on a neighborhood
of $0 \in \C$ are obtained by means of singular loops passing through Siegel singularities of $\fol$ in $\calk$ what leads to
the general form mentioned above. Among these elements we have $f(z) = z^{\lambda}$ with $\lambda > 1$.

In view of the preceding we can suppose without loss of generality that either $h$ is defined on a neighborhood of
$0 \in \C$ or is of the form $h(z) = z^a (1 + u(z))$.

\vspace{0.1cm}

\noindent {\it Claim 1}. $h$ does not preserve the orbits of $\phi_{i\mathcal{X}}$.

\noindent {\it Proof of Claim 1}. Since $h$ does not preserve $\mathcal{X}$ up to a constant factor, it follows
that
$$
h^{\ast} \mathcal{X} = c_h (z) \mathcal{X}
$$
where $c_h (z)$ is a non-constant holomorphic function on its domain. By construction the purely imaginary
flow associated to $h^{\ast} \mathcal{X}$ is contained in $\overline{\Gamma}$ since so is $\phi_{i\mathcal{X}}$.
Denoting by $\mathcal{X}_1$ the (real) vector field associated to this (real) $1$-parameter group, we see that
$\mathcal{X}_1$ is $\R$-linearly independent with $\mathcal{X}$ at points in the (non-empty) open set where
$c_h$ takes values in $\C \setminus \R$. This proves the claim.\qed

Next we have:

\vspace{0.1cm}

\noindent {\it Claim 2}. The vector fields $\mathcal{X}$ and $\mathcal{X}_1$ do not commute. Besides
$[\mathcal{X}, \mathcal{X}_1]$ is not constant.

\noindent {\it Proof of Claim 2}. Recall that $\mathcal{X} = z \partial /\partial z$ whereas $\mathcal{X}_1 =
c_h (z)z \partial /\partial z$. Hence
$$
[\mathcal{X}, \mathcal{X}_1] = -z^2 c_h'(z) \, .
$$
Since we have assumed that $c_h$ is not constant we conclude that $c_h'(z)$ is not identically zero and
the claim follows.\qed

Consider a (``generic'') point $p$ at which $\mathcal{X}, \; \mathcal{X}_1$
are $\R$-linearly independent and where $[\mathcal{X}, \mathcal{X}_1] (p) =(a+ib) \partial /\partial z$ with
$a \neq 0$. The existence of these points follows from Claim~2.
For fixed small reals $s,t$, consider the local diffeomorphism $D^{st}$ fixing $p$ that is obtained as follows:
points close to $p$ are moved by following the (purely imaginary) flow of $\mathcal{X}$ during a time~$s$ and
then we compose this with the purely imaginary flow of $\mathcal{X}_1$ during a time~$t$. Finally we still compose
the mentioned map with the (always purely imaginary) flow of $\mathcal{X}$ (and then $\mathcal{X}_1$) during
a time $s'$ close to $s$ (resp. $t'$ close to $t$) so as to have $p$ as a fixed point of the resulting map $D^{st}$.
Clearly the map $D^{st}$ belongs to $\overline{\Gamma}$ for small $s,t$. It must therefore preserve the transverse measure
over $\Sigma$ identified to a neighborhood of $0\in \C$. Therefore to derive a contradiction proving the statement, it
suffices to show that $p$ is a hyperbolic fixed point provided that $s,t$ are small enough. To check this simply note
that the derivative of $D^{st}$ at $p$ is given by
$$
1 -st [\mathcal{X}, \mathcal{X}_1] (p) + o\,  (s^2 +t^2) = 1 -st (a+ib) + o\,  (s^2 +t^2)\, .
$$
Since $\Vert 1 -st (a+ib) \Vert = 1 -2sta + s^2 t^2 (a^2 +b^2)$ with $a\neq 0$, it follows that the norm of the derivative
of $D^{st}$ at $p$ is different from~$1$ provided that $s,t$ are very small. This finishes the proof of the proposition.\qed

\vspace{0.1cm}

\centerline{\sc {\large Examples}}

Let us close this paper with two classes of non-trivial examples of foliations for which Theorems~A and~B in the Introduction
can immediately be applied. In the sequel consider a foliation $\fol$ on a surface $M$ along with a global meromorphic form $\omega$ defining
$\fol$. Unless otherwise stated we choose $\omega$ so that its divisor of zeros and poles $(\omega)_{\infty} \cup (\omega)_0$
is disjoint from the finite set formed by the singularities of $\fol$. To further simplify the discussion, suppose also that
all singularities of $\fol$ have at least one eigenvalue different from zero. In view of Seidenberg's theorem our construction can
easily be adapted to include foliations with degenerate singularities so that we shall not worry about them here.

Obviously the simplest way to ensure the existence of $\calh$-trajectories with infinite length consists of eliminating the singular points
of $\calh$ that give rise to either future endpoints or to past endpoints for trajectories of $\calh$. Under the above assumptions we have:
\begin{itemize}

\item All points giving rise to past-ends for trajectories of $\calh$ are contained in $(\omega)_{\infty}$.

\item All points giving rise to future-ends for trajectories of $\calh$ are contained in the union of $(\omega)_0$
with the singular set ${\rm Sing}\, (\fol)$ of $\fol$.
\end{itemize}
In addition for a singular point of $\fol$ to yield future ends for $\calh$-trajectories the ratio between its eigenvalues
must belong to $\R_+^{\ast}$.

To apply our theorems we need to have $(\omega)_{\infty} \neq \emptyset$ for otherwise $\omega$ is holomorphic and
hence closed. Thus this particular choice of $\omega$ does not yield an associated foliation $\calh$. Hence a natural idea
is to eliminate the possibility future-ends for the trajectories of $\calh$. With this idea in mind we shall provide our first class
of examples.

\noindent {\sc Example 1}. Foliations with $(\omega)_0 = \emptyset$ and singularities whose eigenvalues have quotients
in $\C \setminus \R_+^{\ast}$.

We note that the class of foliations above include those with singularities in the Siegel domain. Saddle-node singularities
can also be authorized for $\fol$. Suppose that $T$ is a diffuse closed
current invariant by $\fol$. The first important remark to be made about $T$ concerns its co-homology class in $M$. In fact,
under the conditions regarding the singularities of $\fol$ the following holds:
$$
[T].[T] =0
$$
ie. the self-intersection of $T$ vanishes. The proof of the above equation is easy and essentially amounts to the fact that
Siegel singularities cannot contribute non-trivially for the self-intersection of a diffuse current, see for example \cite{marco}. Naturally, as seen in Section~4,
all singularities of $\fol$ lying in the support of $T$ must belong to the Siegel domain since the support of a diffuse current as above
cannot contain either hyperbolic singularities or saddle-nodes, cf. Section~4.
This case is therefore of little interest  for surfaces, such as $\C P (2)$,
whose Picard group is cyclic. However
for surfaces with larger Picard group the condition about the self-intersection of $T$ conveys less information.

For example let $M$ be an affine elliptic $K3$ surface in $\C^3 \subset \C P(3)$ (i.e. the closure of $M$, still denoted by $M$, in
$\C P(3)$ is a $K3$-surface). Indeed, we can choose $M$ to be the usual Fermat's quartic. Next consider a polynomial vector field $X$ tangent to
$M$ and having isolated singularities. The vector field $X$ induces a foliation $\fol$ over $M \subset \C P(3)$ whose singularities satisfy the above conditions modulo
choosing $X$ ``generic''. By exploiting the triviality of the canonical bundle of
$M$, we can easily find a $1$-form $\omega$ on $M$ defining $\fol$ and having empty divisor of zeros. Furthermore the divisor of poles
of $\omega$ coincides with the hyperplane section of $M$ and, in general, $\omega$ is not closed. Now if $T$ is a closed current
invariant by $\fol$, the condition that $[T].[T] =0$
says that $[T]$ is the cohomology class of an elliptic fiber of $M$. However there is {\it a priori}\, no reason to conclude the existence of any compact
leaf for $\fol$ unless Theorem~A is used.

\noindent {\sc Example 2}. Foliations on $\C P(2)$ with singularities whose quotient of eigenvalues belong to $\R_+$.

A disadvantage of the construction employed in Example~1 is that, after all, it depends on the fact that only Siegel singularities
can appear on the support of a diffuse (closed positive) invariant current. Similarly it has been know since long that a foliation of
$\C P(2)$ all of whose singularities are hyperbolic cannot admit a diffuse current as above. To be able to make a significant progress with respect
to these well-known results, it is interesting to allow the foliation to have simple singularities of type ``irrational focus''. In fact, these singularities
may contribute non-trivially to self-intersection of $T$ and they do not contradict the diffuse nature of $T$. The example below is intended to
show that our theorem can be used to provide new results in this context.

To simplify suppose that this happens at a single singularity $p$ with eigenvalues $1,1$ (ie. in local coordinates about $p$ $\fol$ is given by the radial
vector field $x\partial /\partial x + y \partial /\partial y$). The remaining singularities of $\fol$ being hyperbolic, of Siegel type or saddle-nodes.
The existence of $p$ prevents us from concluding that $[T].[T] =0$. Nonetheless we claim the following:

\noindent {\it Claim}. If $\fol$ admits a closed current $T$, then it must also possess an algebraic curve provided that
the degree of $\fol$ is at least~$3$.

\noindent {\it Proof}. Choose affine coordinates so that $p$ belongs to the corresponding line at infinity
$\Delta$. In the affine $\C^2$ we choose a polynomial form $\omega$ representing $\fol$ and such that
its components have only trivial common factors. Viewed as a meromorphic $1$-form on $\C P(2)$, the divisor of zeros of
$\omega$ is empty whereas $\omega$ has poles of order $d-1$ over $\Delta$ where $d$ stands for the degree of
$\fol$. Thanks to Theorem~A to prove the claim it suffices to show that all the trajectories of the resulting foliation
$\calh$ have infinite length. In turn, for this, it is enough to show that $p$ does not provide future endpoints for
these trajectories. To check the latter claim, consider local coordinates $(x,y)$ about $p$, $\{ x=0\} \subset
\Delta$, where $\omega = v(x,y) x^{1-d} (xdy - y dx)$ with $v(0,0) \neq 0$. If we blow-up $\fol$ at $p$, the new
foliation $\widetilde{\fol}$ has no longer singularities over the exceptional divisor which, in fact, is transverse
to $\widetilde{\fol}$. In standard $(x,t)$ coordinates for the blow-up the pull-back of $\omega$ is given
by
$$
v(0,0) x^{3-d} dt \, .
$$
In particular the pull-back of $\omega$ does not have zeros over the exceptional divisor since $v(0,0) \neq 0$
and $d \geq 3$. Therefore the initial singularity $p$ does not provide future endpoints for the trajectories
of $\calh$. The claim is proved.\qed

Naturally the eigenvalues of $\fol$ at $p$ may be supposed to have only the form $1, \lambda$, where $\lambda \in
\R_+^{\ast}$. Besides, if the degree of $\fol$ is at least~$4$, we can allow the existence of two (rather than one) irrational
focus singularities for $\fol$. Several other combinations of these ideas can be used to provide new results on the
structure of invariant curves for foliations as above.

\begin{flushleft}

{\sc Julio C. Rebelo}\\
rebelo@math.univ-toulouse.fr\\

\vspace{0.1cm}

1 - Universit\'e de Toulouse , UPS, INSA, UT1, UTM\\
Institut de Math\'ematique de Toulouse\\
F-31062 Toulouse, FRANCE\\

\vspace{0.1cm}

2 - CNRS, Institut de Math\'ematique de Toulouse UMR 5219\\
F-31062  Toulouse, FRANCE

\end{flushleft}

\end{document}